\input amstex
\documentstyle{amsppt}
\mag=\magstep1
\pageheight{20.5cm}
\vcorrection{-0.2cm}
\pagewidth{15cm}
\hcorrection{-1.0cm}
\document
\topmatter
\title Ricci Flow Unstable Cell Centered at an Einstein 
Metric on the Twistor Space of Positive Quaternion K\"ahler Manifolds 
of Dimension $\geq 8$ \endtitle
\rightheadtext{Ricci Flow Unstable Cell}
\author Ryoichi Kobayashi \endauthor
\affil Graduate School of Mathematics, Nagoya University \endaffil
\thanks{Th first author is thankful to professor Bogdan Alexandrov for 
pointing out an error in the earlier version of this paper.}\endthanks
\abstract{We construct a 2-parameter family $\Cal F^{\text{\rm Z}}$ 
of Riemannian metrics on the twistor space $\Cal Z$ 
of a positive quaternion K\"ahler manifold $M$ satisfying 
the following properties : (1) the family $\Cal F^{\text{\rm Z}}$ 
contains an Einstein metric $g^{\text{\rm Z}}$ and its scalings, 
(2) the family $\Cal F^{\text{\rm Z}}$ is closed under the 
operation of making the convex sums, 
(3) the Ricci map $g \mapsto \text{\rm Ric}(g)$ defines a dynamical system 
on the family $\Cal F^{\text{\rm Z}}$, 
(4) the Ricci flow starting at any metric in the family $\Cal F^{\text{\rm Z}}$ 
stays in $\Cal F^{\text{\rm Z}}$ and is an ancient solution having the Einstein 
metric.$g^{\text{\rm Z}}$ as its asymptotic soliton. 

This means that the family $\Cal F^{\text{\rm Z}}$ is a 2-dimensional ``unstable cell'' 
w.r.to the Ricci flow which is ``centered'' at the Einstein metric $g^{\text{\rm Z}}$. 
We apply the estimates for the covariant derivative of the curvature tensor 
under the Ricci flow to this ``unstable cell'' and settle the LeBrun-Salamon 
conjecture : any irreducible positive quaternion 
K\"ahler manifold is isometric to one of the Wolf spaces.}
\endabstract
\endtopmatter

\def\Z{\Bbb Z} \def\H{\Bbb H} \def\C{\Bbb C} \def\R{\Bbb R} \def\V{\Cal V} 
\def\Sp{\text{\rm Sp}} \def\SO{\text{\rm SO}} \def\SU{\text{\rm SU}} \def\Un{\text{\rm U}}
\def\sp{\text{\rm sp}}  \def\so{\text{\rm so}} 
\def\Om{\Omega} \def\PP{\Bbb P} \def\FS{\text{\rm FS}} 
\def\ov{\overline} \def\wt{\widetilde}  \def\ep{\varepsilon}
\def\Ric{\text{\rm Ric}} \def\Rm{\text{\rm Rm}} 
\def\tr{\text{\rm tr}} \def\p{\partial} \def\Vol{\text{\rm Vol}} \def\Diff{\text{\rm Diff}} 
\def\euc{\text{\rm euc}} \def\loc{\text{\rm loc}} \def\can{\text{\rm can}} 
 \def\G{\text{\rm G}}  \def\wcalD{\wt{\Cal D}} 
\def\ZZ{\text{\rm Z}} \def\wcalZ{\widetilde{\Cal Z}} \def\wcalF{\widetilde{\Cal F}}

\NoBlackBoxes

\noindent
{\bf \S1. Introduction.}
\medskip

In this paper, we construct a new Einstein metric together with a two parameter family of metrics 
containing it on the twistor space of a positive quaternion K\"ahler manifold and describe the behavior 
of this family under the Ricci flow. 

Because the background of this attempt is rather heavy, we start with the explanation 
what ``Ricci flow unstable cell'' means. 
In his beautiful paper [P], G. Perelman introduced his $\Cal W$-functional 
$$\Cal W(g_{ij},f,\tau)=\int_M[\tau(R+|\nabla f|^2)+f-n](4\pi\tau)^{-\frac n2}e^{-f}dV$$
where $M$ is an $n$-dimensional smooth (closed) manifold, $g_{ij}$ is a Riemannian 
metric on $M$, $f$ is a smooth function on $M$, $dV$ is the Riemannian volume form 
of the metric $g_{ij}$ and $\tau$ is a positive scaling parameter. 
The remarkable property of the $\Cal W$-functional proved in [P] is the following. 
Consider the system
$$\left\{\aligned
& \p_tg_{ij}(t,x)=-2R_{ij}(g(t,x))\\
& \p_t u(t,x)=-\triangle_x u(t,x)+Ru(t,x)\qquad (u:=(4\pi\tau)^{-\frac n2}e^{-f})\\ 
& \p_t\tau(t)=-1\endaligned\right.$$
of evolution equations on $M$, where $\triangle$ is the Laplacian w.r.to the metric 
$g_{ij}(t,x)$ acting on functions and $R$ is the scalar curvature of the metric $g_{ij}(t,x)$. 
Note that the second equation is the conjugate heat equation 
and therefore $\int_MudV$ is independent of $t$ (usually normalized to be $1$). 
The short time existence of this system of evolution equations is a consequence of 
the short time existence of the Ricci flow equation (see [P,\S1]). 
The $\Cal W$-functional is monotone non-decreasing along its solution. In fact we have 
the monotonicity formula
$$\p_t\Cal W=2\int_M \tau|R_{ij}+\nabla_i\nabla_jf-\frac1{2\tau}g_{ij}|^2udV\,\,.$$
The origin of this monotonicity formula is the folllowing. 
If we apply the family of $t$-dependent diffeomorphisms $\{\phi_t\}$ obtained by diagonally 
integrating the $t$-dependent vector field $-\nabla f$ ($\nabla$ being the Levi-Civita connection 
w.r.to the metric $g_{ij}(t,x)$) to the above system of equations, we get the system
$$\left\{\aligned
& \p_tg_{ij}(t,x)=-2(R_{ij}(t,x)+\nabla_i\nabla_jf(t,x))\\
& \p_tf(t,x)=-\triangle_x f(t,x)-R(t,x)+\frac{n}{2\tau(t)}\\
& \p_t \tau(t)=-1\endaligned\right.$$
which is precisely the $L^2$-gradient flow of the $\Cal W$-functional under the constraint 
that $(4\pi\tau)^{-\frac n2}e^{-f}dV$ is a fixed measure $dm$ on $M$ independent of $t$. 
In this situation, the $W$-functional should be interpreted as
$$\Cal W^m(g_{ij},f,\tau)=\int_M[\tau(R+|\nabla f|^2)+f-n]dm\,\,.$$
In particular the symmetry under which the functional $\Cal W^m$ in invariant is the 
subgroup of $\Diff_0(M)$ consisting of $dm$-preserving diffeomorphisms. 
The action of other diffeomorphisms is introduced as follows. 
For $\phi\in\Diff_0(M)$ we define $f^{\phi}$ 
by setting $dm=(4\pi\tau)^{-\frac n2}e^{f^{\phi}}dV_{\phi^*g}$. If $\phi$ preserves $dm$ then 
$f^{\phi}$ is just $\phi^*f$ and the converse is also true. The action of $\Diff(M)_0$ is 
thus defined on the configuration space $\{(g_{ij},f,\tau)\}$ on which the functional 
$\Cal W^m$ is defined. 

Perelman's $\Cal W$-functional is a ``coupling" of the {\bf logarithmic Sobolev functional}
\footnote{\,\,The logarithmic Sobolev inequality on the $n$-dimensional Euclidean space 
$\R^n$ is the following. Let $f=f(x)$ satisfies the constraint 
$\int_{\R^n}(4\pi\tau)^{-\frac2n}e^{-f}dV_{\euc}=1$. Then we have
$$\int_{\R^n}[\tau|\nabla f|^2+f-n](4\pi\tau)^{-\frac2n}e^{-f}dV_{\euc}\geq 0$$
where the equality holds iff $f(x)=\frac{|x|^2}{4\tau}$. } 
and the {\bf Hilbert-Einstein functional}\footnote{\,\,The Hilbert-Einstein functional is 
$\displaystyle \int_MRdV_g$ for a closed Riemannian manifold $(M,g)$ and the 
critical points are Einstein metrics.}. 
Suppose that there exists a critical point of the $\Cal W^m$-functional, 
which corresponds to a Ricci soliton 
$$R_{ij}+\nabla_i\nabla_jf-\frac1{2\tau}g_{ij}=0\,\,.$$
This is interpreted, at time $t=-1$ ($\tau=1$), as the initial condition for the Ricci flow equation 
(the solution satisfies the above equation of the Ricci soliton, which evolves under 
a 1-parameter group of diffeomorphisms of $M$ generated by the ``gradient" vector field 
$\nabla f$). 
Perelman [P] showed that the Ricci soliton is characterized by the equality case 
of the logarithmic Sobolev inequality in the following way. 
Let $g_{ij}(-1)$ satisfy the above equation at time $t=-1$ and $g_{ij}(t)$ 
the corresponding solution of the Ricci flow, i.e., the Ricci soliton with initial metric $g_{ij}(-1)$. 
Then the logarithmic Sobolev inequality on $(M,g_{ij}(t))$ introduced in [P] is
$$\split
& \quad \Cal W(g_{ij}(t),\widetilde{f},-t) \geq \Cal W(g_{ij}(t),f(t),-t)=
\inf_{\widetilde{f}:\int_M(4\pi(-t))^{-\frac n2}e^{-\widetilde{f}}dV_{g(t)}=1}\Cal W(g_{ij}(t),\widetilde f,-t)\\
& =:\mu(g_{ij}(t),-t)=\mu(g_{ij}(-1),1)\endsplit$$
where $\widetilde{f}$ is any smooth function on $M$ satisfying 
$\int_M(4\pi(-t))^{\frac n2}e^{-\widetilde{f}}dV_{g(t)}=1$. 
This observation gives us an important information on the behavior of the $\Cal W$-functional 
at a critical point (i.e., the Ricci soliton). We look at the Hessian of the $\Cal W^m$-functional 
at the critical point. The $\Cal W^m$-functional is invariant under the group of all $dm$-preserving 
diffeomorphisms and therefore this action corresponds to the zeros of the Hessian. 
On the other hand, the action of the diffeomorphisms which do not preserve $dm$ may be 
given by the following way. Let $\phi$ be such a diffeomorphism. 
Introduce $f^{\phi}$ by setting  $dm=(4\pi\tau)^{-\frac n2}e^{-f^{\phi}}dV_{\phi^*g}$ and define 
$\phi^*(g,f,\tau)=(\phi^*g,f^{\phi},\tau)$. Then we have
$$\Cal W^m(\phi^*(g,f,\tau))=\int_M[\tau(R_{\phi^*g}+|\nabla f^{\phi}|_{\phi^*g}^2)+f^{\phi}-n]
\underbrace{(4\pi\tau)^{-\frac n2}e^{-f^{\phi}}dV_{\phi^*g}}_{dm}$$
and therefore the $\Cal W^m$-functional increases in the direction of the action of 
the diffeomorphisms which do not preserve $dm$, which follows from the logarithmic Sobolev 
characterization of the Ricci soliton. 
This implies that the tangent space of the configuration space $\{(g,f,\tau)\}$ of the functional 
$\Cal W^m$ decomposes 
into three subspaces $V_0$, $V_+$ and $V_-$. Here, $V_0$ corresponds to the action of the 
$dm$-preserving doffeomorphisms (Hess $=0$), $V_+$ corresponds to the action of the 
diffeomorphisms which do not preserve $dm$ (Hess $\geq 0$) and finally $V_-$ 
(Hess $\leq 0$) corresponds to the rest\footnote{\,\, This is very similar to the behavior 
of the Hilbert-Einstein functional under the Yamabe problem.}. 

The above discussion means that we can speak of the stable/unstable cell w.r.to 
the Ricci flow in the space of Riemannian metrics on a given manifold 
(with respect to the $L^2$-gradient flow of the functional $\Cal W^m$). 

Perelman announced the following convergence theorem of the K\"ahler Ricci flow 
$\p_t g_{i\ov j}(t,x)=-R_{i\ov j}(g(t,x))+g_{i\ov j}(t,x)$ with the initial K\"ahler metric $g_{i\ov j}(0,x)$ 
representing $c_1(M)$ where $M$ is a Fano manifold $M$ (see Tian-Zhu [T-Z]). 

\proclaim{Theorem 1.1 (Perelman, Tian-Zhu)} If $M$ admits a K\"ahler-Einstein metric, then, 
for any initial K\"ahler metric 
$g_{i\ov j}(0,x)$ representing $c_1(M)$, the solution to the K\"ahler-Ricci flow converges to 
a K\"ahler-Einstein metric in the sense of Cheeger-Gromov. \endproclaim

This implies that the normalized K\"ahler-Ricci flow starting at a K\"ahler metric in $c_1(M)$ 
on a Fano manifold admitting a K\"ahler-Einstein metric forms a stable cell w.r.to 
the dynamical system defined by the Ricci flow solution (in short, we call this a 
{\it Ricci flow stable cell}). 

It is then interesting to construct an example of a {\it Ricci flow unstable cell} centered 
at an Einstein metric on an Einstein manifold. Theorem 1.1 implies that, on a positive 
K\"ahler-Einstein manifold with $b_2=1$, one can find a Ricci flow 
unstable cell centered at a K\"ahler-Einstein metric (if any) only from 
a solutions of the non K\"ahler Ricci flow, or one should find 
another non-K\"ahler Einstein metric together with a Ricci flow unstable cell. 

This paper is an attempt toward this direction. 
A Ricci flow unstable cell centered at an Einstein metric (if any) 
necessarily consists of {\it ancient solutions} (a Ricci flow solution is said to 
be an ancient solution [H] if it is defined on a time interval $(-\infty,T]$ where $T\in\R$) 
whose asymptotic soliton (in the sense of [P, \S11]) is the Einstein metric. 
Because there is no guarantee for evolution equations (such as the Ricci flow equation) 
to be solved in the past direction even locally, 
the existence of the ancient solution must be a miracle which comes from 
a very special geometric situation. 
It is best explained by Perelman's result : 
ancient solutions appear from the rescaling procedure of the singularities developed 
in the Ricci flow in finite time (this is a consequence of Perelman's Local Non-Collapsing 
Theorem [P, \S4 and \S7]). 
In this paper we pick up the {\it non holomorphic twistor fibration of the twistor space} $\Cal Z$ 
{\it of positive quaternion K\"ahler manifolds} $M$ as a candidate of a special geometric 
situation admitting a Ricci flow unstable cell centered at an Einstein metric. 
The strategy is to think of the non holomorphic 
$\PP^1$-fibration $\pi:\Cal Z\to M$, imagine the collapse of $\Cal Z$ to $M$ 
(i.e., the $\PP^1$-fiber collapses) or to $\PP^1$ (i.e., the base $M$ collapses) as a ``singularity" 
developed in the Ricci flow on $\Cal Z$ in finite time and then try to construct the corresponding 
``ancient solutions" by suitable rescaling procedure. 
However, this attempt do not seem to work well if one stick to the collapse along 
the twistor fibration. In \S2, we study the family $\Cal F^{\can}$ of the canonical deformation 
metrics and show that $\Cal F^{\can}$ is ``stable'' under the Ricci flow. It turns out 
that the Ricci flow solution converges to the K\"ahler-Einstein metric if we start the Ricci flow 
at any canonical deformation metric sufficiently close to the K\"ahler-Einstein metric. 
The purpose of \S2 is to explain why the canonical deformation metrics do not form 
the Ricci flow unstable cell centered at the K\"ahler-Einstein metric. 
In \S3, we introduce an idea to ``kill'' the basic reason why the family of the canonical deformation 
metric is not a Ricci flow unstable cell and modify the construction of the canonical deformation metric. 
As a result, we construct a new Einstein metric $g_{\sqrt{\frac1{n+2}}}^{\ZZ}$ and a new 
family of $\Cal F^{\ZZ}$ metrics (containing $g^{\ZZ}_{\sqrt{\frac1{n+2}}}$) on the twistor space 
and show that the family $\Cal F^{\ZZ}$ constitutes a Ricci flow unstable cell centered at 
the Einstein metric $g^{\ZZ}_{\sqrt{\frac1{n+2}}}$. 

An oriented $(4n)$-dimensional Riemannian manifold $(M^{4n},g)$ is said to be a {\it quaternion 
K\"ahler manifold}, if its holonomy group is contained in $\Sp(n)\Sp(1)$ $(\subset \SO(4n))$ 
(for a precise definition, see \S2). 
A locally irreducible (in the sense of the local de Rham decomposition defined from the 
irreducible decomposition of the tangent space under the action of the local holonomy group) 
quaternion K\"ahler manifold is necessarily Einstein. 
Therefore we can classify locally irreducible quaternion K\"ahler manifolds 
by the sign of the scalar curvature into three classes. 
A complete locally irreducible quaternion K\"ahler manifold is called positive 
(resp. locally hyper-K\"ahler, negative), if its scalar curvature is positive (resp. zero, negative). 
The locally hyper-K\"ahler  property is equivalent to the absence of the $\Sp(1)$-component 
in the $\Sp(n)\Sp(1)$-holonomy. 
An irreducible positive quaternion K\"ahler manifold turns out to be a simply connected positive 
Einstein manifold. From here on, we restrict our attention to locally irreducible quaternion 
K\"ahler manifolds (and therefore we just say ``quaternion K\"ahler" for ``locally 
irreducible quaternion K\"ahler"). 

The principal bundle of oriented orthonormal frames of a quaternion K\"ahler manifold 
is reduced to a principal $\Sp(n)\Sp(1)$-bundle $\Cal P$ (the holonomy reduction). 
Associated to a quaternion K\"ahler manifold is the twistor space $\Cal Z$, 
which is constructed from $\Cal P$ and the $\Sp(1)$-part of the holonomy group $\Sp(n)\Sp(1)$ 
by putting 
$$Z:=\Cal P\times_{\Sp(n)\Sp(1)}\Bbb P^1=\Cal P/\Sp(n)\Sp(1)\cap U(2n)$$
(this is a $\Bbb P^1$-bundle over $M$). 
The horizontal distribution w.r.to the Levi-Civita connection on the twistor fibration 
$\pi:Z \rightarrow M$ is canonically defined and thus the twistor space $\Cal Z$ 
is equipped with a natural almost complex structure (which turns out to be integrable) 
and two kinds of family of metrics which is ``compatible" with the twistor space construction. 
The first one is the family of the {\it canonical deformation metrics} defined by adding scaled 
standard metrics of the $\Bbb P^1$-fiber and the fixed base metric on $M$ by using the horizontal 
distribution (we will use the description by Chow and Yang in [C-Y]). 
The second one (called the family of {\it Z-metrics}) will be introduced in this paper (see \S 3), 
whose construction is based on the canonically defined a horizontal $(4n-2)$-dimensional 
distribution $\Cal D'$  and the ``associated'' non-horizontal $(4n)$-dimensional distribution 
$\Cal D$ on $\Cal Z$. The goal of this paper is to describe the behavior under the Ricci flow 
of the family of Z-metrics (and comparison with that of the family of canonical deformation 
metrics). 

Typical examples of positive quaternion K\"ahler manifolds are the ``Wolf spaces" [W], i.e., 
$(4n)$-dimensional positive quaternion K\"ahler symmetric spaces. Wolf [W] proved that 
for each compact simple Lie group $G$, there is a Wolf space $G/H$. The classical ones are
$$\split
&\H\PP^n=\frac{\Sp(n+1)}{\Sp(n)\times\Sp(1)},\,\,
\text{\rm Gr}_2(\C^n)=\frac{\SU(n)}{\text{\rm S}(\text{\rm U}(n-2)\times\text{\rm U}(2))},\,\,\\
& \wt{\text{\rm Gr}}_4(\R^{n})=\frac{\SO(n)}{\SO(n-4)\times SO(4)}\endsplit$$
and there are exceptional cases. The noncompact dual of the Wolf spaces are examples 
of negative quaternion K\"ahler manifolds. There exist many other examples of noncompact 
negative quaternion K\"ahler manifolds which are not symmetric (e.g., Alexeevskii [A], Galicki [G]). 
Moreover, we remark that Galicki-Lawson's quaternion K\"ahler reduction method ([G-L]) 
produces many examples of positive quaternion 
K\"ahler orbifolds which are not symmetric\footnote{\,\,the arguments in \S2 for the twistor 
space of positive quaternion K\"ahler manifolds cannot be generalized to the orbifold case 
(see Remark 3.5 in \S3).}. 

\proclaim{Conjecture 1.2 ([L-S])} Any irreducible positive quaternion K\"ahler manifold 
is isometric to one of the Wolf spaces. 
\endproclaim

The twistor space of a Wolf space is a generalized flag manifold 
(rational homogeneous manifold) of the form $G/P$, which is a Fano manifold. 
More precisely, the twistor space of a positive quaternion K\"ahler manifold is a 
``contact Fano manifold" (see (2-13)). LeBrun [LeB] proved that a contact Fano manifold $\Cal Z$ 
is realized as the twistor space of some positive quaternion K\"ahler manifold $M$ 
if and only if $\Cal Z$ admits a K\"ahler-Einstein metric. 

\proclaim{Conjecture 1.3 ([L-S])} Any contact Fano manifold admits a K\"ahler-Einstein metric.
\endproclaim

The affirmative answers to Conjectures 1.2 and 1.3 combined with the result of [LeB] imply 
that a compact complex manifold $\Cal Z$ is the twistor space of a Wolf space 
if and only if $\Cal Z$ is a contact Fano manifold. 
On the other hand, using complex algebraic geometry of contact Fano manifold, 
LeBrun and Salamon [L-S] proved that there exist only finitely many poisitive quaternion 
K\"ahler manifolds with fixed dimension. 
\medskip

We give a brief description on the basic properties of the family of Z-metrics and their 
construction. 

The family $\Cal F^{\ZZ}$ of Z-metrics satisfies the following properties : 
(1) the family $\Cal F^{\ZZ}$ contains an Einstein metric $g^{\ZZ}$ which is 
different from the well-known two Einstein metrics in the family of 
canonical deformation metrics, 
(2) the family $\Cal F^{\ZZ}$ is closed under 
positive scalar multiples and summation, 
(3) the family $\Cal F^{\ZZ}$ is stable under 
the Ricci map $g\mapsto \Ric(g)$ (i.e., the Ricci map sends $\Cal F^{\ZZ}$ to itself), 
(4) any Ricci flow solution starting with initial metric in the family $\Cal F^{\ZZ}$ 
stays in $\Cal F^{\ZZ}$ and is an ancient solution whose asymptotic soliton is the 
Einstein metric $g^{\ZZ}$. 

The construction of the family $\Cal F^{\ZZ}$ is based on the 
canonically defined real $(4n)$-dimensional non-horizontal 
distribution $\Cal D=\{\Cal D_z\}_{z\in\Cal Z}$. The distribution $\Cal D'$ 
is constructed as follows (see \S3). 
There canonically exists a horizontal $(4n-2)$-dimensional distribution 
$\Cal D'=\{\Cal D'_z\}_{z\in\Cal Z}$. Here, $\Cal D_z'$ is $J(z)$-invariant 
where $J(z)$ is the orthogonal complex structure corresponding to $z$. 
Let $L_{J(z)}=(\Cal D_z')^{\perp_{\Cal H_z}}$ where $\Cal H=\{\Cal H_z\}_{\in\Cal Z}$ 
is the horizontal distribution. 
Let $L_{J(z)}$ be the ``diagonal'' in the complex 2-dimensional 
distribution $\{\V_z\oplus L_{J(z)}'\}_{z\in\Cal Z}$ where $\V_z$ is the 
vertical distribution. The $(4n)$-dimensional non-horizontal 
distribution $\Cal D$ on the twistor space $\Cal Z$ is defined by 
$\Cal D_z=\Cal D'_z\oplus L_{J(z)}'$, which is still transversal to 
the twistor fibration $\Cal Z\rightarrow M$. Set $\Cal D^{\perp}
=\{\Cal D_z^{\perp}\}_{z\in\Cal  Z}$. This is a 2-dimensional distribution 
which is invariant w.r.to the canonical complex structure of the twistor space $\Cal Z$. 
Here, the orthogonal complement of $\Cal D_z$ is taken w.r.to 
the basic canonical deformation metric $g^{\can}_1$ (see \S 2). 
The family of Z-metrics is defined in the following way. Let $\{\xi_i\}_{i=0}^3$ 
be a local frame of $\Cal D$ (here we are using the column $n$-vector notation 
identifying $\Cal D_z$ with $T_{\pi(z)}M=\H^n$, as in [CY]) satisfying a condition 
that $\{d\pi(\xi_i)\}_{i=0}^3$ gives an oriented orthonormal basis 
of $(T_{\pi(z)}M,g_{\pi(z)})$. Let $\{X^i\}_{i=0}^3$ denote the dual coframe 
which we extend to be zero on $\Cal D^{\perp}$. 
Let $\{\alpha_1\}_{i=1}^3$ be the $\sp(1)$-part 
of the Levi-Civita connection form of the original quaternion K\"ahler metric and 
$\{\alpha_1,\alpha_3\}$ the orthonormal coframe dual to the orthonormal frame 
along the twistor line. $\cong \PP^1=\Sp(1)/\Un(1)$. Set 
$\hat\alpha_j=\alpha_j-\sum_{i=0}^3\alpha_1(\xi_i)X^i$ ($j=1,3$). 
The Z-metric with parameter $\lambda$ is defined as declaring that 
$\{\lambda\hat\alpha_1,\lambda\hat\alpha_3,X^0,X^1,X^2,X^3\}$ 
forms an oriented orthonormal coframe : 
$$
g^{\ZZ}_{\lambda}=\lambda^2(\hat\alpha_1^2+\hat\alpha_3^2)+{}^tX^0\cdot X^0
+{}^tX^0\cdot X^1+{}^tX^0\cdot X^1+{}^tX^2\cdot X^2+{}^tX^3\cdot X^3\,\,.
$$
For curvature computation we must choose a good frame satisfying certain 
jet conditions at one point, which we briefly describe in the following. 
A point $z\in\Cal Z$ uniquely defines an orthogonal complex structure $J$ 
on $\Cal H_z\cong (T_mM,g_m)$ and its real $(2n)$-dimensional $J$-complex subspace 
spanned by $\xi_0$ and $J\xi_0$, where $\xi_0$ is determined uniquely modulo 
$\SO(2)$-rotation and should be understood as a column 
$n$ vector (and therefore $\{\xi_0,I\xi_0,J\xi_0,K\xi_0\}_{\text{\rm span}}$ is 
$(4n)$-dimensional)\footnote{\,\,The expression $\xi^i$ (resp. its dual $X^i$) for 
$i=0,1,2,3$ should be understood as raw (resp, column) $n$ vectors 
if they are used in the quaternion K\"ahler context related to the base manifold $M$ 
(e.g. in such a notation $\sum_{i=0}^3{}^tX^i\cdot X^i$ which represents 
the quaternion K\"ahler metric on $M$) or just a 1-form obtained by taking 
the sum in other occasion (e.g., in such a notation 
$\alpha_1-\alpha_1(\xi_1)X^1-\alpha_1(\xi_3)X^3$ 
which represents a certain modification of the $\sp(1)$-part of the 
connection form). For this notation is used in [C-Y]. See \S2 of this paper.}. 
We extend $\xi_0$ to a ``unit length'' germ at $m\in M$ 
so that $\nabla\xi_0$ satisfies certain condition. 
Moreover we extend $J$ to a section germ of $\Cal S$ 
so that $\nabla J=0$ at $m$, where $\Cal S$ is the 3-dimensional subbundle 
of $\text{\rm End}^{\text{\rm skew}}(TM)$ 
defining the quaternion K\"ahler structure of $M$. 
Next, we extend $I$ and $K$ to section germs of $\Cal S$ so that (i) $\{I,J,K\}$ 
constitutes the oriented orthonormal frame germ of $\Cal S$ (so $\{I,J,K\}$ 
satisfies the quaternion relations) and (ii) the oriented orthonormal frame 
germ $\sigma$ of the twistor fibration $\Cal Z \rightarrow M$ defined by 
$\{\xi_0,I\xi_0,J\xi_0,K\xi_0\}$ satisfies the condition 
$$(d\sigma)_m(T_mM)=\Cal D_z$$
(of course this is defined by a section germ of 
$\Cal P\rightarrow M$ composed with the projection $\Cal P \rightarrow \Cal Z$). 
Note that such a pair $\{I,K\}$ is defined uniquely modulo $\SO(2)$-rotation. 
This is the ``condition'' which should be satisfied by $I$ and $K$ 
in the construction of our Z-metrics\footnote{\,\,Important comparison : 
The construction of 
Z-metrics is characterized by the condition $(d\sigma)_m(T_mM)=\Cal D_z$. 
On the other hand the construction of the canonical deformation metrics is 
characterized by the condition $\nabla I=0$ and $\nabla K=0$ at $m$. }. 
Let the triple $\{\alpha_i\}_{i=1}^3$ be the $\sp(1)$-part of the connection form 
defined on the holonomy reduction $\Cal P$ of the oriented orthnormal frames. 
Let the orthogonal complex structure $J$ of $T_mM$ defined by an element of 
$\Cal S_m$ represented by $(0,1,0)$. The pair $\{\alpha_1,\alpha_3\}$ then 
represents the infinitesimal deformation of the orthogonal complex structures 
$J$ at $z\in\Cal Z$ and defines the induced horizontal subspace at $z$. 
Set $\xi_1=I\xi_0$, $\xi_2=J\xi_0$, $\xi_3=K\xi_0$. Let $\{X^0,X^1,X^2,X^3\}$ 
be the coframe field dual to $\{\xi_0,\xi_1,\xi_2,\xi_3\}$ (understood as canonical 
1-forms defined on $\Cal P$). 
The family $\Cal F^{\ZZ}$ of metrics on $\Cal Z$ is defined by declaring that 
$$
\split
& \hat\alpha_1:=\underbrace{\alpha_1-\alpha_1(\xi_0)X^0-\alpha_1(\xi_1)X^1
-\alpha_1(\xi_2)X^2-\alpha_1(\xi_3)X^3}_{\text{\rm annihilates $\Cal D_z$}}\,\,,\\
& \hat\alpha_3:=\underbrace{\alpha_3-\alpha_3(\xi_0)X^0-\alpha_3(\xi_1)X^1
-\alpha_3(\xi_2)X^2-\alpha_3(\xi_3)X^3}_{\text{\rm annihilates $\Cal D_z$}}\,\,,\\
& \underbrace{X^i \,\,(i=0,1,2,3)}_{\text{\rm dual to orthonormal basis of $\Cal D_z$}}
\endsplit
$$
forms an oriented orthonormal frame and therefore the Z-metric 
$g^{\ZZ}_{\lambda}$ is defined as
$$
\split
\rho\,g_{\lambda}^{\ZZ}&:=
\rho\,\biggl[
\lambda^2\,
\biggl\{\biggl(\alpha_1\underbrace{-\alpha_1(\xi_0)X^0}_{\text{\rm invisible}}
-\alpha_1(\xi_1)X^1\underbrace{-\alpha_1(\xi_2)X^2}_{\text{\rm invisible}}
-\alpha_1(\xi_3)X^3\biggr)^2\\
& \quad +\biggl(\alpha_3\underbrace{-\alpha_3(\xi_0)X^0}_{\text{\rm invisible}}
-\alpha_3(\xi_1)X^1\underbrace{-\alpha_3(\xi_2)X^2}_{\text{\rm invisible}}
-\alpha_3(\xi_3)X^3\biggr)^2\biggr\}\\
& \quad 
+\,{}^tX^0\cdot X^0+{}^tX^1\cdot X^1+{}^tX^2\cdot X^2+{}^tX^3\cdot X^3
\biggr]
\endsplit
$$
which is expressed in the orthonormal coframe germ at $z\in\Cal Z$ (one point)
\footnote{\,\,Although two 1-forms 
$\alpha_1-\alpha_1(\xi_0)X^0-\alpha_1(\xi_1)X^1-\alpha_1(\xi_2)X^2-\alpha_1(\xi_3)X^3$ 
and 
$\alpha_3-\alpha_3(\xi_0)X^0-\alpha_3(\xi_1)X^1-\alpha_3(\xi_2)X^2-\alpha_3(\xi_3)X^3$ 
are non orthogonal w.r.to the canonical deformation metric, we declare that they are 
orthogonal w.r.to the Z-metrics at $z\in\Cal Z$. 
Moreover, although the terms $\alpha_1(\xi_0)X^0$ and so on are ``invisible'' 
at $z\in\Cal Z$ we cannot ignore them because these become significant when we 
compute the Levi-Civita connection and curvature form by differentiation.}. 

We show that the subspace $\Cal F^{\ZZ}$ is foliated by the trajectories of 
the Ricci flow solution, i.e., the Ricci flow whose initial metric belongs 
to the subspace $\Cal F^{\ZZ}$ stays in $\Cal F^{\ZZ}$ 
as long as the solution exists and is an ancient solution in the sense of Hamilton [H]. 
The trajectory consisting of positive scalar multiples of the Einstein metric 
$g^{\ZZ}_{\sqrt{\frac1{n+2}}}$ is the asymptotic soliton of all other Ricci flow 
trajectories in $\Cal F^{\ZZ}$. 
In particular, the Einstein metric $g^{\ZZ}_{\sqrt{\frac1{n+2}}}$ 
is an ``unstable fixed point" under the dynamical system defined by 
the trajectories of the Ricci flow solutions. 
The Ricci flow solution starting at a metric in $\Cal F^{\ZZ}$ with $\lambda^2>\frac1{n+2}$ 
Gromov-Hausdorff converges (after appropriate scalings) in finite time to the 
Carnot-Carath\'eodory metric defined on the $(4n)$-dimensional distribution $\Cal D$ 
in $\Cal Z$, which isometrically projects to the original quaternion K\"ahler metric 
on the base manifold $M$. We then apply Bando-Shi estimate ([B], [Sh1,2]) 
for $\nabla \Rm$ under the Ricci flow to these 
ancient solutions to prove a limit formula which implies the LeBrun-Salamon conjecture 
that any irreducible positive quaternion K\"ahler manifold $(M,g)$ 
is isometric to one of the Wolf spaces. 
The technical part of the proofs of all results in this paper is based on 
the moving frame computation on the twistor space (see [C-Y]). 
We use Alexeevskii's curvature formula (Theorem 2.3 in this paper) (see [A] and [S]) 
in an essential way. 

Summing up, our main results are the following:

\proclaim{Theorem 1.4} (1) The Ricci tensor of the Z-metric 
$$g^{\ZZ}_{\lambda}=\lambda^2(\hat\alpha_1^2+\hat\alpha_3^2)
+\sum_{i=0}^3{}^tX^i\cdot X^i$$
is given by the formula
$$\Ric_{\lambda}^{\ZZ}=\frac4{\lambda^2}\lambda^2(\hat\alpha_1^2+\hat\alpha_3^2)
+(4n+8)\sum_{i=0}^3{}^tX^i\cdot X^i\,\,.$$
In particular $\rho g_{\lambda}^{\ZZ}$ is Einstein if and only if $\lambda^2=\frac1{n+2}$. 
\medskip

(2) The Ricci flow preserves the family $\R_+\cdot \Cal F^{\ZZ}$ of Z-metrics. 
The space $\R_+\cdot \Cal F^{\ZZ}$ is foliated by the trajectories 
of the Ricci flow solutions and these are all ancient solutions. 
Suppose that the Ricci flow $g(t)$ in $\R_+\cdot \Cal F^{\ZZ}$ is defined on $(-\infty,T)$. 
Then, modulo scaling, the solution $g(t)$ converges to the Einstein metric 
$g^{\ZZ}_{\sqrt{\frac1{n+2}}}$ as $t\to-\infty$. Moreover, modulo scaling, the limit $t\to T$
 corresponds to the ``collapse'' of $g^{\ZZ}_{\lambda}$ where the $X^i$-direction shrinks faster. 
\endproclaim

Applying Bando / Shi estimate to the ancient solutions in Theorem 1.4 (2), we have the 
limit formula
$$\lim_{\lambda\to\infty}|\nabla^{g^{\ZZ}_{\lambda}} 
\Rm^{g^{\ZZ}_{\lambda}}|_{g^{\ZZ}_{\lambda}}=0\,\,.$$

We will show that this limit formula implies the positive answer to the LeBrun-Salamon 
conjecture (Conjecture 1.2) claiming that any locally irreducible compact positive quaternion 
K\"ahler manifold is isometric to one of the Wolf spaces. 
\medskip

Finally, we remark that the methods used in this paper is a generalization of those used in 
[K-O] in the case of self-dual positive Einstein 4-manifolds. 
\bigskip
\noindent
\bigskip
\noindent
{\bf \S 2. Moving Frame Description of Quaternion K\"ahler Manifolds and their 
Twistor Spaces. Canonical Deformation Metrics.}
\medskip

Throughout this paper we will use the Einstein summation convention. 

We denote the quaternions by $\H$ and identify $\R^{4n}=\H^n$. Then $\H$ operates 
on $\R^{4n}=\H^n$ from the right which makes $\R^{4n}=\H^n$ a right $\H$-module. 
We then define a subgroup $\Sp(n)$ of $\SO(4n)$ as
$$\Sp(n)=\{A\in \SO(4n)\,|\,\text{\rm $A$ is $\H$-linear}\}\,\,.$$
The image in $\SO(4n)$ of the right action of the group $\Sp(1)$ of unit quaternions 
on $\H^n$ forms a subgroup of $\SO(4n)$ which we denote by $\Sp(1)$ by abuse of 
notations. 
Then we define the subgroup $\Sp(n)\Sp(1)$ of $\SO(4n)$ to be the product 
of the subgroups $\Sp(n)$ and $\Sp(1)$ in $\SO(4n)$. 
If $n>1$, the group $\Sp(n)\Sp(1)$ is a proper subgroup of $\SO(4n)$, while if $n=1$, 
the group $\Sp(1)\Sp(1)$ coincides with $\SO(4)$. 
For this reason, we assume $n>1$ from here on.

\proclaim{Definition 2.1} A $4n$-dimensional Riemannian manifold is 
quaternion K\"ahler if its holonomy group is contained in $Sp(n)Sp(1)$. 
\endproclaim

A $4n$-dimensional Riemannian manifold $(M^{4n},g)$ is quaternion K\"ahler if 
and only if the principal $\SO(4n)$-bundle $\Cal F$ of oriented orthonormal frames 
reduces to an $\Sp(n)\Sp(1)$-bundle $\Cal P$ (it is called the {\it holonomy reduction} 
of the oriented orthonormal frame bundle). 
Let $(e_A)_{A=1}^{4n}\in\Cal P$ 
be an orthonormal frame of $M$ at $m$. Then we have an identification
$$M_m \rightarrow \H^n\,\,,\,\,
x^Ae_A \mapsto (\,x^a+ix^{n+a}+jx^{2n+a}+kx^{3n+a}\,)_{a=1}^{n}\,\,.$$
Therefore each tangent space $M_m$ becomes a right $\H$-module. A local section 
$(e^A)$ of $\Cal P \rightarrow M$ on an open set $U\subset M$ 
defines a right $\H$-module structure on $TU$. However, this does not induce a global 
right $\H$-module structure on $TM$ because of the existence of the $\Sp(1)$ part in 
$\Sp(n)\Sp(1)$. 

The right action of $i$ and $j$ on $\R^{4n}=\H^n$ are given by the matrices
$$\pmatrix 0&-1&0&0\\ 1&0&0&0\\ 0&0&0&1\\ 0&0&-1&0\endpmatrix\quad
\text{\rm and}\quad 
\pmatrix 0&0&-1&0\\ 0&0&0&-1\\ 1&0&0&0\\ 0&1&0&0\endpmatrix$$
where $0$ (resp. $1$) denotes $n\times n$ zero (resp. identity) matrix. 
Therefore, the Lie algebra $\sp(n)$ is computed as
$$
\pmatrix A_0&-A_1&-A_2&-A_3&\\ A_1&A_0&-A_3&A_2\\ A_2&A_3&A_0&-A_1\\ 
A_3&-A_2&A_1&A_0\endpmatrix$$
where $A_0=-{}^tA_0$ and $A_{\lambda}={}^tA_{\lambda}$ ($1\leq\lambda\leq 4$) are 
$n\times n$ matrices. Similarly, the Lie algebra of the subgroup $\Sp(1)$ of $\SO(4n)$ 
is computed as 
$$
\pmatrix 0&-a_1&-a_2&-a_3\\ a_1&0&a_3&-a_2\\ a_2&-a_3&0&a_1\\ a_3&a_2&-a_1&0
\endpmatrix$$
where $a_1$, $a_2$, $a_3$ are $n\times n$ scalar matrices. 
Therefore the Lie algebra $\sp(n)\sp(1)$ is expressed as
$$
\pmatrix A_0&-A_1-a_1&-A_2-a_2&-A_3-a_3&\\ A_1+a_1&A_0&-A_3+a_3&A_2-a_2\\ 
A_2+a_2&A_3-a_3&A_0&-A_1+a_1\\ A_3+a_3&-A_2+a_2&A_1-a_1&A_0\endpmatrix
\tag2-1$$
where $A_0=-{}^tA_0$ and $A_{\lambda}={}^tA_{\lambda}$ ($1\leq\lambda\leq 4$) are 
$n\times n$ matrices and $a_1$, $a_2$, $a_3$ are $n\times n$ scalar matrices. 
In the following moving frame computation, we will use $\alpha_1$ instead of $a_i$ 
($i=1,2,3$) to represent $\sp(1)$-valued 1-forms. 
\medskip

Let $(M^{4n},g)$ be a (locally irreducible) quaternion K\"ahler manifold, i.e., its 
holonomy group is subgroup of $\SO(4n)$ contained in $\Sp(n)\Sp(1)$. 
Let $\Cal P\rightarrow M$ be the holonomy reduction of the bundle of oriented orthonormal 
frames of $M$. Then $\Cal P$ is a principal $\Sp(n)\Sp(1)$ bundle and we say a point of $\Cal P$ 
(i.e., an orthonormal frame at some point of $M$) a quaternion orthonormal frame. 
The moving frame description of a quaternion K\"ahler manifold is the following. 
Each local quaternion orthonormal frame field $e=(e_A)_{A=1}^{4n}$ 
on an open set $U\subset M$ defines a section 
$e:U\rightarrow \Cal P$ and therefore its dual coframe field $\theta=(\theta^A)_{A=1}^{4n}$ 
is pulled back to $e(U)\subset \Cal P$ via the restriction $\pi|_{e(U)}:e(U)\rightarrow U$ 
of the projection $\Cal P\rightarrow M$ to $e(U)$. 
We thus get a canonical system of 1-forms $(\theta^A)_{A=1}^{4n}$ on $\Cal P$ 
such that at each $e=(e_A)_{A=1}^{4n}\in\Cal P$, $(\theta^A)_{A=1}^{4n}$ coincides 
with its dual coframe pulled back to $\Cal P$. The Levi-Civita connection $\nabla$ 
of $(M,g)$ defines a unique right 
$\Sp(n)\Sp(1)$-invariant $\sp(n)\sp(1)$-valued 1-form 
$\Gamma=(\Gamma^A_{\,\,\,B})$ on $\Cal P$ 
satisfying
$$\nabla_X e_B=e_A\,\Gamma^A_{\,\,\,B}(X)$$
for each local quaternion orthonormal frame $(e_A)_{A=1}^{4n}$ and a local 
vector field $X$ on $M$. Here, 
$\Gamma=(\Gamma^A_{\,\,\,B})$ can be written as
$$(\Gamma^A_B)=\pmatrix \Gamma_0&-\Gamma_1-\alpha_1&-\Gamma_2-\alpha_2& 
-\Gamma_3-\alpha_3&\\ \Gamma_1+\alpha_1&\Gamma_0&-\Gamma_3+\alpha_3& 
\Gamma_2-\alpha_2\\ \Gamma_2+\alpha_2&\Gamma_3-\alpha_3&\Gamma_0& 
-\Gamma_1+\alpha_1\\ \Gamma_3+\alpha_3&-\Gamma_2+\alpha_2& 
\Gamma_1-\alpha_1&\Gamma_0\endpmatrix
\,\,=:\Gamma\tag2-2$$
where  $\Gamma_0=-{}^t\Gamma_0$ and $\Gamma_{\lambda}={}^t\Gamma_{\lambda}$ 
($1\leq\lambda\leq 4$) are $n\times n$ matrix valued 1-forms and $\alpha_1$, $\alpha_2$, 
$\alpha_3$ are $n\times n$ scalar matrix valued 1-forms. 
The 1-forms $\theta^A$ and $\Gamma^A_{\,\,\,B}$ satisfy the following first and second 
structure equations:
$$\split
& d\theta^A+\Gamma^A_{\,\,\,B}\wedge\theta^B=0\,\,,\\
& d\Gamma^A_{\,\,\,B}+\Gamma^A_{\,\,\,C}\wedge\Gamma^C_{\,\,\,B}=\Om^A_{\,\,\,B}\,\,,
\endsplit$$
where $\Om^A_{\,\,\,B}$ is the skew symmetric matrix of 2-forms which is identified with 
the curvature 2-form of $(M,g)$ in the following way. 
For each point $e=(e_A)\in\Cal P$ on $m\in M$ we have
$$R(X,Y)e_B=e_K\,\Om^K_{\,\,\,B}(X,Y)$$
for all $X,Y\in M_m$. Therefore the sectional curvature for the 2-plane spanned by 
$\{e_A,e_B\}$ is given by
$$K(e_A,e_B)=g(R(e_A,e_B)e_B,e_A)=g(\Om^K_{\,\,\,B}(e_A,e_B)e_K,e_A)\,\,.$$

We proceed to the description of the twistor space of a quaternion K\"ahler manifold 
$(M^{4n},g)$. A unit quaternion $q$ is pure imaginary if and only if $q^2=-1$ holds. 
Therefore, from the definition of $\Cal P$, a choice of a quatenion orthonormal frame 
$(e_A)_{A=1}^{4n}\in\Cal P$ of $M_m$ canonically defines the identification
$$\multline
\{\text{\rm unit pure imaginary quaternions}\}\\
\overset{\text{\rm right action}}\to\longleftrightarrow\\
\{\text{\rm orthogonal complex structures on $M_m$}\}\,\,.
\endmultline$$
This identification itself depends on the basis $(e_A)\in\Cal P$. However, 
if $q$ is a unit pure imaginary quaternion, then so is $x^{-1}qx$ for any unit 
quaternion $x$ and therefore the set $\PP^1$ of all orthogonal complex structures 
on $M_m$ is independent of the choice of the basis $(e_A)\in\Cal P$. 
The twistor space $\Cal Z$ of $M$ is by definition the fiber bundle over $M$ 
consisting of all orthogonal complex structures of all tangent spaces of $M$.  
Therefore the twistor space is described as the associated bundle
$$\Cal Z=\Cal P\times_{\Sp(n)\Sp(1)}\PP^1$$
where $\Sp(n)\Sp(1)$ operates on the set $\PP^1$ of unit pure imaginary quaternions 
by the trivial action of $\Sp(n)$ and the canonical right action of the group $\Sp(1)$ 
of unit quaternions given by $q \mapsto x^{-1}qx$. 
It follows from the definition of the twistor space that $\Cal Z$ has a canonical almost 
complex structure. Indeed, the Levi-Civita connection of $(M,g)$ induces the horizontal 
distribution on the twistor fibration $\Cal Z \rightarrow M$. 
We then have the canonical complex structure (defined by identifying the set of 
unit imaginary quaternions with $\PP^1$ by the stereo graphic projection) 
on each $\PP^1$-fiber.  Moreover, we associate 
to each horizontal subspace of  $\text{\rm H}(\Cal Z_z)$ ($z\in\Cal Z$ lies 
over $m\in M$), the almost complex structure of $M_m$ represented by $z\in \Cal Z$. 
We have thus defined canonically an almost complex structure on $\Cal Z$ which we call 
the orthogonal complex structure on $\Cal Z$. Moreover the Levi-Civita horizontal 
distribution on the twistor fibration $\Cal Z\rightarrow M$ canonically defines the sum of 
the scaled fiber Fubini-Study metric and the fixed base metric on $M$ which are expressed 
as 
$$g_{\lambda}^{\can}=\lambda^2g_{\FS}+g_M\,\,.\tag2-3$$
This class of metrics are called the canonical deformation metrics. 

Salamon [S] proved the following fundamental result:

\proclaim{Theorem 2.2} (1) The orthogonal almost complex structure on the twistor space 
$\Cal Z$ is integrable. 
\medskip

(2) Suppose that the scalar curvature of $(M,g)$ is positive. Then there is a unique 
scaling of the fiber metric such that the canonical deformation metric on the twistor 
space $\Cal Z$ is K\"ahler-Einstein with positive scalar curvature. 
\endproclaim

We start by recalling the idea of the Cartan formalism of moving frames. 
Let $(N,g)$ be any $n$-dimensional oriented Riemannian manifold 
and $\Cal F\rightarrow N$ the bundle of all oriented orthonormal frames. 
We have the system $\{\theta^1,\dots,\theta^n\}$ of coframes on $\Cal F$ 
which is, at $p\in \Cal F$ lying over $m\in N$, the system of 1-forms dual 
to the orthonormal frame of $N_m$ represented by the point $p\in\Cal F$. 
Given a local frame field on an open set $U\subset N$, we tautologically 
associate the section $U\rightarrow \Cal F$. Thus the local frames which are not 
unique on $N$ becomes a globally defined single valued object on $\Cal F$ 
and moreover the dual object $\{\theta_1,\dots,\theta_n\}$ consists of differential 
1-forms and therefore we have an advantage being able to work functorially on 
differential forms (such as connection forms) on $\Cal F$. 
For instance, the Riemannian metric on $N$ is written as $(\theta^1)^2+\cdots+(\theta^n)^2$ 
and connection form is computed by taking the exterior differential of 
$\{\theta_1,\dots,\theta_n\}$ on $\Cal F$ and so on. 

Now let us return to our original (quaternion K\"ahler) situation. 
A fiber on $m\in M$ of the twistor fibration $\Cal Z \rightarrow M$ is the set of all 
orthogonal complex structures on the tangent space $M_m$ which is 
canonically identified with $\Sp(n)\Sp(1)/\Sp(n)\Sp(1)\cap \Un(2n) \cong \PP^1$. 
Therefore the twistor space is also defined as the orbit space with respect to the 
$\Sp(n)\Sp(1)\cap \Un(2n)$ action on $\Cal P$, i.e.,
$$\Cal Z=\Cal P/\Sp(n)\Sp(1)\cap \Un(2n)\,\,.$$
We construct local sections $\Cal Z \rightarrow \Cal P$ of the principal 
$\Sp(n)\Sp(1)\cap \Un(2n)$-bundle $\Cal P\rightarrow \Cal Z$ in the following way 
(we use these local sections to construct a certain class of metrics on $\Cal Z$). 
Fix a point $m\in M$. 
Let $\PP^1_m\subset \Cal Z$ be the fiber of the twistor fibration over $m$. 
To each $z \in\PP^1_m$ we (locally) associate a quaternion orthonormal frame 
in the fiber of $\Cal P\rightarrow M$ over $m$ so that the frame is ordered in the way 
compatible with respect to the orthogonal complex structure represented by $z$. 
If $z$ varies on $\PP_m^1$ such frames rotates by an element of $\Sp(n)\Sp(1)$ 
and the rotation is unique modulo those by elements of $\Sp(n)\Sp(1)\cap\Un(2n)$. 
This procedure is possible only locally on $\PP^1_m$ because this is equivalent to make 
the (local) section of the principal $\Sp(n)\Sp(1)\cap\Un(2n)$-bundle 
$\Sp(n)\Sp(1) \rightarrow \Sp(n)\Sp(1)/\Sp(n)\Sp(1)\cap \Un(2n)\cong\PP^1$. 
We extend this construction 
locally on a small open set $U\subset M$ containing $m$ in such 
a way that the extended 
object is parallel at $m$ (one point). 
We perform this procedure at each point $m\in M$. 
This way we have constructed local 
sections of the $\Sp(n)\Sp(1)\cap \Un(2n)$-principal bundle 
$\Cal P \rightarrow \Cal Z$. We then pull back the canonical 1-forms
$$X^0,X^1,X^2,X^3$$
and the $\sp(1)$-part of the connection form
$$\alpha_1,\alpha_2,\alpha_3$$
by the above constructed local sections. 
We thus get the system of 1-forms 
$$\{X^0,X^1,X^2,X^3,\alpha_1,\alpha_3\}\,\,$$
locally at 1 point on $Z$. Then the basic canonical deformation metric 
in an expression in terms of the orthnormal coframes is expressed as 
(see [C-Y]) 
$$g_{1}^{\can}:=(\alpha_1^2+\alpha_3^2)
+{}^tX^0X^0+{}^tX^1X^1+{}^tX^2X^2+{}^tX^3X^3\,\,.$$
In the following arguments, we will consider the metrics of the form 
$$g_{\lambda}^{\can}:=\lambda^2(\alpha_1^2+\alpha_3^2)
+{}^tX^0X^0+{}^tX^1X^1+{}^tX^2X^2+{}^tX^3X^3
\tag2-4$$
on $\Cal Z$ (canonical deformation metrics). 
\medskip

In the following discussion, we compute the curvature form of the canonical deformation 
metrics $g_{\lambda}^{\can}$ defined on $\Cal Z$. 
We use the moving frame computation. The scalings of the various 
standard metrics are hidden in the computation. 
To avoid confusion, we fix our scaling convention in the following way : 
\medskip
\noindent
$\bullet$ We fix the scale of the invariant metric of $\H\PP^n$ so that the sectional 
curvatures range in the interval $[1,4]$, i.e., $\Ric(g_{\H\PP^n})=4(n+2)g_{\H\PP^n}$ 
and so $\text{\rm Scal}(g_{\H\PP^n})=16n(n+2)$. Namely we set 
$$\wt S=16n(n+2)$$
from here on. 
\medskip
\noindent
$\bullet$ We fix the scale of the Fubini-Study metric of the $\PP^1$-fiber of the twistor 
fibration and other cases so that the Gaussian curvature is identically $4$. 
\medskip

We will consider the following two types of scalings. 
\medskip
\noindent
$\bullet$ We fix the $\sum_{i=0}^3{}^tX^iX^i$-part and vary the ratio $S/\wt S$. 
In this case the above metric turns out to be K\"ahler (eventually K\"ahler-Einstein) 
on $\Cal Z$ if and only if $S/\wt S=1$. 
We will use this scaling in the computation in the transversal K\"ahler situation. 
\medskip
\noindent
$\bullet$ We normalize the base quaternion K\"ahler metric $g$ so that $S=\wt S$ holds 
and scale the 
$\sum_{i=0}^3{}^tX^iX^i$-part by the scaling parameter 
$\lambda^2$ (so that the ``curvature becomes $\lambda^{-2}$-times the original one in this 
direction"). 
In this case the metric $g_{\lambda}^{\can}$ turns out to be K\"ahler 
(eventually K\"ahler-Einstein) 
on $\Cal Z$ if and only if $\lambda=1$. 
\medskip

From here on until the end of \S2, we describe the moving frames 
on the twistor space $\Cal Z$ with the above introduced Z-metrics 
(and with complex structure in Theorem 2.2 if necessary). 
\medskip

Let $(e_A)\in\Cal P$ be a quaternion orthonormal basis of $M_m$. 
This canonically defines an identification $M_m$ with
$$\R^{4n}=\H^n=\C^{2n}$$
by
$$(x^a+ix^{n+a}+jx^{2n+a}+kx^{3a+n})_{a=1}^n \,\, \longleftrightarrow \,\, 
(x^a+jx^{2n+a},\,\,x^{n+a}+jx^{3n+a})_{a=1}^n\,\,.$$
The multiplication of $j$ from the right on $M_m=\C^{2n}$ induces the canonical 
almost complex structure $J$ of $\C^{2n}$. 
The infinitesimal deformation (``unit velocity vector tangent to a 1-parameter deformation") 
of unit imaginary quaternions at $j$ is expressed 
as $\alpha_1i+\alpha_3k$ (this is because the tangent space of the set of unit 
imaginary quaternions at $q$ is given by the orthogonality condition 
$\{\sigma\,|\,\Re(q\sigma)=0\}$). In this situation we pick a point $z\in\Cal Z$ 
which induces the almost complex structure on $M_m=\C^{2n}$ corresponding 
to $J$. The canonical almost complex structure of $\Cal Z$ at $z$ over $m$ 
is defined by specifying the basis of all $(1,0)$-forms as follows. 
$$\split
& \zeta^0=\alpha_1+i\alpha_3\,\,,\\
& \zeta^a=x^a+ix^{2n+a}\,\,,\\
& \zeta^{n+a}=x^{n+a}+ix^{3n+a}\,\,.
\endsplit
\tag2-5
$$
For moving frame computation, we introduce the column vectors 
$X^0,X^1,X^2$ and $X^3$ by setting $X^0=(x^a)$, $X^1=(x^{n+a})$, 
$X^2=(x^{2n+a})$ and $X^3=(x^{3n+a})$ where $a=1,\dots,n$. 
Then a basis of the $\C$-vector space of all $(1,0)$-forms is simply given by
$$\split
& \zeta^0=\alpha_1+i\alpha_3\,\,,\\
& Z^1=X^0+iX^2\,\,[=(\zeta^a)]\,\,,\\
& Z^2=X^1+iX^3\,\,[=(\zeta^{n+a})]\,\,.
\endsplit$$
The above argument involving $X^j$ ($j=0,1,2,3$) and $\alpha_i$ ($i=1,3$) 
is the local explanation of the construction of the section of the 
$\Sp(n)\Sp(1)\cap \Un(2n)$-principal bundle 
$\Cal P\rightarrow \Cal Z$ which described above in an abstract way. 

We would like to describe the curvature tensor of $\Cal Z$. To do so, we 
compute the derivation formula on $\Cal Z$ in terms of this basis of $(1,0)$-forms 
on $\Cal Z$. Set ${}^tX=({}^tX^0,{}^tX^1,{}^tX^2,{}^tX^3)$ (${}^tX^i$'s being 
row vectors). 
Then the first structure equation
$$dX+\Gamma\wedge X=0$$
is equivalent to
$$\left\{\aligned
& dZ^1+\ov Z^2\wedge\zeta^0+(\Gamma_0+i(\Gamma_2+\alpha_2))\wedge Z^1
+(-\Gamma_1+i\Gamma_3)\wedge Z^2=0\\
& dZ^2-\ov Z^1\wedge\zeta^0+(\Gamma_1+i\Gamma_3)\wedge Z^1
+(\Gamma_0-i(\Gamma_2-\alpha_2))\wedge Z^2=0
\endaligned\right.$$
which in matrix form is expressed as
$$\split
d\pmatrix \zeta^a\\ \zeta^{n+a}\endpmatrix&=
\underbrace{-\pmatrix \Gamma_0+i\Gamma_2 & -\Gamma_1+i\Gamma_3\\ 
\Gamma_1+i\Gamma_3 & \Gamma_0-i\Gamma_2\endpmatrix
\pmatrix \zeta^a\\ \zeta^{n+a}\endpmatrix}_{\text{\rm $\sp(n)$-action}}\\
& \quad \underbrace{-\pmatrix i\alpha_2 & 0\\ 0 & i\alpha_2\endpmatrix
\pmatrix \zeta^a\\ \zeta^{n+a}\endpmatrix
-\pmatrix 0 & -(\alpha_1+i\alpha_3)\\ \alpha_1+i\alpha_3 & 0\endpmatrix
\pmatrix \ov\zeta^a\\ \ov\zeta^{n+a}\endpmatrix}_{\text{\rm $\sp(1)$-action}}
\endsplit\tag2-6$$
What we must do next is to compute $d\zeta^0$. 
We need the second structure equation
$$d\Gamma+\Gamma\wedge\Gamma=\Om\,\,.$$
The expression (2-2) for the connection form implies the following:
$$
\Om=\pmatrix
\Om^0_0&\Om^0_1&\Om^0_2&\Om^0_3\\ \Om^1_0&\Om^0_0&\Om^1_2&\Om^1_3\\ 
\Om^2_0&\Om^2_1&\Om^0_0&\Om^2_3\\ \Om^3_0&\Om^3_1&\Om^3_2&\Om^0_0
\endpmatrix$$
where ${}^t\Om^0_0=-\Om^0_0$, $\Om^{\mu}_{\nu}={}^t\Om^{\mu}_{\nu}
=-\Om^{\nu}_{\mu}$. 
Let $(\mu,\eta,\nu)$ be any cyclic permutation of $(1,2,3)$. 
Then we have
$$\split
\Om^0_0&=d\Gamma_0+\Gamma_0\wedge\Gamma_0-\Gamma_1\wedge\Gamma_1
-\Gamma_2\wedge\Gamma_2-\Gamma_3\wedge\Gamma_3\,\,,\\
\Om^{\mu}_0&=d(\Gamma_{\mu}+\alpha_{\mu})
+(\Gamma_{\mu}+\alpha_{\mu})\wedge\Gamma_0
+(\Gamma_{\eta}-\alpha_{\eta})\wedge(\Gamma_{\nu}-\alpha_{\nu})\\
&\quad +\Gamma_0\wedge(\Gamma_{\mu}+\alpha_{\mu})
+(-\Gamma_{\nu}+\alpha_{\nu})\wedge(\Gamma_{\eta}+\alpha_{\eta})\\
&=\underbrace{d\Gamma_{\mu}+\Gamma_{\mu}\wedge\Gamma_0
+\Gamma_{\eta}\wedge\Gamma_{\nu}+
\Gamma_{0}\wedge\Gamma_{\mu}
-\Gamma_{\nu}\wedge\Gamma_{\eta}}_{\text{\rm $\Gamma$-part}}\\
&\quad \underbrace{+d\alpha_{\mu}
-2\alpha_{\eta}\wedge\alpha_{\nu}}_{\text{\rm scalar part}}\,\,,\\
\Om^{\eta}_{\nu}&=d(-\Gamma_{\mu}+\alpha_{\mu})
+(\Gamma_{\eta}+\alpha_{\eta})\wedge(-\Gamma_{\nu}-\alpha_{\nu})
+(-\Gamma_{\mu}+\alpha_{\mu})\wedge\Gamma_0\\
&+(\Gamma_{\nu}-\alpha_{\nu})\wedge(\Gamma_{\eta}-\alpha_{\eta})
+\Gamma_0\wedge(-\Gamma_{\mu}+\alpha_{\mu})\\
&=\underbrace{-d\Gamma_{\mu}-\Gamma_{\eta}\wedge\Gamma_{\nu}
-\Gamma_{\mu}\wedge\Gamma_0
+\Gamma_{\nu}\wedge\Gamma_{\eta}
-\Gamma_0\wedge\Gamma_{\mu}
}_{\text{\rm $\Gamma$-part}}\\
&\quad \underbrace{+d\alpha_{\mu}
-2\alpha_{\eta}\wedge\alpha_{\nu}}_{\text{\rm scalar part}}\,\,.
\endsplit\tag2-7$$
Therefore, to compute $d\zeta^0$, we need to know the structure of the 
curvature tensor of a quaternion K\"ahler manifold. 
In fact, the curvature tensor of a quaternion K\"ahler manifold is very special as is 
described in the following Alekseevskii's decomposition theorem (see [A] and [S]):

\proclaim{Theorem 2.3} (1) A locally irreducible quaternion K\"ahler manifold is Einstein.
\medskip

(2) The curvature operator of a locally irreducible quaternion K\"ahler manifold 
($M^{4n},g)$ is of the form
$$\Om=(S/\wt S)\,\wt \Om+\Om'$$
where $\wt \Om$ is the curvature operator of $\H\PP^n$, $\wt S$ is the scalar curvature of 
$\H\PP^n$, $S$ is the scalar curvature of $M$ ($S$ is a constant), and
$$\Om'\in\text{\rm Sym}^2(\sp(n)) \subset \text{\rm Sym}^2(\Lambda^2T^*M)\,\,.$$
\endproclaim

As was declared in \S0, we restrict our attention to locally irreducible quaternion 
K\"ahler manifolds and we say just ``quaternion K\"ahler" instead of ``locally irreducible 
quaternion K\"ahler". 

The meaning of Theorem 2.3 is the following. If we compute the curvature operator of 
a quaternion K\"ahler manifold $(M^{4n},g)$ in terms of the quaternion orthonormal 
basis $(X^A)$, then the curvature operator decomposes into the scalar multiple 
of $\wt \Om$ and the remaining part. The first part is $(S/\wt S)\,\wt \Om$ where $\wt \Om$ is 
expressed exactly in the same form as the curvature operator of $\H\PP^n$ where $(X^A)$ 
is regarded as quaternion orthonormal for the canonical metric of $\H\PP^n$. 
The remaining part $\Om'$ then looks like a curvature operator of a hyper-K\"ahler manifold 
(in particular no scalar $\alpha$-part is involved). 

The $\Sp(n)\Sp(1)$-principal bundle $\Cal P$ of quaternion orthonormal frames 
of $\H\PP^n$ coincides with the group $\Sp(n+1)/\Z_2$ and therefore the curvature tensor 
of $\H\PP^n$ is computed from  the Maurer-Cartan equation on $\Sp(n+1)$ 
applied to the description of the quaternion projective space 
$$\H\PP^n=\frac{\Sp(n+1)}{\Sp(n)\times\Sp(1)}=\frac{\Sp(n+1)/\Z_2}{\Sp(n)\Sp(1)}$$
as a symmetric space. The subgroup $\Sp(n)\Sp(1)$ of $\Sp(n+1)$ defined by 
the Lie algebra embedding
$$
\split
& \pmatrix A_0&-A_1-a_1&-A_2-a_2&-A_3-a_3&\\ A_1+a_1&A_0&-A_3+a_3&A_2-a_2\\ 
A_2+a_2&A_3-a_3&A_0&-A_1+a_1\\ A_3+a_3&-A_2+a_2&A_1-a_1&A_0\endpmatrix 
\longmapsto\\
& \pmatrix 0&a_1&-a_3&a_2\\ -a_1&0&a_2&a_3\\ a_3&-a_2&0&a_1\\ -a_2&-a_3&-a_1&0
\endpmatrix
\oplus\pmatrix A_0&-A_1&A_3&-A_2&\\ A_1&A_0&-A_2&-A_3\\ -A_3&A_2&A_0&-A_1\\ 
A_2&A_3&A_1&A_0\endpmatrix\endsplit$$
where $a_{\mu}$'s are regarded as scalar $n\times n$ matrices in the LHS while 
in the RHS these are regarded just as scalars. Indeed, the Lie algebra 
$$\pmatrix 0&-a_1&-a_2&-a_3\\ a_1&0&a_3&-a_2\\ a_2&-a_3&0&a_1\\ a_3&a_2&-a_1&0
\endpmatrix \leftrightarrow a_1R_i+a_2R_j+a_3R_k$$
of $\sp(1)$ considered as a Lie subalgebra of $\so(4n)$ stems from the matrix expression 
of the right action of quaternions on $\H^n$. Let $R_q$ (resp. $L_q$) denote the right 
(resp. left) action of a quaternion $q$ on $\H^n$. The correspondence 
$$R_i\mapsto -L_i\,\,,\,\,R_j\mapsto -L_k\,\,,\,\,R_k\mapsto L_j$$
defines a Lie algebra isomorphism of $\sp(1)$ which converts the right action of $\H$ on 
$\H^n$ to the left action. This Lie algebra isomorphism converts the above expression to 
$$\pmatrix 0&a_1&-a_3&a_2\\ -a_1&0&a_2&a_3\\ a_3&-a_2&0&a_1\\ -a_2&-a_3&-a_1&0
\endpmatrix \leftrightarrow a_1(-L_i)+a_3L_j+a_2(-L_k)\,\,.$$
To obtain a general expression of the Lie algebra $\sp(n+1)$ we put
$$\pmatrix
0&a_1&-a_3&a_2&-{}^tX^0&-{}^tX^1&{}^tX^3&-{}^tX^2\\
-a_1&0&a_2&a_3&-{}^tX^{01}&-{}^tX^{11}&-{}^tX^{21}&-{}^tX^{31}\\
a_3&-a_2&0&a_1&-{}^tX^{02}&-{}^tX^{12}&-{}^tX^{22}&-{}^tX^{32}\\
-a_2&-a_3&-a_1&0&-{}^tX^{03}&-{}^tX^{13}&-{}^tX^{23}&-{}^tX^{33}\\
X^0&X^{01}&X^{02}&X^{03}&A_0&-A_1&A_3&-A_2\\
X^1&X^{11}&X^{12}&X^{13}&A_1&A_0&-A_2&-A_3\\
-X^3&X^{21}&X^{22}&X^{23}&-A_3&A_2&A_0&-A_1\\
X^2&X^{31}&X^{32}&X^{33}&A_2&A_3&A_1&A_0\endpmatrix$$
and determine $X^{\mu\nu}$'s from the commutativity with the right action 
of $i$ and $j$, namely, with the matrices
$$
\pmatrix 0&-1&0&0\\ 1&0&0&0\\ 0&0&0&1\\ 0&0&-1&0\endpmatrix \oplus 
\pmatrix 0&-E_n&0&0\\ E_n&0&0&0\\ 0&0&0&E_n\\ 0&0&-E_n&0\endpmatrix$$
and
$$
\pmatrix 0&0&-1&0\\ 0&0&0&-1\\ 1&0&0&0\\ 0&1&0&0\endpmatrix\oplus
\pmatrix 0&0&-E_n&0\\ 0&0&0&-E_n\\ E_n&0&0&0\\ 0&E_n&0&0\endpmatrix\,\,.$$
It follows that the Lie algebra $\sp(n+1)$ is expressed as
$$\pmatrix
0&a_1&-a_3&a_2&-{}^tX^0&-{}^tX^1&{}^tX^3&-{}^tX^2\\
-a_1&0&a_2&a_3&{}^tX^1&-{}^tX^0&-{}^tX^2&-{}^tX^3\\
a_3&-a_2&0&a_1&-{}^tX^3&{}^tX^2&-{}^tX^0&-{}^tX^1\\
-a_2&-a_3&-a_1&0&{}^tX^2&{}^tX^3&{}^tX^1&-{}^tX^0\\
X^0&-X^1&X^3&-X^2&A_0&-A_1&A_3&-A_2\\
X^1&X^0&-X^2&-X^3&A_1&A_0&-A_2&-A_3\\
-X^3&X^2&X^0&-X^1&-A_3&A_2&A_0&-A_1\\
X^2&X^3&X^1&X^0&A_2&A_3&A_1&A_0\endpmatrix\,\,.$$
Therefore, the following matrix represents a basis of left invariant 1-forms 
on $\Sp(n+1)$ :
$$\wt\Gamma_{\alpha,X}:=\pmatrix
0&\wt\alpha_1&-\wt\alpha_3&\wt\alpha_2&-{}^tX^0&-{}^tX^1&{}^tX^3&-{}^tX^2\\
-\wt\alpha_1&0&\wt\alpha_2&\wt\alpha_3&{}^tX^1&-{}^tX^0&-{}^tX^2&-{}^tX^3\\
\wt\alpha_3&-\wt\alpha_2&0&\wt\alpha_1&-{}^tX^3&{}^tX^2&-{}^tX^0&-{}^tX^1\\
-\wt\alpha_2&-\wt\alpha_3&-\wt\alpha_1&0&{}^tX^2&{}^tX^3&{}^tX^1&-{}^tX^0\\
X^0&-X^1&X^3&-X^2&\wt\Gamma_0&-\wt\Gamma_1&\wt\Gamma_3&-\wt\Gamma_2\\
X^1&X^0&-X^2&-X^3&\wt\Gamma_1&\wt\Gamma_0&-\wt\Gamma_2&-\wt\Gamma_3\\
-X^3&X^2&X^0&-X^1&-\wt\Gamma_3&\wt\Gamma_2&\wt\Gamma_0&-\wt\Gamma_1\\
X^2&X^3&X^1&X^0&\wt\Gamma_2&\wt\Gamma_3&\wt\Gamma_1&\wt\Gamma_0
\endpmatrix\,\,.$$
The Maurer-Cartan equation
$$d\wt\Gamma_{\alpha,X}+\wt\Gamma_{\alpha,X}\wedge\wt\Gamma_{\alpha,X}=0$$
implies the following. Let $(\mu,\eta,\nu)$ be any cyclic permutation of $(1,2,3)$. Then :
$$
\split
& d\wt\Gamma_0-X^0\wedge {}^tX^0-X^1\wedge {}^tX^1
-X^2\wedge {}^tX^2-X^3\wedge {}^tX^3\\
& \quad +\wt\Gamma_0\wedge\wt\Gamma_0
-\wt\Gamma_1\wedge\wt\Gamma_1-\wt\Gamma_2\wedge\wt\Gamma_2
-\wt\Gamma_3\wedge\wt\Gamma_3=0\,\,,\\
& d\wt\alpha_{\mu}-2\wt\alpha_{\eta}\wedge\wt\alpha_{\nu}-{}^tX^{\mu}\wedge X^0
+{}^tX^0\wedge X^{\mu}+{}^tX^{\nu}\wedge X^{\eta}-{}^tX^{\eta}\wedge X^{\nu}=0\,\,,\\
& d\wt\Gamma_{\mu}+\wt\Gamma_{\mu}\wedge\wt\Gamma_0
+\wt\Gamma_0\wedge\wt\Gamma_{\mu}-\wt\Gamma_{\nu}\wedge\wt\Gamma_{\eta}
+\wt\Gamma_{\eta}\wedge\wt\Gamma_{\nu}\\
&\quad -X^{\mu}\wedge{}^tX^0+X^0\wedge{}^tX^{\mu}-X^{\nu}\wedge{}^tX^{\eta}
+X^{\eta}\wedge{}^tX^{\nu}=0\,\,.
\endsplit$$
This and (2-4) implies that the curvature $(\wt\Om^{\mu}_{\nu})$ of $\H\PP^n$ 
is expressed as
$$
\split
\wt\Om^0_0&=d\wt\Gamma_0-\wt\Gamma_1\wedge\wt\Gamma_1
-\wt\Gamma_2\wedge\wt\Gamma_2-\wt\Gamma_3\wedge\wt\Gamma_3\\
&=X^0\wedge{}^tX^0+X^1\wedge{}^tX^1+X^2\wedge{}^tX^2+X^3\wedge{}^tX^3\,\,,\\
\wt\Om^{\mu}_0&=d\wt\Gamma_{\mu}+\wt\Gamma_{\mu}\wedge\wt\Gamma_0
+\wt\Gamma_0\wedge\wt\Gamma_{\mu}-\wt\Gamma_{\nu}\wedge\wt\Gamma_{\eta}
+\wt\Gamma_{\eta}\wedge\wt\Gamma_{\nu}\\
&\quad +d\wt\alpha_{\mu}-2\wt\alpha_{\eta}\wedge\wt\alpha_{\nu}\\
&=X^{\mu}\wedge{}^tX^0-X^0\wedge{}^tX^{\mu}+X^{\nu}\wedge{}^tX^{\eta}
-X^{\eta}\wedge{}^tX^{\nu}\\
&\quad +2({}^tX^{\mu}\wedge X^0+{}^tX^{\eta}\wedge X^{\nu})\,\,,\\
\wt\Om^{\eta}_{\nu}&=-d\wt\Gamma_{\mu}-\wt\Gamma_{\eta}\wedge\wt\Gamma_{\nu}
+\wt\Gamma_{\nu}\wedge\wt\Gamma_{\eta}-\wt\Gamma_0\wedge\wt\Gamma_{\mu}
-\wt\Gamma_{\mu}\wedge\wt\Gamma_0\\
&\quad +d\wt\alpha_{\mu}-2\wt\alpha_{\eta}\wedge\wt\alpha_{\nu}\\
&=-X^{\mu}\wedge{}^tX^0+X^0\wedge{}^tX^{\mu}-X^{\nu}\wedge{}^tX^{\eta}
+X^{\eta}\wedge{}^tX^{\nu}\\
&\quad +2({}^tX^{\mu}\wedge X^0+{}^tX^{\eta}\wedge X^{\nu})\,\,.
\endsplit\tag2-8$$
In particular this implies that the sectional curvatures of $\H\PP^n$ are, for instance, 
$K(e_1,e_a)=K(e_1,e_{n+a})=K(e_1,e_{2n+a})=K(e_1,e_{3n+a})=1$ ($2\leq a\leq n)$ and 
$K(e_1,e_{n+1})=K(e_1,e_{2n+1})=K(e_1,e_{3n+1})=4$ (i.e., $\H\PP^n$ is $1/4$-pinched). 
This implies that $\H\PP^n$ is Einstein with $\Ric(g)=4(n+2)g$. 

We now return to the computation of the structure equation with respect to the class 
of metrics on the twistor space $\Cal Z$ introduced above 
for general quaternion K\"ahler manifolds. 

We give a proof of Theorem 2.2 by following the arguments in [C-Y], 
since the proof of our main theorem is based on the moving frame proof of 
Theorem 2.2 in the technical level. 
Indeed, what is essential in the proof of our main theorem is the comparison 
of the canonical deformation and Z-metrics (which we introduce in \S3) 
on the twistor space $\Cal Z$ in terms of the Cartan formalism of the  moving frames. 
\medskip

We start with the case of $M=\H\PP^n$. We compute in the transversal K\"ahler setting. 
To compute $d\zeta^0$ ($\zeta^0=\alpha_1+i\alpha_3$), we use the formula
$$
d\wt\alpha_{\mu}-2\wt\alpha_{\eta}\wedge\wt\alpha_{\nu}=2({}^tX^{\mu}\wedge X^0
+{}^tX^{\eta}\wedge X^{\nu})\,\,.\tag2-9
$$
This implies
$$
d\wt\zeta^0=-2i\wt\alpha_2\wedge\wt\zeta^0+{}^tZ^2\wedge Z^1-{}^tZ^1\wedge Z^2\,\,.
$$
It follows from (2-3) and (2-6) the structure equation on the twistor space of $\H\PP^n$ :
$$
d\,\pmatrix \wt\zeta^0\\ Z^1\\ Z^2\endpmatrix
=-\pmatrix 2i\wt\alpha_2 & -{}^tZ^2 & {}^tZ^1\\ \ov Z^2 & \wt\Gamma_0
+i\wt\Gamma_2+i\wt\alpha_2 & 
-\wt\Gamma_1+i\wt\Gamma_3\\ -\ov Z^1 & \wt\Gamma_1+i\wt\Gamma_3 & 
\wt\Gamma_0-i\wt\Gamma_2+i\wt\alpha_2
\endpmatrix\,\wedge\,\pmatrix \wt\zeta^0\\ Z^1\\ Z^2\endpmatrix\,\,.\tag2-10$$

\medskip

We then turn to the general case, i.e., $M$ being a quaternion K\"ahler manifold, $\Cal Z$ 
(resp. $\wcalZ$) its twistor space (resp. extended twistor space). We first consider the 
transversal K\"ahler setting. 

We recall that $\zeta^0$, $Z^1$ and $Z^2$ is a basis of $(1,0)$-forms at 
a point of $\Cal Z$ w.r.to the canonical complex structure. It follows from (2-3) that
$$
d\,\pmatrix \zeta^0\\ Z^1\\ Z^2\endpmatrix
=-\pmatrix 2i\alpha_2&*&*\\ \ov Z^2 & \Gamma_0+i\Gamma_2+i\alpha_2 & 
-\Gamma_1+i\Gamma_3\\ -\ov Z^1 & \Gamma_1+i\Gamma_3 & \Gamma_0
-i\Gamma_2+i\alpha_2
\endpmatrix\,\wedge\,\pmatrix \zeta^0\\ Z^1\\ Z^2
\endpmatrix+\cdots$$
and therefore what we have to compute is to express 
$d\alpha_{\mu}-2\alpha_{\eta}\wedge\alpha_{\nu}$ in terms of $Z^{\mu}$'s 
($\zeta^0$ being $\alpha_1+i\alpha_3$). To compute $d\zeta^0$, 
we observe from (2-4) that
$$
\Om^{\mu}_0+\Om^{\eta}_{\nu}=2d\alpha_{\mu}-4\alpha_{\eta}\wedge \alpha_{\nu}\,\,.
\tag2-11$$
holds ($\mu,\eta,\nu$ being any cyclic permutation of $(1,2,3)$). 
Combining (2-8) with Theorem 2.3, we get
$$
\split
& d\alpha_{\mu}-2\alpha_{\eta}\wedge \alpha_{\nu}=\frac12(\Om^{\mu}_0
+\Om^{\eta}_{\nu})\\
& \quad =\frac12(S/\wt S)(\wt\Om^{\mu}_0+\wt\Om^{\eta}_{\nu})
+\frac12({\Om'}^{\mu}_0+{\Om'}^{\eta}_{\nu})\\
& \quad =(S/\wt S)(d\wt\alpha_{\mu}-2\wt\alpha_{\eta}\wedge\wt\alpha_{\nu})\quad 
[\text{\rm because $\Om'$ part does not involve the $\alpha$-part}]\\
& \quad =2(S/\wt S)({}^tX^{\mu}\wedge X^0+{}^tX^{\eta}\wedge X^{\nu})\,\,.
\endsplit
\tag2-12$$
The formula (2-12) implies that the $\SO(2)$-bundle $\wcalZ \rightarrow \Cal Z$ 
defines a Hermitian holomorphic line bundle on the twistor space $\Cal Z$ 
with its curvature form proportional to the K\"ahler-Einstein metric of $\Cal Z$. 
We get from formula (2-12) the formula
$$
\split
d\zeta^0&=d(\alpha_1+i\alpha_3)\\
&=2\alpha_2\wedge\alpha_3
+(d\alpha_1-2\alpha_2\wedge\alpha_3)+2i\alpha_1\wedge\alpha_2
+i(d\alpha_3-2\alpha_1\wedge\alpha_2)\\
&=-2i\alpha_2\wedge\zeta^0+(S/\wt S)({}^tZ^2\wedge Z^1-{}^tZ^1\wedge Z^2)
\endsplit
\tag2-13
$$
which means that the twistor space $\Cal Z$ is a {\it holomorphic contact manifold}. 
It follows that
$$
d\,\pmatrix \zeta^0\\ Z^1\\ Z^2\endpmatrix
=-\pmatrix 2i\alpha_2 & -(S/\wt S){}^tZ^2 & (S/\wt S){}^tZ^1\\ 
\ov Z^2 & \Gamma_0+i\Gamma_2+i\alpha_2 & -\Gamma_1+i\Gamma_3\\ 
-\ov Z^1 & \Gamma_1+i\Gamma_3 & \Gamma_0-i\Gamma_2+i\alpha_2
\endpmatrix\,\wedge\,\pmatrix \zeta^0\\ Z^1\\ Z^2\endpmatrix\,\,.\tag2-13$'$
$$
Let $\Gamma$ denote the matrix in the right hand side of (2-9). 
The formula (2-9) confirms that the almost complex structure defined by 
the basis $\{\zeta^0,Z^1,Z^2\}$ on $\Cal Z$ of the space of $(1,0)$-forms 
is integrable, since the right hand side contains no $(0,2)$-forms. 
Moreover, the matrix in (2-13$'$) becomes skew-Hermitian if and only if the metric of $M$ is 
scaled so that $S/\wt S=1$. This means that the Hermitian metric of $\Cal Z$ with 
the property that the basis  $\{\zeta^0,Z^1,Z^2\}$ of $(1,0)$-forms is unitary, 
is a K\"ahler metric of $\Cal Z$ if and only if the metric of $M$ is scaled so that 
$S/\wt S=1$. Now let the metric $M$ is scaled so that $S/\wt S=1$. 
The second structure equation 
$d\Gamma+\Gamma\wedge\Gamma=\Om$ 
computes the curvature $\Om$ of the K\"ahler metric 
$\zeta^0\wedge \ov\zeta^0
+{}^tZ^1\wedge\ov Z^1+{}^tZ^2\wedge \ov Z^2$. 
The direct computation shows the following : 
$$
\pmatrix
2\zeta^0\wedge\ov\zeta^0+{}^tZ^1\wedge\ov Z^1
& \zeta^0\wedge{}^t\ov Z^1
& \zeta^0\wedge{}^t\ov Z^2\\
+{}^tZ^2\wedge\ov Z^2 & & \\
&&&\\
Z^1\wedge\ov\zeta^0
& \Om^0_0+i\,\Om^2_0
& -\frac12\{\Om^1_0+\Om^3_2-i\,(\Om^3_0+\Om^2_1)\}\\
 & -\ov Z^2\wedge{}^tZ^2+\zeta^0\wedge\ov\zeta^0
& +\ov Z^2\wedge{}^tZ^1\\
&&&\\
Z^2\wedge\ov\zeta^0
& \frac12\{\Om^1_0+\Om^3_2+i\,(\Om^3_0+\Om^2_1)\}
& \Om^0_0+i\,\Om^3_1\\
 & +\ov Z^1\wedge{}^tZ^2
& -\ov Z^1\wedge {}^tZ^1+\zeta^0\wedge\ov\zeta^0
\endpmatrix$$
which is certainly skew-Hermitian. 
The Ricci form of the K\"ahler metric $\zeta^0\wedge \ov\zeta^0
+{}^tZ^1\wedge\ov Z^1+{}^tZ^2\wedge \ov Z^2$ (identified with $g_1^{\can}$) 
is given by
$$
\Ric(\Om)=\tr(\Om)=2(n+1)\biggl\{\zeta^0\wedge\ov\zeta^0
+{}^tZ^1\wedge \ov Z^1+{}^tZ^2\wedge\ov Z^2\biggr\}\,\,.$$
This implies that the K\"ahler metric $\zeta^0\wedge \ov\zeta^0
+{}^tZ^1\wedge\ov Z^1+{}^tZ^2\wedge \ov Z^2$ (this is $g_1^{\can}$) 
is K\"ahler-Einstein. This complete the proof of Theorem 2.2 \qed
\medskip

We are now ready to compute the Levi-Civita connection of the canonical deformation 
metrics. We have, at $z\in\Cal Z$ modulo terms vanishing to order ($\geq 2$), the following 
matrix expression : 
$$
d\pmatrix \lambda\alpha_1\\ \lambda\alpha_3\\ 
X^0\\ X^1\\ X^2\\ X^3\endpmatrix
+\pmatrix 
0&-2\alpha_2&-\lambda {}^tX^1&\lambda {}^tX^0&\lambda {}^tX3&-\lambda {}^tX^2\\
2\alpha_2&0&-\lambda {}^tX^3&\lambda {}^tX^2&-\lambda {}^tX^1&\lambda {}^tX^0\\
\lambda X^1&\lambda X^3
&\Gamma_0&-\Gamma_1^-&-\Gamma_2^+&-\Gamma_3^-\\
-\lambda X^0&-\lambda X^2
&\Gamma_1^-&\Gamma_0&-\Gamma_3^+&\Gamma_2^-\\
-\lambda X^3&\lambda X^1
&\Gamma_2^+&\Gamma_3^+&\Gamma_0&-\Gamma_1^+\\
\lambda X^2&-\lambda X^0
&\Gamma_3^-&-\Gamma_2^-&\Gamma_1^+&\Gamma_0
\endpmatrix
 \wedge 
\pmatrix \lambda\alpha_1\\ \lambda\alpha_3\\ X^0\\ X^1\\ X^2\\ X^3\endpmatrix
=\pmatrix 0\\0\\0\\0\\0\\0\endpmatrix
\tag2-14
$$
where $\Gamma_1^{\pm}=\Gamma_1\pm(\lambda^2-1)\alpha_1$, 
$\Gamma_3^{\pm}=\Gamma_3\pm(\lambda^2-1)\alpha_3$ and 
$\Gamma_2^{\pm}=\Gamma_2\pm(\lambda^2-1)\alpha_2$. 
The skew symmetric matrix in (2-14) is the the connection form of the Levi-Civita connection 
which we denote as $\Gamma_{\lambda}$.
(2-14) is the first structure equation for the canonical deformation metric $g^{\can}_{\lambda}$ 
on $\Cal Z$ and the matrix part should be understood as 1-form germs at $z\in\Cal Z$ 
modulo terms vanishing at $z$ to the order $\geq 2$. Putting $\zeta^0=\alpha_1+i\alpha_3$, 
$Z^1=X^0+iX^2$ and $Z^2=X^1+iX^3$, (2-14) is rewritten as
$$d\pmatrix \lambda\zeta^0\\ Z^1\\ Z^2\endpmatrix+
\pmatrix 2i\alpha_2 & -\lambda{}^tZ^2 & \lambda{}^tZ^1\\
(\lambda-\lambda+\lambda^{-1})\ov{Z}^2 & \Gamma_0+i\Gamma_2+i\alpha_2 & -\Gamma_1+i\Gamma_3\\
-(\lambda+\lambda-\lambda^{-1})\ov{Z}^1 & \Gamma_1+i\Gamma_3 & \Gamma_0-i\Gamma_2+i\alpha_2
\endpmatrix \pmatrix \lambda\zeta^0\\ Z^1\\ Z^2\endpmatrix=\pmatrix 0\\ 0\\ 0\endpmatrix
\tag2-14$'$
$$
The matrix in (2-14) is skew Hermitian if and only if $\lambda^2=1$, i.e., 
the metric $g_{\lambda}^{\can}$ is the basic canonical deformation metric $g_1^{\can}$ 
(this reproves Theorem 2.2). 

The second structure equation 
$d\Gamma_{\lambda}+\Gamma_{\lambda}\wedge\Gamma_{\lambda}=\Om_{\lambda}$ 
computes the curvature form $\Om^{\can}_{\lambda}$. 
Put
$$\Om_{\lambda}=\pmatrix 
{\Om_{\lambda}}^{-2}_{-2}&{\Om_{\lambda}}^{-2}_{-1}&{\Om_{\lambda}}^{-2}_{0}
&{\Om_{\lambda}}^{-2}_{1}&{\Om_{\lambda}}^{-2}_{2}&{\Om_{\lambda}}^{-2}_{3}\\
{\Om_{\lambda}}^{-1}_{-2}&{\Om_{\lambda}}^{-1}_{-1}&{\Om_{\lambda}}^{-1}_{0}
&{\Om_{\lambda}}^{-1}_{1}&{\Om_{\lambda}}^{-1}_{2}&{\Om_{\lambda}}^{-1}_{3}\\
{\Om_{\lambda}}^{0}_{-2}&{\Om_{\lambda}}^{0}_{-1}&{\Om_{\lambda}}^{0}_{0}
&{\Om_{\lambda}}^{0}_{1}&{\Om_{\lambda}}^{0}_{2}&{\Om_{\lambda}}^{0}_{3}\\
{\Om_{\lambda}}^{1}_{-2}&{\Om_{\lambda}}^{1}_{-1}&{\Om_{\lambda}}^{1}_{0}
&{\Om_{\lambda}}^{1}_{1}&{\Om_{\lambda}}^{1}_{2}&{\Om_{\lambda}}^{1}_{3}\\
{\Om_{\lambda}}^{2}_{-2}&{\Om_{\lambda}}^{2}_{-1}&{\Om_{\lambda}}^{2}_{0}
&{\Om_{\lambda}}^{2}_{1}&{\Om_{\lambda}}^{2}_{2}&{\Om_{\lambda}}^{2}_{3}\\
{\Om_{\lambda}}^{3}_{-2}&{\Om_{\lambda}}^{3}_{-1}&{\Om_{\lambda}}^{3}_{0}
&{\Om_{\lambda}}^{3}_{1}&{\Om_{\lambda}}^{3}_{2}&{\Om_{\lambda}}^{3}_{3}
\endpmatrix$$
Then we have
$$\split
& {\Om_{\lambda}}^{-2}_{-2}=0\,\,,\quad {\Om_{\lambda}}^{-1}_{-1}=0\,\,,\\
& {\Om_{\lambda}}^{-1}_{-2}=2d\alpha_2-\lambda^2\,{}^tX^3\wedge X^1
-\lambda^2 \,{}^tX^2\wedge X^0+\lambda^2\,{}^tX^1\wedge X^3
+\lambda^2\,{}^tX^0\wedge X^2\\
&\quad =4\alpha_3\wedge \alpha_1+(4-2\lambda^2)\,{}^tX^3\wedge X^1
+(4-2\lambda^2)\,{}^tX^2\wedge X^0
\endsplit$$
In the rest of the curvature computation we use (2-2) and the first structure equation 
$dX+\Gamma\wedge X=0$. 
A typical computation is
$$\split
{\Om_{\lambda}}^0_{-2}&=\lambda \, dX^1+\lambda X^3\wedge(2\alpha_2)
+\Gamma_0\wedge (\lambda X^1)
+(\Gamma_1+(1-\lambda^2)\alpha_1)\wedge (\lambda X^0)\\
&\quad +(\Gamma_2+\alpha_2)\wedge (\lambda X^3)
-(\Gamma_3+(1-\lambda^2)\alpha_3)\wedge (\lambda X^2)\}\\
& =\lambda\{-(\Gamma_1+\alpha_1)\wedge X^0-\Gamma_0\wedge X^1
+(\Gamma_3-\alpha_3) \wedge X^2-(\Gamma_2-\alpha_2)\wedge X^3\}\\
& \quad -2\lambda \alpha_2\wedge X^3+\lambda \Gamma_0\wedge X^1
+\Gamma_1\wedge X^0+\lambda(1-\lambda^2)\alpha_1\wedge X^0\\
&\quad +\lambda \Gamma_2\wedge X^3+\lambda \alpha_2\wedge X^3
-\lambda \Gamma_3\wedge X^2-\lambda(1-\lambda^2)\alpha_3\wedge X^2\\
& =\lambda^3 X^0\wedge \alpha_1+(2\lambda-\lambda^3) X^2\wedge\alpha_3\,\,.
\endsplit$$
By similar computations as above, we have the following expression of a part of the curvature : 
$$
\split
& {\Om_{\lambda}}^0_{-2}=\lambda^3 X^0\wedge \alpha_1
+(2\lambda-\lambda^3) X^2\wedge \alpha_3\,\,,\quad 
{\Om_{\lambda}}^{0}_{-1}=\lambda^3X^0\wedge \alpha_3
-(2\lambda-\lambda^3) X^2\wedge \alpha_1\,\,,\\
& {\Om_{\lambda}}^1_{-2}=\lambda^3 X^1\wedge \alpha_1
+(2\lambda-\lambda^3)X^3\wedge \alpha_3\,\,,\quad 
{\Om_{\lambda}}^{1}_{-1}=\lambda^3 X^1\wedge \alpha_3
-(2\lambda-\lambda^3) X^3\wedge \alpha_1\,\,,\\
& {\Om_{\lambda}}^2_{-2}=\lambda^3 X^2\wedge \alpha_1
-(2\lambda-\lambda^3) X^0\wedge \alpha_3\,\,,\quad 
{\Om_{\lambda}}^{2}_{-1}=\lambda^3 X^2\wedge \alpha_3
+(2\lambda-\lambda^2) X^0\wedge \alpha_1\,\,,\\
& {\Om_{\lambda}}^3_{-2}=\lambda^3 X^3\wedge \alpha_1
-(2\lambda-\lambda^3) X^1\wedge \alpha_3\,\,, \quad
{\Om_{\lambda}}^{3}_{-1}=\lambda^3  X^3\wedge \alpha_3
+(2\lambda-\lambda^3) X^1\wedge \alpha_1\,\,.
\endsplit$$
In the computation of ${\Om_{\lambda}}^{\mu}_{\nu}$ ($\mu,\nu\geq 0$), we use (2-4) and 
compare ${\Om_{\lambda}}^{\mu}_{\nu}$ and $\Om^{\mu}_{\nu}$. 
We have the following results, which should be understood modulo hyper-K\"ahler 
contribution, i.e., ${\Om_{\lambda}}_{\nu}^{\mu}$ should be understood 
modulo ${\Om'}_{\nu}^{\mu}$ :
$$\split
{\Om_{\lambda}}^0_0&={\Om_{\lambda}}^2_2
=\Om^0_0-\lambda^2(X^1\wedge{}^tX^1-X^3\wedge{}^tX^3)\\
&=X^0\wedge{}^tX^0+X^1\wedge{}^tX^1+X^2\wedge{}^tX^2+X^3\wedge{}^tX^3
-\lambda^2(X^1\wedge{}^tX^1+X^3\wedge {}^tX^3)\,\,,\\
{\Om_{\lambda}}^1_1&={\Om_{\lambda}}^3_3
=\Om^0_0-\lambda^2(X^0\wedge{}^tX^0-X^2\wedge{}^tX^2)\\
&=X^0\wedge{}^tX^0+X^1\wedge{}^tX^1+X^2\wedge{}^tX^2+X^3\wedge{}^tX^3
-\lambda^2(X^0\wedge{}^tX^0+X^2\wedge {}^tX^2)
\endsplit$$
and
$$\split
{\Om_{\lambda}}^1_0&
=d(\Gamma_1+\omega_1)-\lambda^2d\alpha_1
+\lambda^2(X^0\wedge{}^tX^1+X^2\wedge{}^tX^3)\\
&\quad +(\Gamma_1+\alpha_1-\lambda^2\alpha_1)\wedge\Gamma_0
+\Gamma_0\wedge(\Gamma_1+\alpha_1-\lambda^2\alpha_1)\\
&\quad +(-\Gamma_3+\alpha_3-\lambda^2\alpha_3)\wedge(\Gamma_2+\alpha_2)
+(\Gamma_2-\alpha_2)\wedge(\Gamma_3+\alpha_3-\lambda^2\alpha_3)\,\,\\
&=\Om^1_0-\lambda^2(d\alpha_1-2\alpha_2\wedge\alpha_3)
+\lambda^2(X^0\wedge{}^tX^1+X^2\wedge{}^tX^3)\\
&=X^1\wedge{}^tX^0-X^0\wedge{}^tX^1+X^3\wedge{}^tX^2
-X^2\wedge{}^tX^3+2(X^1\wedge{}^tX^0+{}^tX^2\wedge X^3)\\
&\quad +\lambda^2(X^0\wedge{}^tX^1+X^2\wedge{}^tX^3
-2\,{}^tX^1\wedge X^0-2\,{}^tX^2\wedge X^3)\,\,\\
{\Om_{\lambda}}^2_0&
=d(\Gamma_2+\alpha_2)+\lambda^2(X^3\wedge{}^tX^1-X^1\wedge{}^tX^3)\\
&\quad +(\Gamma_2+\alpha_2)\wedge\Gamma_0
+(\Gamma_3-\alpha_3+\lambda^2\alpha_3)\wedge(\Gamma_1+\alpha_1-\lambda^2\alpha_1)\\
&\quad +\Gamma_0\wedge(\Gamma_2+\alpha_2)
+(-\Gamma_1+\alpha_1-\lambda^2\alpha_1)\wedge(\Gamma_3+\alpha_3-\lambda^2\alpha_3)\\
&=\Omega^2_0+\lambda^2(X^3\wedge{}^tX^1-X^1\wedge{}^tX^3)
+2\lambda^2\alpha_3\wedge\alpha_1\\
&=X^2\wedge{}^tX^0-X^0\wedge{}^tX^2+X^2\wedge{}^tX^3
-X^3\wedge{}^tX^1+2({}^tX^2\wedge X^0+{}^tX^3\wedge X^1)\\
&\quad +\lambda^2(X^3\wedge{}^tX^1-X^1\wedge{}^tX^3)
+4\lambda^2\alpha_3\wedge\alpha_1\,\,.
\endsplit
$$
Similarly, modulo ${\Om'}_{\nu}^{\mu}$, we have : 
$$
\split
{\Om_{\lambda}}^3_0 & = \Om^3_0- \lambda^2 (d\alpha_3-2 \alpha_1 \wedge \alpha_2)
+\lambda^2 (-X^2 \wedge {}^tX^1+X^0 \wedge {}^tX^3)\\
& = X^3\wedge {}^tX^0-X^0\wedge {}^tX^3+X^2\wedge {}^tX^1
-X^1\wedge{}^tX^2+2({}^tX^3\wedge X^0+{}^tX^1\wedge X^2)\\
& \quad +\lambda^2(-X^2\,\wedge {}^tX^1+X^0\wedge{}^tX^3
-2\,{}^tX^3\wedge X^0-2\,{}^tX^1\wedge X^2)\,\,,\\
{\Om_{\lambda}}^2_1 & =\Om^2_1+\lambda^2(d\alpha_3-2\alpha_1\wedge\alpha_2)
+\lambda^2(-X^3\wedge {}^tX^0+X^1\wedge{ }^tX^2)\\
& = X^3\wedge {}^tX^0-X^0\wedge {}^tX^3+X^2\wedge {}^tX^1-X^1\wedge {}^tX^2
+2({}^tX^0\wedge X^3+{}^tX^2\wedge X^1)\\
& \quad +\lambda^2(-X^3\wedge {}^tX^0 + X^1\wedge {}^tX^2
+2\,{}^tX^3\wedge X^0+2\,{}^tX^1\wedge X^2)\,\,,\\
{\Om_{\lambda}}^3_1&=\Om^3_1 + \lambda^2 (X^2\wedge{}^tX^0-X^0\wedge{}^tX^2)
+4\lambda^2 \alpha_3 \wedge \alpha_1\\
&=-X^2 \wedge {}^tX^0+X^0\wedge {}^tX^2-X^1\wedge {}^tX^3+X^3\wedge {}^tX^1
+2({}^tX^2\wedge X^0+{}^tX^3 \wedge X^1)\\
& \quad +\lambda^2 (X^2 \wedge {}^tX^0-X^0 \wedge {}^tX^2)
+4\lambda^2 \alpha_3 \wedge \alpha_1\,\,,\\
{\Om_{\lambda}}^3_2&=\Om^3_2 + \lambda^2 (d\alpha_1-2\alpha_2\wedge\alpha_3)
+ \lambda^2 (X^2\wedge {}^tX^3+X^0 \wedge {}^tX^1)\\
& = X^1\wedge {}^tX^0-X^0 \wedge {}^tX^1+X^3 \wedge {}^tX^2-X^2 \wedge {}^tX^3
+2(X^0 \wedge {}^tX^1+{}^tX^3 \wedge X^2)\\
&\quad +\lambda^2(X^2\wedge {}^tX^3+X^0 \wedge {}^tX^1
+2\,{}^tX^1\wedge X^0+2\,{}^tX^2\wedge X^3)\,\,.
\endsplit
$$
Here, we recall that
$$
\split
\Om^{\mu}_0&=\wt\Om^{\mu}_0+{\Om'}^{\mu}_0\\
&=X^{\mu}\wedge{}^tX^0-X^0\wedge{}^tX^{\mu}+X^{\nu}\wedge{}^tX^{\eta}
-X^{\eta}\wedge{}^tX^{\nu}\\
& \quad +2({}^tX^{\mu}\wedge X^0+{}^tX^{\eta}\wedge X^{\nu})+{\Om'}^{\mu}_0\,\,,\\
\Om^{\eta}_{\nu}&=\wt\Om^{\eta}_{\nu}+{\Om'}^{\eta}_{\nu}\\
&=-X^{\mu}\wedge{}^tX^0+X^0\wedge{}^tX^{\mu}-X^{\nu}\wedge{}^tX^{\eta}
+X^{\eta}\wedge{}^tX^{\nu}\\
&\quad +2({}^tX^{\mu}\wedge X^0+{}^tX^{\eta}\wedge X^{\nu})
+{\Om'}^{\eta}_{\nu}\,\,,
\endsplit$$
$(\mu,\eta,\nu)$ being any cyclic permutation of $(1,2,3)$, where the ``hyper-K\"ahler part" 
${\Om'}^{\mu}_{\nu}$ has no contribution to the Ricci tensor. Therefore we can ignore the 
${\Om'}^{\mu}_{\nu}$-part in the computation of the Ricci tensor. 
\medskip

We are now ready to compute the Ricci tensor of the metric $g_{\lambda}^{\can}$ on $\Cal Z$. 
Note that the dependency on the point under consideration is completely hidden in 
the ``hyper-K\"ahler part" $\Om'$ and $\Om'$ has no contribution to the Ricci tensor. 
Therefore, although we do not assume the homogeneity of $M$, we are able to 
compute the Ricci tensor purely Lie theoretically as if we were working 
on the Riemannian homogeneous space. 
Let $\xi_{-2},\xi_{-1},\xi_0,\xi_1,\xi_2,\xi_3$ be the frame of $\Cal Z$ dual to the coframe
$\alpha_1,\alpha_2,X^0,X^1,X^2,X^3$. 
Then
$$\lambda^{-1}\xi_{-2},\lambda^{-1}\xi_{-1},\xi_0,\xi_1,\xi_2,\xi_3$$
is the frame (orthonormal w.r.to the metric $g_{\lambda}^{\can}$) dual to the coframe 
$$\lambda\alpha_1,\lambda\alpha_2,X^0,X^1,X^2,X^3\,\,.$$
We compute the components of the Ricci tensor using the following formula. 
Let $\{e_i\}_{i=1}^n$ and $(\Om^i_j)_{i,j=1,\dots,n}$ be an orthonormal frame and 
the associated curvature form of an $n$-dimensional Riemannian manifold. 
Then we have
$$\Ric(e_i,e_j)=\sum_{k=1}^ng(R(e_i,e_k)e_k,e_j)=\sum_{k=1}^ng(\Om^j_k(e_i,e_k)e_j,e_j)\,\,.$$
Using this formula, the components of the Ricci tensor $\Ric^{\can}_{\lambda}$ are computed as follows. 
If the expression ${\Om_{\lambda}}^{-2}_k(\lambda^{-1}\xi_{-2},\xi_k)\lambda^{-1}\xi_{-2}$ 
in the following computation means to take the summation over all combination of a fixed vector 
in $\xi_{-2}$ and any vector in $\xi_k$ ($0\leq k\leq 3$), we have
$$
\split
& \quad \Ric^{\can}_{\lambda}(\lambda^{-1}\xi_{-2},\lambda^{-1}\xi_{-2})
=g_{\lambda}^{\can}(\Rm^{\can}_{\lambda}(\lambda^{-1}\xi_{-2},\lambda^{-1}\xi_{-1})\lambda^{-1}\xi_{-1},
\lambda^{-1}\xi_{-2})\\
&\quad +\sum_{k=0}^3g_{\lambda}^{\can}(\Rm^{\can}_{\lambda}(\lambda^{-1}\xi_{-2},\xi_k)\xi_k,
\lambda^{-1}\xi_{-2})\\
&=g_{\lambda}^{\can}({\Om_{\lambda}}^{-2}_{-1}
(\lambda^{-1}\xi_{-2},\lambda^{-1}\xi_{-1})\lambda^{-1}\xi_{-2},\lambda^{-1}\xi_{-2})\\
&\quad +\sum_{k=0}^3g_{\lambda}^{\can}({\Om_{\lambda}}^{-2}_k(\lambda^{-1}\xi_{-2},\xi_k)
\lambda^{-1}\xi_{-2},\lambda^{-1}\xi_{-2})\\
&=\frac{4}{\lambda^2}+4n\lambda^2\,\,.
\endsplit
$$
Indeed, 
$\displaystyle {\Om_{\lambda}}^{-2}_{-1}(\lambda^{-1}\xi_{-2},\lambda^{-1}\xi_{-1})
=4\alpha_1\wedge\alpha_3(\lambda^{-1}\xi_{-2},\lambda^{-1}\xi_{-1})
=\frac{4}{\lambda^2}$ and 
\newline
$\displaystyle \sum_{k=0}^3{\Om_{\lambda}}^{-2}_k(\lambda^{-1}\xi_{-2},\xi_k)
=\lambda^3 X^k\wedge \alpha_1(\xi_k,\lambda^{-1}\xi_{-2})=4n\lambda^2$. 
Similarly, we have
$$
\Ric^{\can}_{\lambda}(\lambda^{-1}\xi_{-1},\lambda^{-1}\xi_{-1})=\frac{4}{\lambda^2}+4n\lambda^2\,\,.
$$
Next, under the similar convention in the expression such as 
${\Om_{\lambda}}^{-1}_k(\lambda^{-1}\xi_{-2},\xi_k)\lambda^{-1}\xi_{-1}$, we have
$$
\split
& \quad \Ric^{\can}_{\lambda}(\lambda^{-1}\xi_{-2},\lambda^{-1}\xi_{-1})
=g_{\lambda}^{\can}(\Rm^{\can}_{\lambda}(\lambda^{-1}\xi_{-2},\lambda^{-1}\xi_{-1})\lambda^{-1}\xi_{-1},
\lambda^{-1}\xi_{-1})\\
&\quad +\sum_{k=0}^3g_{\lambda}^{\can}(\Rm^{\can}_{\lambda}(\lambda^{-1}\xi_{-2},\xi_k)\xi_k,
\lambda^{-1}\xi_{-1})\\
&=g_{\lambda}^{\can}({\Om_{\lambda}}^{-1}_{-1}
(\lambda^{-1}\xi_{-2},\lambda^{-1}\xi_{-1})\lambda^{-1}\xi_{-1},\lambda^{-1}\xi_{-1})\\
&\quad +\sum_{k=0}^3g_{\lambda}^{\can}({\Om_{\lambda}}^{-1}_k(\lambda^{-1}\xi_{-2},\xi_k)
\lambda^{-1}\xi_{-1},\lambda^{-1}\xi_{-1})\\
&=0\,\,.
\endsplit$$
Indeed, 
$\displaystyle \quad \sum_{k=0}^3{\Om_{\lambda}}^{-1}_k(\lambda^{-1}\xi_{-2},\xi_k)
 =\lambda^3X^0\wedge \alpha_3
-(2\lambda-\lambda^2)X^2\wedge\alpha_1(\xi_0,\lambda^{-1}\xi_{-2})
 +\lambda^3 X^1\wedge\alpha_3
-(2\lambda-\lambda^2)X^3\wedge\alpha_1(\xi_1,\lambda^{-1}\xi_{-2})
+\lambda^3 X^2\wedge\alpha_3
+(2\lambda-\lambda^2)X^0\wedge\alpha_1(\xi_2,\lambda^{-1}\xi_{-2})
+\lambda^3 X^1\wedge\alpha_3
-(2\lambda-\lambda^2)X^1\wedge\alpha_1(\xi_3,\lambda^{-1}\xi_{-2})=0$. 

Next, under the similar convention in the expression such as 
$\Rm^{\can}_{\lambda}(\lambda^{-1}\xi_{-2},\xi_k)\xi_k$, we have
$$\split
\Ric^{\can}_{\lambda}(\lambda^{-1}\xi_{-2},\xi_0)&=g_{\lambda}^{\can}(\Rm^{\can}_{\lambda}(\lambda^{-1}\xi_{-2},
\lambda^{-1}\xi_{-1})\lambda^{-1}\xi_{-1},\xi_0)\\
&\quad +\sum_{k=0}^3g_{\lambda}^{\can}(\Rm^{\can}_{\lambda}(\lambda^{-1}\xi_{-2},\xi_k)\xi_k,\xi_0)=0\,\,.
\endsplit$$
Indeed, $\displaystyle {\Om_{\lambda}}^0_{-1}(\lambda^{-1}\xi_{-2},\lambda^{-1}\xi_{-1})\xi_0,\xi_0)
=\lambda^3(X^0\wedge\alpha_3-(2\lambda-\lambda^2)X^2\wedge\alpha_1)
(\lambda^{-1}\xi_{-1},\lambda^{-1}\xi_{-2})\xi_0=0$ and 
$\displaystyle {\Om_{\lambda}}^0_k(\lambda^{-1}\xi_{-2},\xi_k)\xi_0=0$. Similarly, we have
$$
\split
& \quad \Ric^{\can}_{\lambda}(\lambda^{-1}\xi_{-2},\xi_1)=\Ric^{\can}_{\lambda}(\lambda^{-1}\xi_{-2},\xi_2)
=\Ric^{\can}_{\lambda}(\lambda^{-1}\xi_{-2},\xi_3)\\
& =\Ric^{\can}_{\lambda}(\lambda^{-1}\xi_{-1},\xi_0)=\Ric^{\can}_{\lambda}(\lambda^{-1}\xi_{-1},\xi_1)
=\Ric^{\can}_{\lambda}(\lambda^{-1}\xi_{-1},\xi_2)\\
& =\Ric^{\can}_{\lambda}(\lambda^{-1}\xi_{-1},\xi_3)=0\,\,.
\endsplit
$$
We use (2-6) and (2-10) to compute the remaining components of the Ricci tensor 
(under the similar convention in the expression such as 
$\Rm^{\can}_{\lambda}(\xi_0,\xi_0')\xi_0'$ and $\Rm^{\can}_{\lambda}(\xi_0,\xi_k)\xi_k$ 
which is explained more precisely below) : 
$$
\split
\Ric^{\can}_{\lambda}(\xi_0,\xi_0)&
=g_{\lambda}^{\can}(\Rm^{\can}_{\lambda}(\xi_0,\lambda^{-1}\xi_{-2})\lambda^{-1}\xi_{-2},\xi_0)
+g_{\lambda}^{\can}(\Rm^{\can}_{\lambda}(\xi_0,\lambda^{-1}\xi_{-1})\lambda^{-1}\xi_{-1},\xi_0)\\
&\quad +g_{\lambda}^{\can}(\Rm^{\can}_{\lambda}(\xi_0,\xi_0')\xi_0',\xi_0)
+\sum_{k=1}^3g_{\lambda}^{\can}(\Rm^{\can}_{\lambda}(\xi_0,\xi_k)\xi_k,\xi_0)\\
&=4n+8-4\lambda^2\,\,.\endsplit
$$
Indeed, 
$[\Rm^{\can}_{\lambda}(\xi_0,\lambda^{-1}\xi_{-2})\lambda^{-1}\xi_{-2}]_{\xi_0}
={\Om_{\lambda}}^0_{-2}(\xi_0,\lambda^{-1}\xi_{-2})\xi_0=\lambda^3X^0\wedge\alpha_1(\xi_0,\lambda^{-1}\xi_{-2})\xi_0=\lambda^2\xi_0$,  $[\Rm^{\can}_{\lambda}(\xi_0,\lambda^{-1}\xi_{-1})\lambda^{-1}\xi_{-1}]_{\xi_0}
={\Om_{\lambda}}^0_{-1}(\xi_0,\lambda^{-1}\xi_{-1})\xi_0
=\lambda^3X^0\wedge\alpha_3(\xi_0,\lambda^{-1}\xi_{-1})\xi_0=\lambda^2\xi_0$. 
In the following computation, as in the above computation, 
we fix a vector in $\xi_0$ and sum over all vectors in $\xi_0'$ 
(the same as $\xi_0$ but the fixed vector in $\xi_0$ removed) or those in $\xi_k$. 
In other wards,  the pair $(\xi_0,\xi'_0)$ which appears in the following expression, e.g., 
$\Rm^{\can}_{\lambda}(\xi_0,\xi'_0)\xi'_0$ (resp $\Rm^{\can}_{\lambda}(\xi_0,\xi_k)\xi_k$) means to 
take the sum over all combination of the fixed vector in $\xi_0$ and vectors in $\xi_0$ 
except the fixed one (resp. vectors in $\xi_k$). 
We then have 
$$[\Rm^{\can}_{\lambda}(\xi_0,\xi_0')\xi_0']_{\xi_0}
={\Om_{\lambda}}^0_{0'}(\xi_0,\xi_0')\xi_0
=(n-1)\xi_0$$
and furthermore we have
$$
\split
&\quad [\sum_{k=1}^3\Rm^{\can}_{\lambda}(\xi_0,\xi_k)\xi_k]_{\xi_0}
=\sum_{k=1}^3{\Om_{\lambda}}^0_k(\xi_0,\xi_k)\xi_0\\
&=\{X^1\wedge{}^tX^0-X^0\wedge{}^tX^1+X^3\wedge{}^tX^2
-X^2\wedge{}^tX^3+2(X^1\wedge{}^tX^0+{}^tX^2\wedge X^3)\\
&\quad +\lambda^2(X^0\wedge{}^tX^1+X^2\wedge{}^tX^3
-2\,{}^tX^1\wedge X^0-2\,{}^tX^2\wedge X^3)\}(\xi_1,\xi_0)\xi_0\\
&\quad +\{X^2\wedge{}^tX^0-X^0\wedge{}^tX^2+X^2\wedge{}^tX^3
-X^3\wedge{}^tX^1+2({}^tX^2\wedge X^0+{}^tX^3\wedge X^1)\\
&\quad +\lambda^2(X^3\wedge{}^tX^1-X^1\wedge{}^tX^3)
+4\lambda^2\alpha_3\wedge\alpha_1\}(\xi_2,\xi_0)\xi_0\\
& \quad +\{X^3\wedge {}^tX^0-X^0\wedge {}^tX^3+X^2\wedge {}^tX^1
-X^1\wedge{}^tX^2+2({}^tX^3\wedge X^0+{}^tX^1\wedge X^2)\\
& \quad +\lambda^2(-X^2\,\wedge {}^tX^1+X^0\wedge{}^tX^3
-2\,{}^tX^3\wedge X^0-2\,{}^tX^1\wedge X^2)\}(\xi_3,\xi_0)\xi_0\\
&=\{(n+3-3\lambda^2)+(n+3)+(n+3-3\lambda^2)\}\xi_0\\
&=(3n+9-6\lambda^2)\xi_0
\endsplit$$
where $[\cdots]_{\xi_i}$ in the above expression means to take the component 
of the fixed vector 
in $\xi_0$ of $\cdots$. Therefore we have
$$\Ric^{\can}_{\lambda}(\xi_0,\xi_0)=2\lambda^2+(n-1)+3n+9-6\lambda^2=4n+8-4\lambda^2\,\,.$$

Furthermore, we have
$$
\split
\Ric^{\can}_{\lambda}(\xi_0,\xi_1)&=g_{\lambda}^{\can}(\Rm^{\can}_{\lambda}(\xi_0,\lambda^{-1}\xi_{-2})
\lambda^{-1}\xi_{-2},\xi_1)
+g_{\lambda}^{\can}(\Rm^{\can}_{\lambda}(\xi_0,\lambda^{-1}\xi_{-1})\lambda^{-1}\xi_{-1},\xi_1)\\
&\quad +\sum_{k=1}^3g_{\lambda}^{\can}(\Rm^{\can}_{\lambda}(\xi_0,\xi_k)\xi_k,\xi_1)=0\,\,.
\endsplit$$
Indeed, $[\Rm^{\can}_{\lambda}(\xi_0,\lambda^{-1}\xi_{-2})\xi_{-2}]_{\xi_1}
={\Om_{\lambda}}^1_{-2}(\xi_0,\lambda^{-1}\xi_{-2})\xi_1
=(\lambda^3X^1\wedge\alpha_1+(2\lambda-\lambda^2)X^3\wedge\alpha_3)
(\xi_0,\lambda^{-1}\xi_{-2})\xi_1=0$, 
$[\Rm^{\can}_{\lambda}(\xi_0,\lambda^{-1}\xi_{-2})\xi_{-2}]_{\xi_1}
={\Om_{\lambda}}^1_{-1}(\xi_0,\lambda^{-1}\xi_{-1})\xi_1
=(\lambda^3X^1\wedge\alpha_3-(2\lambda-\lambda^2)X^3\wedge\alpha_1)
(\xi_0,\lambda^{-1}\xi_{-1})\xi_1=0$ 
and $\Rm^{\can}_{\lambda}(\xi_0,\xi_k)\xi_k={\Om_{\lambda}}^1_k(\xi_0,\xi_k)\xi_1=0$. 

Similarly we have
$$
\split
& \Ric^{\can}_{\lambda}(\xi_1,\xi_1)=\Ric^{\can}_{\lambda}(\xi_2,\xi_2)
=\Ric^{\can}_{\lambda}(\xi_3,\xi_3)=4n+8-4\lambda^2\,\,,\\
& \Ric^{\can}_{\lambda}(\xi_0,\xi_2)=\Ric^{\can}_{\lambda}(\xi_0,\xi_3)=\Ric^{\can}_{\lambda}(\xi_1.\xi_2)
=\Ric^{\can}_{\lambda}(\xi_1,\xi_3)=\Ric^{\can}_{\lambda}(\xi_2,\xi_3)=0\,\,.
\endsplit
$$
To sum up, we have

\proclaim{Proposition 2.4} The Ricci tensor of the metric 
$$g^{\can}_{\lambda}=\lambda^2(\alpha_1^2+\alpha_3^2)+\sum_{i=0}^3{}^tX^i\cdot X^i$$ 
on the twistor space $\Cal Z$ is given by the formula
$$
\Ric^{\can}_{\lambda}=(4\lambda^{-2}+4n\lambda^2)\lambda^2(\alpha_1^2+\alpha_3^2)
+(4n+8-4\lambda^2)
({}^tX^0\cdot X^0+{}^tX^1\cdot X^1+{}^tX^2\cdot X^2+{}^tX^3\cdot X^3)\,\,.
\tag2-15
$$
In particular, the canonical deformation metric $g_{\lambda}^{\can}$ is Einstein if and only if 
$\lambda^2=1$ or $\lambda^2=\frac1{n+1}$. If $\lambda^2=1$ then $g_{\lambda}^{\can}$ is 
K\"ahler-Einstein and if $\lambda^2=\frac1{n+1}$ then $g^{\can}_{\lambda}$ is 
Einstein, Hermitian but not K\"ahler. 
\endproclaim

\bigskip
\noindent
{\bf \S3. Z-metrics.}
\par
\bigskip

In the following discussion, we introduce new metrics (called Z-metrics 
denoted by $g^{\ZZ}_{\lambda}$) on the 
twistor space of positive quaternion K\"ahler manifolds and compute 
the curvature form. 
We use the moving frame computation. The scalings of the various 
standard metrics are hidden in the computation. 
To avoid confusion, we fix our scaling convention in the following way : 
\medskip
\noindent
$\bullet$ We fix the scale of the invariant metric of $\H\PP^n$ so that the sectional 
curvatures range in the interval $[1,4]$, i.e., $\Ric(g_{\H\PP^n})=4(n+2)g_{\H\PP^n}$ 
and so $\text{\rm Scal}(g_{\H\PP^n})=16n(n+2)$. Namely we set 
$$\wt S=16n(n+2)$$
from here on. 
\medskip
\noindent
$\bullet$ We fix the scale of the Fubini-Study metric of the $\PP^1$-fiber of the twistor 
fibration and other cases so that the Gaussian curvature is identically $4$. 
\medskip
\noindent
$\bullet$ We fix the scale of the invariant metric of the $\Sp(1)$-fiber of the extended 
twistor fibration $\wcalZ\rightarrow M$ so that the sectional curvature is identically $1$. 
\medskip

We will consider the following partial scaling similar to the case of canonical 
deformation metrics in \S2. 
\medskip
\noindent
$\bullet$ We normalize the base quaternion K\"ahler metric $g$ so that $S=\wt S$ holds 
and scale the 
$(\alpha_1-\alpha_1(\xi_1)X^1-\alpha_1(\xi_3)X^3)^2
+(\alpha_3-\alpha_3(\xi_1)X^1-\alpha_3(\xi_3)X^3)^2$-part (ignoring ``invisible'' 
part) by the parameter 
$\lambda^2$ (so that the ``curvature becomes $\lambda^{-2}$-times the original one in this 
direction"). 
It turns out that the Z-metric $g^{\ZZ}_{\lambda}$ is Einstein (but not K\"ahler) on $\Cal Z$ 
if and only if $\lambda=\frac1{n+2}$. 
\medskip

Next we proceed to a construction of a new family of Riemannian metrics
$$\R_+\cdot \Cal F^{\ZZ}=\{\rho g_{\lambda}^{\ZZ}\}_{\lambda>0,\rho>0}$$
(the family of Z-metrics) on the twistor space $\Cal Z$. 
Our strategy is to modify the construction of the canonical deformation metrics 
so that kill the $-\lambda^2$-term in before the $\sum_{i=0}^3{}^tX^i\cdot X^i$-term in the 
formula of $\Ric^{\can}_{\lambda}$ in Proposition 2.4. Indeed, the $-\lambda^2$-term 
in question constitutes the basic reason why the family of canonical deformation metrics 
is {\it not} a Ricci flow unstable cell. The simplest example of a Ricci flow unstable cell 
arises in the Ricci flow solutions on the product space $S^n\times S^n$ where initial metric 
an independently scaled  constant curvature metrics (if the initial 
metric has slightly different curvature, then the $S^n$-factor with larger curvature 
will extinct earlier and before its extinction the metric of the total product space 
is ``far'' from being Einstein !). 
Therefore we try to modify the construction 
of the canonical deformation metric so that the resulting metric becomes ``closer'' to 
the ``product'' metric (of $(M,g)$ and $(\PP^1,g_{\FS})$, in a ``weak'' sense). 
The construction should be ``canonically'' based on the definition of the 
twistor space in a way somewhat different from the family of canonical deformation metrics. 
So we begin the construction of $\Cal F^{\ZZ}$ with a useful interpretation of the twistor space. 
We identify $\H$ with $\C^2$ by writing $x_0+ix_1+jx_2+kx_3=x_0+jx_2+i(x_1+jx_3)$. 
Let $x+jy\in\H$ and $u+jv\in\H$. The manipulation 
$$(x+jy)(u+jv)=(xu-\ov yv)+j(yu+\ov xv)
=\pmatrix x&-\ov y\\ y&\ov x\endpmatrix\pmatrix u\\ v\endpmatrix$$
implies that the right action of the unit quaternions is identified with the standard holomorphic 
action of $\SU(2)$ to $\C^2$ and the left action of the unit quaternions is not holomorphic. 
The right action of the unit quaternion $u+jv$ is holomorphic if and only if $v=0$, i.e., 
$u+jv$ is the unit complex number. The local decomposition 
$\SO(4)=\Sp(1)_{\text{\rm l}}\times \Sp(1)_{\text{\rm r}}$ 
implies that the standard action of $\SO(4)$ on $\R^4$ decomposes into the left and right actions 
of $\Sp(1)$ (denoted by $\Sp(1)_{\text{\rm l}}$ and $\Sp(1)_{\text{\rm r}}$). 
Therefore the space of the orthogonal complex structures 
of $\H$ is identified with $\SO(4)/\Un(2)\cong \Sp(1)_{\text{\rm r}}/\Un(1)=\PP^1$. 
On the other hand, the standard action of $\Sp(1)_{\text{\rm l}}=\SU(2)$ on $\C^2$ implies that 
the space of all complex lines in $\C^2$ is identified with $\Sp(1)_{\text{\rm l}}/\Un(1)
=\SU(2)/\text{\rm S}(\Un(1)\times \Un(1))=\PP^1$. The identification of  $\PP^1$ is with 
the space of all orthogonal complex structures on $\H$ is given by 
$$\split
& \text{\rm orthogonal complex structure $J$ on $\H$}\\
&\quad  \leftrightarrow 
\text{\rm $\Un(1)$-subgroup of $\Sp(1)_{\text{\rm r}}$ which acts on $\C^2_J$ 
holomorphically}
\endsplit$$
and the identification of $\PP^1$ is with 
the space of all complex lines in $\C^2$ is given by
$$
\split
& \text{\rm complex line $L$ in $\C^2$}\\
&\quad  \leftrightarrow 
\text{\rm $\Un(1)$-subgroup of $\Sp(1)_{\text{\rm l}}$ which fixes $[L]\in \PP^1$}\,\,.
\endsplit$$
Identifying $\Sp(1)_{\text{\rm l}}$ and $\Sp(1)_{\text{\rm r}}$, we can 
transfer a $\Un(1)$-subgroup in $\Sp(1)_l$ to a $\Un(1)$-subgroup 
in $\Sp(1)_r$.  A $\Un(1)$-subgroup in $\Sp(1)_r$ determines an axis of rotation 
of $\PP^1(\C)$ and therefore two lines in $\C^2$ fixed by the $\Un(1)$-action. 
As the twistor line $\PP^1(\C)$ is simply connected, we can choose one of them 
globally and therefore we have the following canonical correspondence
$$
\split
& \text{\rm orthogonal complex structure $J$ on $\H$}\\
& \leftrightarrow 
\text{\rm a complex line $L_J$ of $\C^2_J$}\,\,.
\endsplit
$$
Similarly, $\Sp(n)\Sp(1)/\Sp(n)\Sp(1)\cap\Un(2n)=\PP^1$ is 
identified with the space of all orthogonal complex structures of $\H^n$ 
and also to the space of all (simultaneous choice of) complex lines 
in each $\H$-line in $\H^n$. Therefore, applying the above correspondence to 
each $\H$-lines in $\H^n$, we have the following canonical correspondence
$$
\split
& \text{\rm orthogonal complex structure $J$ on $\H^n$}\\
&\quad \leftrightarrow 
\text{\rm a real $(4n-2)$-dimensional complex subspace $D'_J$ of $\C^{2n}_J$}\\
&\quad \text{\rm which cuts out a complex line from each $\H$-line in $\H^n$}\\
& \quad \leftrightarrow \text{\rm a complex line $L_J=(D'_J)^{\perp}$ of $\C^{2n}_J$}\,\,.
\endsplit
$$
Consider the standard twistor fibration 
$$\pi : (\PP^{2n+1}(\C),\PP^{2n-1}(\C)) \rightarrow (\PP^n(\H),\PP^{n-1}(\H))\,\,.$$
Pick a point $m\in\PP^n(\H)-\PP^{n-1}(\H)$. Then the fiber $\PP^1_m=\pi^{-1}(m)$ 
is identified with the moduli space of the orthogonal complex structures of 
$T_m\PP^n(\H)\cong \H^n$. Pick a point $z\in\PP^1_m$ and let $J(z)$ be the 
corresponding orthogonal complex structure. 
The hyperplane $D_z$ of $\PP^{2n+1}(\C)$ generated by 
$\PP^{2n-1}(\C)=\pi^{-1}(\PP^{n-1}(\H))$ 
(inverse image of the hyperplane at infinity) and $z\in\PP^1_m$ 
gives an identification of the affine part $\H^n\cong \PP^n(\H)-\PP^{n-1}(\H)$ 
with $\C^{2n}_{J(z)}$, i.e., $\H^n$ with the orthogonal complex structure $J(z)$. 
Assume that the point $m$ under consideration is the origin of $\H^n$. 
Then the tangent space of $D''_z$ coincides with the horizontal subspace $H_z$ 
at $z\in\PP^1_m$. 
At each point $z\in\PP^1_m$ we have a 2-dimensional subspace 
$L_{J(z)}:=(D'_{J(z)})^{\perp}$ of $D_z=H_z$ which is a complex line w.r.to the 
orthogonal complex structure $J(z)$. 
If we consider $\H$-lines 
$L\cong \PP^1(\H)$ in $\PP^n(\H)$ passing through the origin $m\in\H^n=\PP^n(\H)-\PP^{n-1}(\H)$ 
and the sub twistor fibrations $\PP^3(\C)\cong \pi^{-1}(L) \rightarrow L\cong \PP^1(\H)$ 
over these $\H$-lines, the collection of $J(z)$-complex lines from each $\H$-line 
in $\H^n$ constitutes a real $(4n-2)$-dimensional subspace $D'_z$ which is 
a complex subspace of $\C^{2n}_{J(z)}$. 
From this consideration, we see that 
there exists a column $n$ vector $\xi_0\in D_z'$ uniquely modulo $\SO(2)$-rotation 
so that $\{\xi_0,J(z)\xi_0\}_{\text{\rm span}}$ is a real $(2n)$-dimensional complex 
subspace in $D_z'$. 
We observe that there is a ``canonical'' way constructing a real $(4n)$-dimensional 
non-horizontal subspace $D_z$ of $T_z\PP^{2n+1}(\C)$ which is complex w.r.to $J(z)$ 
by using the orthogonal complement $L_{J(z)}=(D'_{J(z)})^{\perp}$. In fact, $D_z$ 
is constructed as a direct sum of $D'_z$ and a complex line $L_{J(z)}'$ 
in the linear subspace ($\cong \C^2$) of $T_z\PP^{2n+1}(\C)$ spanned 
by $T_z\PP^1_m$ and $L_{J(z)}$. 
We would like to ``canonically'' construct $L_{J(z)}'$ 
by choosing unit vectors $u \in T_z\PP^1_m$ and $v \in L_{J(z)}$ and putting 
$L_{J(z)}':=\C\cdot (u+v)$. Here, $u$ and $v$ are chosen from the tangent bundle 
$\Cal O(2)$ of $\PP^1$. Therefore defining such $L_{J(z)}'$ is equivalent 
to defining a non-trivial embedding 
$\Cal O(2) \rightarrow \Cal O(2)\oplus \Cal O(2)$, i.e., defining a graph of 
a non-zero section of $\text{\rm Hom}(\Cal O(2),\Cal O(2)=\Cal O(-2)\otimes \Cal O(2)=\Cal O$. 
Therefore $L_{J(z)}'$ is defined  modulo unit complex numbers, i.e., a choice of 
a unit non-zero section of $\Cal O$. So, we choose $1\in H^0(\PP^1,\Cal O)$ 
to define $L_{J(z)}'$ canonically. 
The $S^1$-fiber in the extended twistor space 
$\wcalZ=\frac{\Sp(n+1)/\Z_2}{\Sp(n)\Sp(1)\cap\SU(2n)}$ sitting over a point $z\in\Cal Z$ of 
a twistor line equip $L_{J(z)}$ (also $L_{J(z)}'$) with additional information, i.e., 
an oriented orthonormal frame. 
Indeed, the $S^1$-bundle 
$\wcalZ \rightarrow \frac{\Sp(n+1)/\Z_2}{\Sp(n)\Sp(1)\cap\Un(2n)}$ induces 
the Hopf fibration $S^3 \rightarrow S^2$ over each twistor line in $\Cal Z$ and 
therefore the $S^1$-fiber over $z$ corresponds to 
the rotation by unit complex numbers of an oriented orthonormal frame of 
$L_{J(z)}$ (also $L_{J(z)}'$). 
Therefore we have a well-defined
$$L_{J(z)}':=\C\cdot (u+v)$$
which is ``canonical'' from Riemannian view point (indeed, this is also ``canonical'' from 
the construction of Z-metrics which is explained in the following). 
We get a ``canonical'' non-horizontal subspace $D_z$ as an orthogonal direct sum
$$
D_z:=D_z'\oplus L_{J(z)}'\,\,.
$$

We proceed to a general case. Let $(M^{4n},g)$ ($n\geq 2$) 
be a positive quaternion K\"ahler manifold. 
Let $\Cal H$ be the horizontal distribution of the twistor fibration, i.e., the $(4n)$-dimensional 
distribution consisting of the horizontal subspaces $\Cal H_z$ at $z\in\Cal Z$. 
We consider the $\G(4n-2,4n)$-bundle ($\G(p,q)$  being the Grassmannian 
of $p$-dimensional subspaces in $\R^q$) associated to the 
holonomy reduction $\Cal P\rightarrow M$ of the oriented orthonormal frame bundle. 
As the structure group of the oriented orthonormal frame bundle 
reduces from $\SO(4n)$ to $\Sp(n)\Sp(1)$, the associated 
$\G(4n-2,4n)$-bundle also reduces to a smaller bundle with fiber 
$\Sp(n)\Sp(1)/\Sp(n)\Sp(1)\cap \text{\rm S}(O(2)\times O(4n-2)) 
\cong \Sp(n)\Sp(1)/\Sp(n)\Sp(1)\cap \Un(2n) \cong \PP^1$ 
which turns out to be isomorphic to the twistor fibration $\Cal Z \rightarrow M$. 
Let $m\in M$ and $z\in\PP^1_m=\pi^{-1}(m)$, $\pi\,:\,\Cal Z\rightarrow M$ 
being the twistor fibration. The holonomy reduction of the associated $\G(4n-2,4n)$-bundle 
corresponds to the association
$$\PP^1_m \ni z \longmapsto \Cal D_z'\in \G(4n-2,H_z)$$
where $\Cal D_z'$ is the same as in the model case 
$\PP^{2n+1}(\C)\rightarrow \PP^n(\H)$, i.e., $\Cal D_z'$ is a real $(4n-2)$-dimensional 
subspace of $H_z\cong \C^{2n}_{J(z)}$ which is complex w.r.to $J(z)$ and 
specifies the complex line from each $\H$-lines in $H_z$ ($J(z)$ is the orthogonal 
complex structure of $T_mM$ determined by $z\in\PP^1_m$). 
Let $\Cal L_{J(z)}$ be the orthogonal complement in $\Cal H_z$ 
of $\Cal D_z'\subset \Cal H_z$. Then $L_{J(z)}$ is a complex line in $\C^{2n}_{J(z)}$ 
($\Cal H_z$ equipped with the orthogonal complex structure $J(z)$). 
We consider the subspace ($\cong\C^2$) of $T_z\Cal Z$ spanned by 
$T_z\PP^1_m\cong\C$ and $\Cal L_{J(z)}(\cong\C)\subset \Cal H_z$. 
If we choose unit vectors $u\in T_z\PP^1_m$ and $v\in\Cal L_{J(z)}$ 
in a suitable way, then
$$\Cal L_{J(z)}':=\C\cdot (u+v)$$
is well-defined as before. Indeed, we define the extended twistor space by 
$$\wcalZ:=\Cal P/\Sp(n)\Sp(1)\cap\SU(n)\,\,.$$
Then the ambiguity by multiplication by unit complex numbers in the choice of $u$ 
and $v$ is interpreted as the rotation along the $S^1$-fiber of $\wcalZ\rightarrow\Cal Z$. 
Therefore we can work on $\wcalZ$ without ambiguity and 
the resulting vector $u+v$ on $\Cal Z$ is defined modulo multiplication by unit 
complex numbers. This means that the complex line $\C\cdot(u+v)$ is well-defined. 
The $(4n)$-dimensional non-horizontal distribution $\Cal D:=\{\Cal D_z\}_{z\in\Cal Z}$ 
on $\Cal Z$ is now defined by
$$\Cal D_z:=\Cal D'_z\oplus \Cal L_{J(z)}'\,\,.$$
We choose a column $n$ vector $\xi_0$ of $\Cal D'_z$ (representing an $n$-dimensional 
subspace of $T_mM$). Then $\{\xi_0,J(z)\xi_0\}_{\text{\rm span}}$ is 
a $J(z)$-invariant $(2n)$-dimensional subspace of $\Cal D'_z$. 
In this situation, $J(z)$ is the unique orthogonal complex structure 
of $T_mM$ with the above property (i.e., any (2n)-dimensional subspace of $T_mM$ 
spanned by $\{X_1,qX_1\}$ where $X_1\in D'_z$ and $q$ is 
any orthogonal complex structure of $T_mM$ not equal to $\pm J(z)$ 
is not contained in $D_z$). 
The correspondence $z\mapsto \Cal D_z'$ defines 
a $(4n-2)$-dimensional horizontal distribution $\Cal D'$ 
on the twistor space $\Cal Z$. 
Each point $z\in Z$ represents an orthogonal complex structure $J(z)$ ($J(z)=J$, say) 
and a $(4n-2)$-dimensional subspace $D'_z \subset T_mM$ ($z$ lies over $m\in M$) 
which contains the $(2n)$-dimensional $J$-complex subspace spanned 
by $\{\xi_0,J\xi_0\}$. 
Let $\Cal S\subset \text{\rm End}^{\text{\rm skew}}(TM)$ be the 3-dimensional 
sub-bundle over $M$ which defines the quaternion K\"ahler structure of $(M,g)$. 
Let $\{I,J,K\}$ ($J=J(z)$ under consideration) be an orthonormal basis of $\Cal S_m$. 
Then $\{\xi_0,I\xi_0,J\xi_0,K\xi_0\}$ is an orthonormal basis of $T_mM$. 
We extend $\xi_0$ to a vector field germ at $m \in M$ so that the contribution 
from the Levi-Civita connection to $\nabla \xi_0$ satisfies certain condition 
to be specified later. 
If we further extend the triple $\{I,J,K\}$ as orthonormal frame 
germ of the bundle $\Cal S\rightarrow M$, we get an orthonormal frame germ 
$\{\xi_0,I\xi_0,J\xi_0,K\xi_0\}$. 
First we extend $J$ to a unit section germ of the bundle $\Cal S \rightarrow M$ 
so that $\nabla J$ vanishes in $\xi_0$- and $J\xi_0$-directions 
at $m$ (this is the condition on the $\sp(1)$-part of the Levi-Civita connection). 
Then we use the distribution $\Cal D$ 
which we have just constructed above to extend $I$ and $K$. 
Namely we extend $I$ and $K$ 
in the following way. The section germ $\{\xi_0,I\xi_0,J\xi_0,K\xi_0\}$ of 
the holonomy reduction 
$\Cal P\rightarrow M$ of the oriented orthonormal frame bundle composed 
with the projection 
$\Cal P\rightarrow \Cal Z$ defines a section germ $\sigma$ of 
the twistor fibration $\pi\,:\,\Cal Z\rightarrow M$ at $m$ passing through 
$z\in\PP^1_m$. 
Since $\sigma$ is a map germ $(M,m)\rightarrow (\Cal Z,z)$ 
which is linearly non degenerate 
at $m$, we can define the $(4n)$-dimensional subspace 
$(d\sigma)_m(T_mM)\subset T_z\Cal Z$. 
We choose the extension of $I$ and $K$ so that the relation
$$
(d\sigma)_m(T_mM)=\Cal D_z\tag3-1
$$
holds. 
Since local orthonormal frame fields of $M$ one to one correspond to local sections 
of the bundle $\Cal P\rightarrow M$, the requirement (3-1) is realized by certain 
extensions of $I$ and $K$ such that the triple $\{I,J,K\}$ constitutes an orthonormal 
frame germ of $\Cal S \rightarrow M$ at $m$. Now we describe the procedure 
precisely. 
Given a section germ $J$ satisfying $J^2=-1$ 
of the bundle $\Cal S\rightarrow M$, the extension of $I$ and $K$ so that 
the triple $\{I,J,K\}$ constitutes an oriented orthonormal frame germ of $\Cal S$ 
(which necessarily satisfies the quaternion relations) is unique modulo 
(non constant) $\SO(2)$-rotation. 
This way we get an orthonormal frame germ 
$\{\xi_0,I\xi_0,J\xi_0,K\xi_0\}$. This triple defines a section germ 
of the holonomy reduction $\Cal P \rightarrow M$ of the oriented orthonormal 
frame bundle at $m$. The composition with the projection 
$\Cal P\rightarrow \Cal Z$ defines a section germ $\sigma$ of the twistor bundle 
$\Cal Z \rightarrow M$. The section germ $\sigma$ defines a section germ 
(denoted by the same symbol $\sigma$) 
from (the image of $\sigma$ in) $\Cal Z$ to $\Cal P$
and defines the pull-back of differential forms on $\Cal P$ 
to those on the twistor space $\Cal Z$ (taking values only on tangent vectors of 
the image of $\sigma$ in $\Cal Z$). 
The requirement that the section germ $\sigma$ satisfies the 
condition (3-1) is equivalent to requiring the following conditions: 
(i) the section germ $\xi_0$ at $m$ satisfies 
the condition that the contribution to $\nabla\xi_0$ at $m$ 
from the components $\Gamma_0$, $\Gamma_2$, $\alpha_1$ and $\alpha_3$ of 
the Levi-Civita connection of $g$ vanish (are ``invisible'') 
in $\xi_0$- and $J\xi_0$-directions\footnote{\,\,Although 
performing the projection $\Cal P\rightarrow \Cal Z$ significantly decreases 
the information on the $\Gamma$-part, the vanishing of $\Gamma_0$ 
and $\Gamma_2$ at $m$ in the $\xi_0$- and $J\xi_0$-directions 
is still a very important information on the choice of $\xi_0$, because 
this reflects the geometry behind the construction of the distribution $\Cal D'$.}, 
(ii) the unit section germ $J$ of the bundle $\Cal S\rightarrow M$ at $m$ satisfies 
the condition that the contribution to $\nabla J$ at $m$ vanishes in 
$\xi_0$- and $J\xi_0$-directions, (iii) the $\sp(1)$-part of the Levi-Civita connection 
in $\nabla I$ and $\nabla K$ at $m$ vanishes in $\xi_0$- and $J\xi_0$-directions,  
(iv) the original $\sp(1)$-part of the connection form defined on $\Cal Z$ 
should be replaced by the pull-back of the $\sp(1)$-part of the connection form 
defined on $\Cal P$ via the map germ $\Cal Z \rightarrow \Cal P$ 
constructed from the section germ (composed with the projection 
$\Cal P\rightarrow \Cal Z$) $\sigma : (M,m) \rightarrow (\Cal Z,z)$ 
satisfying (3-1). However, what we get by this procedure coincides with 
the original one. 

As we will work on the twistor space $\Cal Z$, we must define $\xi_i$'s as 
tangent vectors at $z$ of $\Cal Z$ and extend them as a vector field germs 
on $\Cal Z$ (rather than vector fields germs  on $M$). 
We recall that $\xi_i$'s, as tangent vectors at $z\in\Cal Z$,  
are defined as the image under the differential $d\sigma$ of the orthonormal section 
germ $\sigma$ determined by $\{\xi_0,I\xi_0,J\xi_0,K\xi_0\}$. 
Therefore if we extend $\xi_0$ as a section of the distribution $\{D_z'\}_{z\in\PP_m^1}$ 
in the $\PP^1_m$-fiber direction, the extended object defines a desired vector field 
germ $\xi_0$ on $\Cal Z$. 
It follows from the expression (2-2) that the derivation formula satisfied by the extended 
$\xi_0$ is expressed as
$$
\nabla\xi_0=\Gamma_0\otimes\xi_0+(\Gamma_1+\alpha_1)\otimes \xi_1
+(\Gamma_2+\alpha_2)\otimes\xi_2+(\Gamma_3+\alpha_3)\otimes\xi_3
\tag3-2
$$
where $\nabla$ is the connection of the distribution $\Cal D$ on $\Cal Z$ 
obtained by pulling back the Levi-Civita connection of $g$ on $M$ 
via the local map $\sigma$ from $\Cal Z$ to $\Cal P$ (in the $\PP^1$-fiber direction 
the formula (3-2) describes the rotation by the $\sp(1)$-part $\alpha_1$ and $\alpha_3$ 
of the Levi-Civita connection of $g$). 
Moreover, all 1-forms appearing as coefficients in the above formula are obtained by pulling 
back the connection form of the Levi-Civita connection of $g$. The connection form is 
originally defined on $\Cal P$ and pulled back via the section germ $\sigma$ to $\Cal Z$, 
where $\sigma$ is regarded as a map from $\Cal Z$ to $\Cal P$ by identifying $TM$ 
and $d\sigma(TM)$ and the extension of $\xi_0$ in the fiber $\PP^1_m$ direction. 
Once we extend $\xi_0$ to a vector field germ on $\Cal Z$ at $z$, we automatically 
get the quadruple $\{\xi_0,I\xi_0,J\xi_0,K\xi_0\}$ which defines a map germ from $\Cal Z$ to 
$\Cal P$ at $z$. 
We note the followings : (i) the components $\Gamma_0$ and $\Gamma_2$ from the 
$\Gamma$-part (the $\Sp(n)$-component of the Levi-Civita 
connection of $g$) are ``invisible'' in $\xi_0$- and $J\xi_0$-directions 
(but ``visible'' in the $\xi_1$- and $\xi_3$- directions) at $m$ and 
no condition is imposed on  the $\Gamma$-part, 
(ii) the components $\alpha_1$ and $\alpha_3$ from the $\sp(1)$-part of 
the Levi-Civita connection of $g$ are  ``invisible'' in $\xi_0$- and $J\xi_0$-directions 
(but ``visible'' in the $\xi_1$- and $\xi_3$- directions). In addition to these 
conditions, we note that 
(iii) the component $\alpha_2$ is ``invisible'' at $m$ in every direction 
(this is a consequence of the definitions of the twistor space 
and the triple $\{\alpha_1,\alpha_2,\alpha_3\}$). 

The reason why no condition is imposed on the $\Gamma$-part at this stage 
is that the information on the $\Gamma$-part of the Levi-Civita connection form 
of $g$ is lost after performing the projection $\Cal P \rightarrow \Cal Z$. 
This means that, although the equations $\alpha_1=\alpha_3=0$ defines 
the horizontal subspace in the twistor space $\Cal Z$, the $\Gamma$-part of the 
connection form which appears in $\nabla\xi_0$ depends on 
the original local section of the orthonormal frame bundle $\Cal P\rightarrow M$ 
by which the $\Gamma$-part is pulled back to $\Cal Z$ from $\Cal P$ 
and therefore $\Gamma_i$'s do not necessarily vanish at horizontal vectors 
in the twistor fibration $\Cal Z\rightarrow M$. Later, we need to impose a 
certain condition on the $\Gamma$-part in question,  in order to construct an Einstein metric 
together with Ricci flow unstable cell in \S4 and effectively apply 
Bando-Shi's gradient estimate in \S5. The condition we need is the following : 
$$
\Gamma_1\,\,\text{\rm and}\,\, \Gamma_3\,\,\text{\rm  in (3-2) vanish at}\,\,m. \tag3-2$'$
$$
This is the condition on the $\Gamma$-part of $\nabla\xi_0$ where $\xi_0$ is 
the extension to a vector field germ on $M$ at $m$, which can be certainly realized. 
\medskip
 
To sum up, 
\newline\noindent
(a) We extend a horizontal vector $\xi_0$ chosen at $z\in\Cal Z$ to a 
vector field germ at $z\in\Cal Z$ in the following way : 
\newline\noindent
(a1) We first extend $\xi_0$ to a vector field germ on $M$ at $m$ so that 
the $\sp(1)$-part of the connection form vanishes at $m$ in the $\xi_0$- and 
$J\xi_0$-directions (i.e., $\alpha_1(\xi_i)=\alpha_3(\xi_i)=0$ ($i=0,2$). Here, we automatically 
$\alpha_2=0$ at $m$ and we impose no condition on the $\sp(n)$-part. 
\newline\noindent
(a2) Then we extend $\xi_0$ as a section of the bundle $\{\Cal D_z'\}_{z\in\PP_m^1}$ in 
the fiber $\PP^1_m$-direction. 
\newline\noindent
(b) We extend $J$ (the orthogonal complex structure of $T_mM$ at $z\in\PP^1_m$) 
to a section $\text{\rm End}^{\text{\rm skew}}(\Cal D)|_{\PP_m^1}$ in the following way : 
\newline\noindent
(b1) We first extend $J$ to a section germ of $\Cal S \rightarrow M$ so that $\nabla J=0$ 
holds at $m$. 
\newline\noindent
(b2) Then we extend $J$ along the $\PP^1_m$ fiber tautologically, i.e., so that 
$\nabla J=2\alpha_3\otimes I-2\alpha_1\otimes K$ holds. 
\newline\noindent
(b3) Then we extend $I$ and $K$ along all directions of $\Cal Z$ so that 
$\nabla I=-2\alpha_3\otimes J+2\alpha_2\otimes K$ 
and $\nabla K=-2\alpha_2\otimes I+2\alpha_1\otimes J$ hold, where 
$\alpha_1$ and $\alpha_3$ are the $\sp(1)$-part of the connection form, 
which is identified with the pulled back to $\Cal Z$ via the map germ 
$\Cal Z\rightarrow \Cal P$ constructed from the section germ (composed 
with the projection $\Cal P\rightarrow \Cal Z$) $\sigma:(M,m)\rightarrow 
(\Cal Z,z)$ satisfying (3-1). 
\newline\noindent
(b4) In \S5, we will use the condition (3-2$'$). 
\medskip

We would  like to use the above constructed $(4n)$-dimensional 
{\it non-horizontal} distribution $\Cal D$ on $\Cal Z$ to construct the family $\Cal F^{\ZZ}$ 
of Z-metrics on $\Cal Z$. This attempt works most symmetric way on the canonical 
$\SO(2)$-extension $\wcalZ$ of $\Cal Z$. Define the extended twistor space $\wcalZ$ by
$$\wcalZ=\Cal P/\Sp(n)\Sp(1)\cap\SU(2n)\,\,.$$
The meaning of the $\SO(2)$-bundle $\wcalZ\rightarrow \Cal Z$ is the following. 
This $\SO(2)$-bundle structure defines a Hermitian holomorphic negative line bundle 
on the complex manifold $\Cal Z$ with its curvature form proportional to 
the K\"ahler-Einstein form on $\Cal Z$ (see a statement just after the formula (2-8) in 
old version)). 

We construct a family $\wcalF^{\ZZ}$ of metrics on $\wcalZ$ 
and the family $\Cal F^{\ZZ}$ of Z-metrics on $\Cal Z$ is defined uniquely 
so that the projection $\wcalZ\rightarrow \Cal Z$ is a Riemannian submersion. 

We start the construction of $\wcalF^{\ZZ}$ by finding a basis of the 3-dimensional 
space spanned by $\alpha_i$'s ($i=1,2,3$) modified by $X^j$'s ($j=0,1,2,3$) which 
annihilate the natural lift (defined later) to $\wcalZ$ of $\Cal D_z$. 

The $\sp(1)$-component
$$\{\alpha_1,\alpha_2,\alpha_3\}$$
of the connection form of the Levi-Civita connection of the original quaternion K\"ahler metric 
$g$ fits into the formula
$$\nabla \pmatrix I&J&K\endpmatrix=\pmatrix I&J&K\endpmatrix \otimes 
\pmatrix 0 & 2\alpha_3 & -2 \alpha_2\\
-2\alpha_3 & 0 & 2\alpha_1\\
2\alpha_2 & -2\alpha_1&0\endpmatrix\,\,.$$
This formula makes sense as the derivation formula of the oriented orthonormal section 
germ of the bundle $\Cal S(\Cal D)$ of $\text{\rm End}^{\text{\rm skew}}(\Cal D)$ 
w.r.to the connection induced from the connection just as in (3-1), where the sub-bundle 
$\Cal S(\Cal D)\subset \text{\rm End}^{\text{\rm skew}}$ is induced from 
$\Cal S \subset \text{\rm End}^{\text{\rm skew}}$ just in the same sense as (3-1). 
From here on, the connection $\nabla$ which will appear in all derivation formulae 
should be understood in the same way. 
Moreover, as $\xi_1=I\xi_0$ and so on, the above derivation formula implies
$$\split
\nabla \pmatrix \xi_1&\xi_2&\xi_3\endpmatrix&
=\pmatrix \xi_1&\xi_2&\xi_3 \endpmatrix \otimes 
\pmatrix 0 & 2\alpha_3 & -2 \alpha_2\\
-2\alpha_3 & 0 & 2\alpha_1\\
2\alpha_2 & -2\alpha_1 & 0\endpmatrix\\
&\quad +\pmatrix I\nabla\xi_0& J\nabla\xi_0 & K\nabla\xi_0\endpmatrix
\endsplit
\tag3-3$$
and $\nabla\xi_0$ {\it etc.} should be understood as
$$
\left\{\aligned
& \nabla\xi_0=\Gamma_0\xi_0+(\Gamma_1+\alpha_1)\xi_1
+(\Gamma_2+\alpha_2)\xi_2+(\Gamma_3+\alpha_3)\xi_3\,\,\\
& I\nabla\xi_0=\Gamma_0\xi_1-(\Gamma_1+\alpha_1)\xi_0
+(\Gamma_2+\alpha_2)\xi_3-(\Gamma_3+\alpha_3)\xi_2\,\,\\
& J\nabla\xi_0=\Gamma_0\xi_2-(\Gamma_1+\alpha_1)\xi_3
-(\Gamma_2+\alpha_2)\xi_0+(\Gamma_3+\alpha_3)\xi_1\,\,\\
& K\nabla\xi_0=\Gamma_0\xi_3+(\Gamma_1+\alpha_1)\xi_2
-(\Gamma_2+\alpha_2)\xi_1-(\Gamma_3+\alpha_3)\xi_0\,\,.
\endaligned\right.
\tag3-4
$$
Although these derivation formulae involve terms which are ``invisible'' at $z\in\Cal Z$, 
we must take these terms into account, because their derivatives become significant 
in the curvature computation. 

We introduce the system of 1-forms $\{X^0,X^1,X^2,X^3\}$ on $\Cal Z$ 
which should be understood as a column $n$ vector (i.e., each $X^i$ 
represents an $n$-dimensional subspace in 
the cotangent space of $M$) by requiring that the system 
$\{X^0,X^1,X^2,X^3\}$ forms a basis dual to $\{\xi_0,\xi_1=I\xi_0,\xi_2=J\xi_0,\xi_3=K\xi_0\}$  
in $\wcalD$ at each point of $\wcalZ$ and annihilates the orthogonal complement 
$(\wcalD_{\wt z})^{\perp}$ w.r.to the basic canonical deformation metric $g_1^{\can}$ 
on $\wcalZ$. 
Of course the system of 1-forms $\{X^0,X^1,X^2,X^3\}$ is obtained from 
the map germ $\sigma$ and the extension of $\{I,J,K\}$ in the fiber $\PP^1_m$ 
direction by the pull-back of the canonical 1-forms on the holonomy reduction $\Cal P$ 
of the oriented orthonormal frames of $(M,g)$. 
Note that we can define the basic canonical deformation metric on $\wcalZ$ 
just in the same way as on $\Cal Z$. To do this we just replace the Fubini-Study 
metric on $\PP^1$ with curvature $4$ by the bi-invariant metric on $\SU(2)$ with 
sectional curvature identically $1$. 
The dual version of the above formulae can be described in terms of the covariant 
derivative (or exterior derivative) of $X^i$'s
$$
\left\{\aligned
& \nabla \pmatrix X^1\\ X^2\\ X^3\endpmatrix
+ \pmatrix 0 & 2\alpha_3 & -2 \alpha_2\\
-2\alpha_3 & 0 & 2\alpha_1\\
2\alpha_2 & -2\alpha_1 & 0\endpmatrix \otimes 
\pmatrix X^1\\ X^2\\ X^3 \endpmatrix 
+\pmatrix I\nabla X^0\\ J\nabla X^0\\ K\nabla X^0\endpmatrix
=\pmatrix 0\\ 0\\ 0\endpmatrix\\
& \nabla X^0=-\Gamma_0\otimes X^0+(\Gamma_1+\alpha_1)\otimes X^1
+(\Gamma_2+\alpha_2)\otimes X^2+(\Gamma_3+\alpha_3)\otimes X^3
\endaligned\right.
$$
or equivalently
$$
\left\{
\aligned
& dX^1-2\alpha_2\wedge X^3+2\alpha_3\wedge X^2+IdX^0=0\,\,,\\
& dX^2-2\alpha_3\wedge X^1+2\alpha_1\wedge X^3+JdX^0=0\,\,,\\
& dX^3-2\alpha_1\wedge X^2+2\alpha_2\wedge X^1+KdX^0=0\,\,,\\
& dX^0=-\Gamma_0\wedge X^0+(\Gamma_1+\alpha_1)\wedge X^1
+(\Gamma_2+\alpha_2)\wedge X^2+(\Gamma_3+\alpha_3)\wedge X^3
\endaligned\right.
\tag3-5
$$
(all coefficients are ``invisible'' at $z\in\Cal Z$). 
These formulae hold on $\Cal P$ before the projection onto the (extended) twistor space. 
Therefore we can speak of the Sasakian structure on $\wcalZ$ and its transversal 
K\"ahler structure. 
The merit of considering $\wcalZ$ is that we have a canonical choice of 
the triple $\{I,J,K\}$ ($J=J(z)$, $z\in\Cal Z$) along the $\SO(2)$-fiber over $z$, 
while in the above construction on $\Cal Z$ the choice of $\{I,K\}$ is determined 
only modulo $\SO(2)$-rotation. Once $\xi_0$ is chosen as in the above discussion, 
we get an orthonormal frame germ $\{\xi_0,I\xi_0,J\xi_0,K\xi_0\}$ at $m$. 
By this procedure we have defined a (germ of) section of the extended 
twistor bundle $\wcalZ \rightarrow M$ which is horizontal in 
the $\{\xi_0,J\xi_0\}_{\text{\rm span}}$-direction 
but {\it not} horizontal in the $\{I\xi_0,K\xi_0\}_{\text{\rm span}}$-direction. 
We have thus defined a $(4n)$-dimensional distribution $\wcalD$ on $\wcalZ$ 
which inherits the same property as $\Cal D$ originally defined on $\Cal Z$. 

We are now ready to define the Z-metric on the extended twistor space $\wcalZ$. 
We note that the system of 1-form germs
$$
\left\{\aligned & \alpha_1-\alpha_1(\xi_0)X^0-\alpha_1(\xi_1)X^1
-\alpha_1(\xi_2)X^2-\alpha_1(\xi_3)X^3\\
& \alpha_2-\alpha_2(\xi_0)X^0-\alpha_2(\xi_1)X^1-\alpha_2(\xi_2)X^2-\alpha_2(\xi_3)X^3\\
& \alpha_3-\alpha_3(\xi_0)X^0-\alpha_3(\xi_1)X^1-\alpha_3(\xi_2)X^2-\alpha_3(\xi_3)X^3
\endaligned\right.
$$
at $\wt z\in\wcalZ$ (respectively
$$
\left\{\aligned  & \alpha_1-\alpha_1(\xi_0)X^0-\alpha_1(\xi_1)X^1
-\alpha_1(\xi_2)X^2-\alpha_1(\xi_3)X^3\\
& \alpha_3-\alpha_3(\xi_0)X^0-\alpha_3(\xi_1)X^1-\alpha_3(\xi_2)X^2-\alpha_3(\xi_3)X^3
\endaligned\right.
$$
at $z\in \Cal Z$) annihilates the $(4n)$-dimensional distribution $\wcalD$ on $\wcalZ$ 
(resp. $\Cal D$ on $\Cal Z$). However, these 1-form germs are not mutually 
orthogonal at $\wt z\in\wcalZ$ (resp. at $z\in\Cal Z$) w.r.to the canonical deformation 
metrics. 
Indeed, although we have $\alpha_1(\xi_i)=\alpha_2(\xi_i)=\alpha_3(\xi_i)=0$ 
($i=0,2$) at $\wt z\in\wcalZ$ and $\alpha_2=0$ for every direction at $z\in\Cal Z$, 
we have $\alpha_i(\xi_j)\not=0$ for $i,j=1,3$ (resp. for $i=1,2,3$ and $j=1,3$) 
at $z\in\Cal Z$ (resp. $\wt z\in\wcalZ$) which make these 1-forms mutually 
non orthogonal w.r.to $g^{\can}_{\lambda}$. 
Here we should write ${}^t\alpha_k(\xi_i)X^i$ instead of $\alpha_k(\xi_i)X^i$ because 
$\alpha_k(\xi_i)$ should be interpreted as a row vector which makes ``inner product'' 
with a column vector $X^i$. However, for the brevity, we have omit the transposition sign. 
We define the family of Z-metrics $\wcalF^{\ZZ}$ on $\wcalZ$ in the following way. 
For the definition, we use the basic canonical deformation metric.$g^{\can}_1$. 
Before writing down the explicit form of the Z-metrics, we must check the effect that 
the vectors $\xi_1$ and $\xi_3$ are not horizontal. We recall that both 
$\xi_1$ and $\xi_3$ decompose into the orthogonal sum of column $(n-1)$-vector 
in $\Cal D_z'$ (interpreted as an $(n-1)$-dimensional subspace) and a unit vector 
in $L_{J(z)}'$. The length of these components w.r.to the basic canonical deformation 
metric $g_1^{\can}$ are $n-1$ and $\sqrt{2}$ respectively. Therefore the quadruple
$$\biggl\{\xi_0\,\,,\,\,\sqrt{\frac{n}{n+1}}\xi_1\,\,,\,\,\xi_2\,\,,\,\,\sqrt{\frac{n}{n+1}}\xi_3\biggr\}$$
is an orthonormal basis of $\Cal D_z$ w.r.to the basic canonical deformation metric. 
We introduce normalized version of $\xi_i$ ($i=1,3$) and $X^i$ ($i=1,3$) by 
putting
$$
\xi_1':=\sqrt{\frac{n}{n+1}}\xi_1\,\,,\,\,\xi_3':=\sqrt{\frac{n}{n+1}}\xi_3\,\,,\,\,
X^1_n:=\sqrt{\frac{n+1}{n}}X^1\,\,,\,\,X^3_n:=\sqrt{\frac{n+1}{n}}X^3\,\,.
$$
Then $\{\xi_0,\xi_1',\xi_2,\xi_3'\}$ is an orthonormal basis of $\Cal D_z$ w.r.to the basic 
canonical deformation metric $g^{\can}_1$ and the quadruple
$\{X^0,X^1_n,X^2,X^2_n\}$ is the dual basis of $(\Cal D_z')^*$. 
We extend these 1-forms so that these vanish on vectors in $(\Cal D_z)^{\perp}$ 
where the orthogonal complement $\perp$ is taken w.r.to the 
basic canonical deformation metric $g_1^{\can}$. 
We thus have a frame germ
$$
\{\lambda(\alpha_1-\sum_{i=0}^3\alpha_1(\xi_i)X^i)\,\,,\,\,
\lambda(\alpha_2-\sum_{i=0}^3\alpha_2(\xi_i)X^i)\,\,,\,\,
\lambda(\alpha_3-\sum_{i=0}^3\alpha_3(\xi_i)X^i)\,\,,\,\,
X^0\,\,,\,\,X^1_n\,\,,\,\,X^2\,\,,\,\,X^3_n\}\,\,.
$$
One would like to define a family of metrics on $\wcalZ$ by declaring that 
the above frame germ (after applying some orthogonalization procedure 
to the first three 1-forms w.r.to the basic canonical deformation metric) 
being orthonormal coframe (modulo total scaling parameter). 
However, this attempt of defining new metrics does not fit with our strategy 
of constructing new metrics on $\Cal Z$ (or $\wcalZ$) which is ``closer'' 
to the independently scaled product metrics of $S^n\times S^n$. 
Instead, we define the family of Z-metrics by declaring that the frame germ
$$
\{\lambda(\alpha_1-\sum_{i=0}^3\alpha_1(\xi_i)X^i)\,\,,\,\,
\lambda(\alpha_2-\sum_{i=0}^3\alpha_2(\xi_i)X^i)\,\,,\,\,
\lambda(\alpha_3-\sum_{i=0}^3\alpha_3(\xi_i)X^i)\,\,,\,\,
X^0\,\,,\,\,X^1\,\,,\,\,X^2\,\,,\,\,X^3\}\,\,.
$$
(without normalization on $X^1$ and $X^3$) being 
orthonormal (modulo total scaling parameter $\rho$) : 
$$
\split
\rho\,\wt g_{\lambda}^Z&:=
\rho\,\biggl[\lambda^2\,
\biggl\{\biggl(\alpha_1-\alpha_1(\xi_0)X^0-\alpha_1(\xi_1)X^1
-\alpha_1(\xi_2)X^2-\alpha_1(\xi_3)X^3\biggr)^2\\
& \quad +\biggl(\alpha_2-\alpha_2(\xi_0)X^0-\alpha_2(\xi_1)X^1-\alpha_2(\xi_2)X^2
-\alpha_2(\xi_3)X^3\biggr)^2\\
& \quad +\biggl(\alpha_3-\alpha_3(\xi_0)X^0-\alpha_3(\xi_1)X^1
-\alpha_3(\xi_2)X^2-\alpha_3(\xi_3)X^3\biggr)^2\biggr\}\\
& \quad 
+\,{}^tX^0\cdot X^0+{}^tX^1\cdot X^1+{}^tX^2\cdot X^2+{}^tX^3\cdot X^3
\biggr]\,\,.
\endsplit
$$
The family $\Cal F^{\ZZ}$ of Z-metrics on the twistor space $\Cal Z$ is defined 
by declaring that the $\SO(2)$-bundle $\wcalZ \rightarrow \Cal Z$ 
is a Riemannan submersion w.r.to Z-metrics on $\wcalZ$ and $\Cal Z$ : 
$$
\split
\rho\,g_{\lambda}^Z&:=
\rho\,\biggl[\lambda^2\,
\biggl\{\biggl(\alpha_1-\alpha_1(\xi_0)X^0-\alpha_1(\xi_1)X^1
-\alpha_1(\xi_2)X^2-\alpha_1(\xi_3)X^3\biggr)^2\\
& \quad +\biggl(\alpha_3-\alpha_3(\xi_0)X^0-\alpha_3(\xi_1)X^1
-\alpha_3(\xi_2)X^2-\alpha_3(\xi_3)X^3\biggr)^2\biggr\}\\
&\quad 
+\,{}^tX^0\cdot X^0+{}^tX^1\cdot X^1+{}^tX^2\cdot X^2+{}^tX^3\cdot X^3
\biggr]\,\,.
\endsplit
$$
Here $\lambda$ is a positive partial scaling parameter. This expression should be understood 
as an expression in terms of 1-form germs at $z\in\Cal Z$ under consideration. 
These expressions define Z-metrics by specifying the oriented orthonormal coframe 
at one point $\wt z\in\wcalZ$ and $z\in\Cal Z$. Moreover, we regard these as metric germs 
at $\wt z\in\wcalZ$ or $z\in\Cal Z$ and compute their curvature form (we included ``invisible'' 
terms in the definition of the Z-metrics because these terms are significant in the computation 
of the Levi-Civita connection and the curvature form by differentiation). 
Even if the above expression is in the nice form only at one point ($z$ and $\wt z$) 
under consideration, the curvature computation regarding these metrics as germs is justified 
(we show this later). 
In these expressions, we note that the meaning of $X^i$'s ($i=0,1,2,3$) 
in the definition of the Z-metrics is {\it not the same} as that in the definition 
of the canonical deformation metrics. 
In the definition of the canonical deformation metrics, the corresponding orthonormal basis 
consists of vertical and horizontal vectors, while in the definition of Z-metrics the corresponding 
orthonormal basis does not consist of vertical / horizontal vectors. 
Indeed, $\alpha_1-\alpha_1(\xi_1)X_1-\alpha_1(\xi_3)X^3$ and 
$\alpha_3-\alpha_3(\xi_1)X^1-\alpha_3(\xi_3)X^3$ (resp.  $\xi_1$ 
and $\xi_3$) are not vertical (resp. not horizontal)
\footnote{\,\,However, $X^i$'s ($i=0,1,2,3$ have the {\it same} meaning in the sense that 
in either cases they stem from the canonical 1-forms defined on $\Cal P$.}, where we ignored 
``invisible'' parts at the point under consideration. 

We set 
$$
\left\{\aligned
& \hat\alpha_1:=\alpha_1-\alpha_1(\xi_0)X^0-\alpha_1(\xi_1)X^1
-\alpha_1(\xi_2)X^2-\alpha_1(\xi_3)X^3\,\,,\\
& \hat\alpha_3;=\alpha_3-\alpha_3(\xi_0)X^0-\alpha_3(\xi_1)X^1
-\alpha_3(\xi_2)X^2-\alpha_3(\xi_3)X^3\,\,.
\endaligned\right.
$$
These are regarded as 1-form germs at a point $z\in\Cal Z$ under consideration. 
We have from (2-12) the formulae of covariant derivatives w.r.to the Levi-Civita connection 
of the basic canonical deformation metric $g_1^{\can}$:
$$
\left\{\aligned
& \nabla\alpha_1=2\alpha_2\otimes\alpha_3-{}^tX^0\otimes X^1
+{}^tX^1\otimes X^0+{}^tX^2\otimes X^3-{}^tX^3\otimes X^2\,\,,\\
& \nabla\alpha_3=-2\alpha_2\otimes\alpha_1-{}^tX^0\otimes X^3
+{}^tX^3\otimes X^0+{}^tX^1\otimes X^2-{}^tX^2\otimes X^1\,\,,
\endaligned\right.
\tag3-6
$$
where $\nabla$ is the Levi-Civita connection of $g$ and the convention is that 
the Levi-Civita connection form appears like 
$(\text{\rm connection form}\otimes \text{\rm section})$ in the above formulae 
(we just replace $\otimes$ by $\wedge$ to get (2-12) from the above covariant derivative 
formula in this convention). 

We prepare some useful formulae. 
Applying $K$ and $I$ to the moving frame formula for $dX^0$ 
(note that $IX^0=-X^1$, $JX^0=-X^2$, $KX^0=-X^3$ and so on), 
we have
$$
\left\{\aligned
& dX^0=-\Gamma_0\wedge X^0+(\Gamma_1+\alpha_1)\wedge X^1
+(\Gamma_2+\alpha_2)\wedge X^2+(\Gamma_3+\alpha_3)\wedge X^3\,\,,\\
& IdX^0=(\Gamma_1+\alpha_1)\wedge X^0+\Gamma_0\wedge X^1
+(\Gamma_3+\alpha_3)\wedge X^2-(\Gamma_2+\alpha_2)\wedge X^3\,\,,\\
& JdX^0=(\Gamma_2+\alpha_2)\wedge X^0-(\Gamma_3+\alpha_3)\wedge X^1
+\Gamma_0\wedge X^2+(\Gamma_1+\alpha_1)\wedge X^3\,\,,\\
& KdX^0=(\Gamma_3+\alpha_3)\wedge X^0+(\Gamma_2+\alpha_2)\wedge X^1
-(\Gamma_1+\alpha_1)\wedge X^2+\Gamma_0\wedge X^3\,\,.
\endaligned\right.
\tag3-7
$$
From (3-4), we have ($i=1,3$): 
$$
\left\{\aligned
& \alpha_i(\nabla\xi_0)=\alpha_i(\xi_0)\Gamma_0+\alpha_i(\xi_1)(\Gamma_1+\alpha_1)
+\alpha_i(\xi_2)(\Gamma_2+\alpha_2)+\alpha_i(\xi_3)(\Gamma_3+\alpha_3)\,\,,\\
& \alpha_i(I\nabla\xi_0)=\alpha_i(\xi_1)\Gamma_0-\alpha_i(\xi_0)(\Gamma_1+\alpha_1)
+\alpha_i(\xi_3)(\Gamma_2+\alpha_2)-\alpha_i(\xi_2)(\Gamma_3+\alpha_3)\,\,,\\
& \alpha_i(J\nabla\xi_0)=\alpha_i(\xi_2)\Gamma_0-\alpha_i(\xi_3)(\Gamma_1+\alpha_1)
-\alpha_i(\xi_0)(\Gamma_2+\alpha_2)+\alpha_i(\xi_1)(\Gamma_3+\alpha_3)\,\,\\
& \alpha_i(K\nabla\xi_0)=\alpha_i(\xi_3)\Gamma_0+\alpha_i(\xi_2)(\Gamma_1+\alpha_1)
-\alpha_i(\xi_1)(\Gamma_2+\alpha_2)-\alpha_i(\xi_0)(\Gamma_3+\alpha_3)\,\,.
\endaligned\right.
\tag3-8
$$
We recall that $\Gamma_0$ and $\Gamma_2$ are not ``invisible'' in the $\xi_1$- and $\xi_3$-
directions and moreover $\Gamma_1$ is invisible in the $\xi_0$- and $\xi_2$-directions but 
not invisible in the $\xi_1$ and $\xi_2$-directions. 
From (3-7) and (3-8), we have, modulo terms vanishing to order $\geq 2$ at $z\in\Cal Z$, the 
following formulae: 
$$
\split
& \alpha_3(\xi_1)IdX^0+\alpha_3(\xi_3)KdX^0\\
&=\alpha_3(\nabla\xi_0)\wedge X^0
+\alpha_3(I\nabla\xi_0)\wedge X^1
+\alpha_3(J\nabla\xi_0)\wedge X^2
+\alpha_3(K\nabla\xi_0)\wedge X^3\\
& \quad +\{\alpha_3(\xi_0)(\Gamma_1+\alpha_1)+\alpha_3(\xi_2)(\Gamma_3+\alpha_3)\}\wedge X^1\\
&\quad +\{-\alpha_3(\xi_2)(\Gamma_1+\alpha_1)+\alpha_3(\xi_0)(\Gamma_3+\alpha_3)\}\wedge X^3\\
&\quad -\alpha_3(\xi_0)\Gamma_0\wedge X^0-\alpha_3(\xi_2)\Gamma_2\wedge X^0
-\alpha_3(\xi_2)\Gamma_0\wedge X^2-\alpha_3(\xi_0)\Gamma_2\wedge X^2\,\,,\\
& \alpha_1(\xi_3)KdX^0+\alpha_1(\xi_1)IdX^0\\
&=\alpha_1(\nabla\xi_0)\wedge X^0
+\alpha_1(I\nabla\xi_0)\wedge X^1
+\alpha_1(J\nabla\xi_0)\wedge X^2
+\alpha_1(K\nabla\xi_0)\wedge X^3\\
& \quad 
+\{\alpha_1(\xi_0)(\Gamma_1+\alpha_1)+\alpha_1(\xi_2)(\Gamma_3+\alpha_3)\}\wedge X^1\\
&\quad +\{-\alpha_1(\xi_2)(\Gamma_1+\alpha_1)+\alpha_1(\xi_0)(\Gamma_3+\alpha_3)\}\wedge X^3\\
&\quad -\alpha_1(\xi_0)\Gamma_0\wedge X^0-\alpha_1(\xi_2)\Gamma_2\wedge X^0
-\alpha_1(\xi_2)\Gamma_0\wedge X^2-\alpha_1(\xi_0)\Gamma_2\wedge X^2
\,\,.
\endsplit
\tag3-9
$$
Moreover we have
$$
\split
& dX^0=(\Gamma_1+\alpha_1)\wedge X^1+(\Gamma_3+\alpha_3)\wedge X^3\,\,,\\
& IdX^0=(\Gamma_1+\alpha_1)\wedge X^0+(\Gamma_3+\alpha_3)\wedge X^2\,\,,\\
& JdX^0=-(\Gamma_3+\alpha_3)\wedge X^1+(\Gamma_1+\alpha_1)\wedge X^3\,\,,\\
& KdX^0=(\Gamma_3+\alpha_3)\wedge X^0-(\Gamma_1+\alpha_1)\wedge X^2
\endsplit
\tag3-10
$$
at $z\in\Cal Z$ modulo ``invisible'' terms. 
Now we compute the derivation formulae for $\hat\alpha_1$ and $\hat\alpha_3$ 
modulo terms vanishing to order $\geq 2$ at $z\in\Cal Z$. 
We have
$$
\split
d\hat\alpha_1&=d\biggl(\alpha_1-\sum_{i=9}^3\alpha_1(\xi_i)X^i\biggr)\\
&=d\alpha_1-\sum_{i=0}^3\biggl\{(\nabla\alpha_1)(\xi_i)\wedge X^i
+\alpha_1(\nabla\xi_i)\wedge X^i+\alpha_1(\xi_i)dX^i\biggr\}
\endsplit
$$
where $\nabla$ is the Levi-Civita connection of the basic canonical deformation 
metric $g^{\can}_1$. 
Applying (2-12), (3-6),  (3-3) and (3-5) to $d\alpha_1$, $\nabla\alpha_1$, $\nabla\xi_i$ 
and $dX^i$, we compute the right hand side as
$$
\split
& 2\alpha_2\wedge\alpha_3+2({}^tX^1\wedge X^0+{}^tX^2\wedge X^3)\\
& \quad  -(2\alpha_2\otimes\alpha_3-{}^tX^0\otimes X^1
+{}^tX^1\otimes X^0+{}^tX^2\otimes X^3-{}^tX^3\otimes X^2)(\xi_0)\wedge X^0\\
& \quad  -\alpha_1(\nabla\xi_0)\wedge X^0\\
& \quad -\alpha_1(\xi_0)dX^0\\
& \quad  -(2\alpha_2\otimes\alpha_3-{}^tX^0\otimes X^1
+{}^tX^1\otimes X^0+{}^tX^2\otimes X^3-{}^tX^3\otimes X^2)(\xi_1)\wedge X^1\\
& \quad  -\alpha_1(2\alpha_2\otimes\xi_3-2\alpha_3\otimes\xi_2+I\nabla\xi_0)\wedge X^1\\
& \quad  -\alpha_1(\xi_1)(2\alpha_2\wedge X^3-2\alpha_3\wedge X^2-IdX^0)\\
& \quad  -2(2\alpha_2\otimes\alpha_3-{}^tX^0\otimes X^1
+{}^tX^1\otimes X^0+{}^tX^2\otimes X^3-{}^tX^3\otimes X^2)(\xi_2)\wedge X^2\\
& \quad  -\alpha_1(2\alpha_3\otimes\xi_1-2\alpha_1\otimes\xi_3+J\nabla\xi_0)\wedge X^2\\
& \quad  -\alpha_1(\xi_2)(2\alpha_3\wedge X^1-2\alpha_1\wedge X^3-JdX^0)\\
& \quad  -(2\alpha_2\otimes\alpha_3-{}^tX^0\otimes X^1
+{}^tX^1\otimes X^0+{}^tX^2\otimes X^3-{}^tX^3\otimes X^2)(\xi_3)\wedge X^3\\
& \quad  -\alpha_1(2\alpha_1\otimes\xi_2-2\alpha_2\otimes\xi_1+K\nabla\xi_0)\wedge X^3\\
& \quad  -\alpha_1(\xi_3)(2\alpha_1\wedge X^2-2\alpha_2\wedge X^1-KdX^0)\,\,.
\endsplit
$$
Executing the cancellation and applying (18), we can reduce this, modulo terms vanishing 
at $z \in \Cal Z$ to order $\geq 2$, to the following : 
$$
\split
& 2\alpha_2\wedge\alpha_3-2\alpha_3(\xi_0)\alpha_2\wedge X^0
-2\alpha_3(\xi_1)\alpha_2\wedge X^1
-2\alpha_3(\xi_2)\alpha_2\wedge X^2-\alpha_3(\xi_3)\alpha_2\wedge X^3\\
& \quad -\alpha_1(\xi_0)\{(\Gamma_1+\alpha_1)\wedge X^1+(\Gamma_3+\alpha_3)\wedge X^3\}\\
&\quad +\alpha_1(\xi_2)\{-(\Gamma_3+\alpha_3)\wedge X^1+(\Gamma_1+\alpha_1)\wedge X^3\}\\
&\quad -\alpha_1(\nabla\xi_0)\wedge X^0-\alpha_1(I\nabla\xi_0)\wedge X^1
-\alpha_1(J\nabla\xi_0)\wedge X^2-\alpha_1(K\nabla\xi_0)\wedge X^3\\
& \quad +\alpha_1(\xi_3)KdX^0+\alpha_1(\xi_1)IdX^0
\,\,.
\endsplit
$$
Using the definition of $\hat\alpha_3$ and applying (17), this becomes
$$
\split
& \quad 2\alpha_2\wedge \hat\alpha_3
-\alpha_1(\xi_0)\{(\Gamma_1+\alpha_1)\wedge X^1+(\Gamma_3+\alpha_3)\wedge X^3\}\\
&\quad +\alpha_1(\xi_2)\{-(\Gamma_3+\alpha_3)\wedge X^1+(\Gamma_1+\alpha_1)\wedge X^3\}\\
&\quad 
+\{\alpha_1(\xi_0)(\Gamma_1+\alpha_1)+\alpha_1(\xi_2)(\Gamma_3+\alpha_3)\}\wedge X^1\\
&\quad +\{-\alpha_1(\xi_2)(\Gamma_1+\alpha_1)+\alpha_1(\xi_0)(\Gamma_3+\alpha_3)\}\wedge X^3\\
&\quad -\alpha_1(\xi_0)\Gamma_0\wedge X^0-\alpha_1(\xi_2)\Gamma_2\wedge X^0
-\alpha_1(\xi_2)\Gamma_0\wedge X^2-\alpha_1(\xi_0)\Gamma_2\wedge X^2
\endsplit
$$
and finally we get the derivation formula for $\hat\alpha_1$ (modulo terms vanishing 
at $z$ to order $\geq 2$): 
$$
\split
d\hat\alpha_1&=2\alpha_2\wedge \hat\alpha_3\\
&\quad -\alpha_1(\xi_0)\Gamma_0\wedge X^0-\alpha_1(\xi_2)\Gamma_2\wedge X^0
-\alpha_1(\xi_2)\Gamma_0\wedge X^2-\alpha_1(\xi_0)\Gamma_2\wedge X^2
\,\,.
\endsplit
\tag3-11
$$
Next we compute $d\hat\alpha_3$. Computing just as in the same way, we have
$$
\split
d\hat\alpha_3&=d\biggl(\alpha_3-\sum_{i=0}^3\alpha_3(\xi_i)X^i\biggr)\\
&=d\alpha_3-\sum_{i=0}^3\biggl\{(\nabla\alpha_3)(\xi_i)X^i+\alpha_3(\nabla\xi_i)X^i
+\alpha_3(\xi_i)dX^i\biggr\}\\
&=-2\alpha_2\wedge\alpha_1+2({}^tX^3\wedge X^0+{}^tX^1\wedge X^2)\\
& \quad  -(-2\alpha_2\otimes\alpha_1-{}^tX^0\otimes X^3
+{}^tX^3\otimes X^0+{}^tX^1\otimes X^2-{}^tX^2\otimes X^1)(\xi_0)\wedge X^0\\
& \quad  -\alpha_3(\nabla\xi_0)\wedge X^0\\
&\quad -\alpha_3(\xi_0)dX^0\\
& \quad  -(-2\alpha_2\otimes\alpha_1-{}^tX^0\otimes X^3
+{}^tX^3\otimes X^0+{}^tX^1\otimes X^2-{}^tX^2\otimes X^1)(\xi_1)\wedge X^1\\
& \quad  -\alpha_3(2\alpha_2\otimes\xi_3-2\alpha_3\otimes\xi_2+I\nabla\xi_0)\wedge X^1\\
& \quad  -\alpha_3(\xi_1)(2\alpha_2\wedge X^3-2\alpha_3\wedge X^2-IdX^0)\\
& \quad  -2(-2\alpha_2\otimes\alpha_1-{}^tX^0\otimes X^3
+{}^tX^3\otimes X^0+{}^tX^1\otimes X^2-{}^tX^2\otimes X^1)(\xi_2)\wedge X^2\\
& \quad  -\alpha_3(2\alpha_3\otimes\xi_1-2\alpha_1\otimes\xi_3+J\nabla\xi_0)\wedge X^2\\
& \quad  -\alpha_3(\xi_2)(2\alpha_3\wedge X^1-2\alpha_1\wedge X^3-JdX^0)\\
& \quad  -(-2\alpha_2\otimes\alpha_1-{}^tX^0\otimes X^3
+{}^tX^3\otimes X^0+{}^tX^1\otimes X^2-{}^tX^2\otimes X^1)(\xi_3)\wedge X^3\\
& \quad  -\alpha_3(2\alpha_1\otimes\xi_2-2\alpha_2\otimes\xi_1+K\nabla\xi_0)\\
& \quad  -\alpha_3(\xi_3)(-2\alpha_1\wedge X^2+2\alpha_2\wedge X^1-KdX^0)\,\,.
\endsplit$$
Executing cancellation and using (18), the right hand side reduces to
$$
\split
&-2\alpha_2\wedge\alpha_1+2\alpha_1(\xi_0)\alpha_2\wedge X^0
+2\alpha_1(\xi_1)\alpha_2\wedge X^1
+2\alpha_1(\xi_2)\alpha_2\wedge X^2+\alpha_1(\xi_3)\alpha_2\wedge X^3\\
&\quad -\alpha_3(\xi_0)\{(\Gamma_1+\alpha_1)\wedge X^1+(\Gamma_3+\alpha_3)\wedge X^3\}\\
&\quad +\alpha_3(\xi_2)\{(-(\Gamma_3+\alpha_3)\wedge X^1+(\Gamma_1+\alpha_1)\wedge X^3\}\\
&\quad -\alpha_3(\nabla\xi_0)\wedge X^0-\alpha_3(I\nabla\xi_0)\wedge X^1
-\alpha_3(J\nabla\xi_0)\wedge X^2-\alpha_3(K\nabla\xi_0)\wedge X^3\\
& \quad +\alpha_3(\xi_1)IdX^0+\alpha_3(\xi_3)KdX^0
\,\,.
\endsplit
$$
We compute just in the same way as above using the definition of $\hat\alpha_1$ and applying 
(3-9), (3-10) to get the derivation formula of $\hat\alpha_3$ (modulo terms vanishing 
at $z$ to order $\geq 2$): 
$$\split
d\hat\alpha_3&=-2\alpha_2\wedge\hat\alpha_1
-\alpha_3(\xi_0)\{(\Gamma_1+\alpha_1)\wedge X^1+(\Gamma_3+\alpha_3)\wedge X^3\}\\
&\quad +\alpha_3(\xi_2)\{-(\Gamma_3+\alpha_3)\wedge X^1+(\Gamma_1+\alpha_1)\wedge X^3\}\\
&\quad +\{\alpha_3(\xi_0)(\Gamma_1+\alpha_1)+\alpha_3(\xi_2)(\Gamma_3+\alpha_3)\}\wedge X^1\\
&\quad +\{-\alpha_3(\xi_2)(\Gamma_1+\alpha_1)+\alpha_3(\xi_0)(\Gamma_3+\alpha_3)\}\wedge X^3\\
&\quad -\alpha_3(\xi_0)\Gamma_0\wedge X^0-\alpha_3(\xi_2)\Gamma_2\wedge X^0
-\alpha_3(\xi_2)\Gamma_0\wedge X^2-\alpha_3(\xi_0)\Gamma_2\wedge X^2\\
&=-2\alpha_2\wedge\hat\alpha_1\\
&\quad -\alpha_3(\xi_0)\Gamma_0\wedge X^0-\alpha_3(\xi_2)\Gamma_2\wedge X^0
-\alpha_3(\xi_2)\Gamma_0\wedge X^2-\alpha_3(\xi_0)\Gamma_2\wedge X^2\,\,.
\endsplit
\tag3-12
$$
Let us compare the derivation formulae (3-11) and (3-12) of $\hat\alpha_i$ ($i=1.3$) 
with those (2-12) of $\alpha_i$ ($i=1,3$). The remarkable difference is the following. 
In (2-12), the term such as $X^{\mu}\wedge X^{\nu}$ appears (geometrically this means 
that the curvature form of the canonical Hermitian metric of the $S^1$-bundle 
$\wcalZ\rightarrow \Cal Z$ is proportional to the K\"ahler-Einstein metric $g^{\can}_1$ 
on $\Cal Z$). On the other hand, although there are similar terms 
(such as $-\alpha_3(\xi_0)\Gamma_0\wedge X^0$) in (3-11) and (3-12), 
these terms appear with coefficients ``invisible'' at $z\Cal Z$ (such as $\alpha_3(\xi_0)$). 
We are thus tempted to think that the distribution $\Cal D$ 
on the twistor space $\Cal Z$ defined by the equation $\hat\alpha_1=\hat\alpha_3=0$ 
is much ``closer to being  integrable'' than is the horizontal distribution $\Cal H$ of the 
twistor fibration $\Cal Z\rightarrow M$ (w.r.to 
the Levi-Civita connection of the quaternion K\"ahler metric $g$ on $M$). 
In other words, we are tempted to think that the distribution 
$\Cal D=\{\hat\alpha_1=\hat\alpha_3=0\}$ on the 
twistor space $\Cal Z$ equipped with the Z-metrics $g^{\ZZ}_{\lambda}$ 
constructed from the distribution $\Cal D$ would be much ``closer'' to the product structure 
on, e.g., $S^n\times S^n$ equipped with the independently scaled 
product constant curvature metrics than is the horizontal distribution $\Cal H$ 
equipped with the canonical deformation metrics. 

We are now ready to compute the Levi-Civita connection of the Z-metrics. 
We recall that the $\Gamma$-part of the Levi-Civita connection of $g$ satisfies the condition 
that $\Gamma_0$ and $\Gamma_2$ are ``invisible'' in $\xi_0$- and $J\xi_0$-directions 
(this stems from the construction of the distribution $\Cal D$ on $\Cal Z$). 
This implies that we can write
$$
\split
& \Gamma_0=p\,X^1+q\,X^3\,\,,\\
& \Gamma_2=r\,X^1+s\,X^3\,\,.
\endsplit
$$
Here, these expressions should be understood in the column $n$ vector notation 
under the identification $\R^{4n}=\H^n$. We write down the first structure equation 
of the metric $g^{\ZZ}_{\lambda}$ on $\Cal Z$ w.r.to the moving frame 
$$\{\lambda\hat\alpha_1\,\,,\,\,\lambda\hat\alpha_3\,\,,\,\,X^0\,\,,\,\,X^1\,\,,\,\,X^2\,\,,\,\,X^3\}$$
at $z\in\Cal Z$. We recall that Z-metric $g^{\ZZ}_{\lambda}$ is defined by declaring 
that the above moving frame is an oriented orthonormal frame. 
We write this as a column vector as in the first structure equation in \S2, i.e., the connection 
matrix $\Gamma^{\ZZ}_{\lambda}$ is defined by
$$d\,\pmatrix \lambda\hat\alpha_1\\ \lambda\hat\alpha_3\\ X^0\\ X^1\\ X^2\\ X^3\endpmatrix
+\Gamma_{\lambda}^{\ZZ}\wedge 
\pmatrix \lambda\hat\alpha_1\\ \lambda\hat\alpha_3\\ X^0\\ X^1\\ X^2\\ X^3\endpmatrix
=\pmatrix 0\\0\\0\\0\\0\\0\endpmatrix\,\,.$$
For this purpose, we put, for $i=1,3$ : 
$$\split
& a_i:=\frac12\{\alpha_i(\xi_0)p+\alpha_i(\xi_2)r\}\,\,,\\
& b_i:=\frac12\{\alpha_i(\xi_0)q+\alpha_i(\xi_2)s\}\,\,,\\
& c_i:=\frac12\{\alpha_i(\xi_2)p+\alpha_i(\xi_0)r\}\,\,,\\
& d_i:=\frac12\{\alpha_i(\xi_2)q+\alpha_i(\xi_0)s\}\,\,.
\endsplit
$$
Note that these are all ``invisible'' at $z\in\Cal Z$. 
It follows from the derivation formulae (3-11) and (3-12) that the connection matrix $\Gamma^{\ZZ}_{\lambda}$ 
under question is written as
$$
\pmatrix
0&-2\alpha_2&\lambda(a_1X^1+b_1X^3)&-\lambda(a_1X^0+c_1X^2)&\lambda(c_1X^1+d_1X^3)&-\lambda(b_1X^0+d_1X^2)\\
2\alpha_2&0&\lambda((a_3X^1+b_3X^3)&-\lambda(a_3X^0+c_3X^2)&\lambda(c_3X^1+d_3X^3)&-\lambda(b_3X^0+d_3X^2)\\
{}&{}&                              {}&{-\Gamma_1-\alpha_1}&{}&{-\Gamma_3-\alpha_3}\\
{}&{}&                              {}&{-\lambda^2(a_1\alpha_1+a_3\alpha_3}&{}&{-\lambda^2(b_1\alpha_1+b_3\alpha_3}\\
{-\lambda(a_1X^1}&{-\lambda(a_3X^1}&{\Gamma_0}&{-a_1\alpha_1(\xi_1)X^1}&{-\Gamma_2+\alpha_2}&{-b_1\alpha_1(\xi_1)X^1}\\
{+b_1X^3)}&{+b_3X^3)}&{+0}&{-a_1\alpha_1(\xi_3)X^3}&{+0}&{-b_1\alpha_1(\xi_3)X^3}\\
{}&{}&                                {}&{-a_3\alpha_3(\xi_1)X^1}&{}&{-b_3\alpha_3(\xi_1)X^1}\\
{}&{}&                                {}&{-a_3\alpha_3(\xi_3)X^3)}&{}&{-b_3\alpha_3(\xi_3)X^3)}\\
{}&{}&                                {\Gamma_1+\alpha_1}&{}&{-\Gamma_3+\alpha_3}&{}\\
{}&{}&                                {+\lambda^2(a_1\alpha_1+a_3\alpha_3}&{}&{+\lambda^2(c_1\alpha_1+c_3\alpha_3}&{}\\
{\lambda(a_1X^0}&{\lambda(a_3X^0}&{-a_1\alpha_1(\xi_1)X^1}&{\Gamma_0}&{-c_1\alpha_1(\xi_1)X^1}&{\Gamma_2-\alpha_2}\\
{+c_1X^2)}&{+c_3X^2)}&{-a_1\alpha_1(\xi_3)X^3}&{+0}&{-c_1\alpha_1(\xi_3)X^3}&{+0}\\
{}&{}&                                {-a_3\alpha_3(\xi_1)X^1}&{}&{-c_3\alpha_3(\xi_1)X^1}&{}\\
{}&{}&                                {-a_3\alpha_3(\xi_3)X^3)}&{}&{-c_3\alpha_3(\xi_3)X^3)}&{}\\
{}&{}&                                {}&{\Gamma_3-\alpha_3}&{}&{-\Gamma_1+\alpha_1}\\
{}&{}&                                {}&{-\lambda^2(c_1\alpha_1+c_3\alpha_3}&{}&{-\lambda^2(d_1\alpha_1+d_3\alpha_3}\\
{-\lambda(c_1X^1}&{-\lambda(c_3X^1}&{\Gamma_2+\alpha_2}&{-c_1\alpha_1(\xi_1)X^1}&{\Gamma_0}&{-d_1\alpha_1(\xi_1)X^1}\\
{+d_1X^3)}&{+d_3X^3}&{+0}&{-c_1\alpha_1(\xi_3)X^3}&{+0}&{-d_1\alpha_1(\xi_3)X^3}\\
{}&{}&                                {}&{-c_3\alpha_3(\xi_1)X^1}&{}&{-d_3\alpha_3(\xi_1)X^1}\\
{}&{}&                                {}&{-c_3\alpha_3(\xi_3)X^3)}&{}&{-d_3\alpha_3(\xi_3)X^3)}\\
{}&{}&                                {\Gamma_3+\alpha_3}&{}&{\Gamma_1-\alpha_1}&{}\\
{}&{}&                                {+\lambda^2(b_1\alpha_1+b_3\alpha_3}&{}&{+\lambda^2(d_1\alpha_1+d_3\alpha_3}&{}\\
{\lambda(b_1X^0}&{\lambda(b_3X^0}&{-b_1\alpha_1(\xi_1)X^1}&{-\Gamma_2+\alpha_2}&{-d_1\alpha_1(\xi_1)X^1}&{\Gamma_0}\\
{+d_1X^2)}&{+d_3X^2)}&{-b_1\alpha_1(\xi_3)X^3}&{+0}&{-d_1\alpha_1(\xi_3)X^3}&{+0}\\
{}&{}&                                {-b_3\alpha_3(\xi_1)X^1}&{}&{-d_3\alpha_3(\xi_1)X^1}&{}\\
{}&{}&                                {-b_3\alpha_3(\xi_3)X^3)}&{}&{-d_3\alpha_3(\xi_3)X^3)}&{}
\endpmatrix
\tag3-13
$$
The skew-symmetric matrix $\Gamma^{\ZZ}_{\lambda}$ in (3-13) 
is the connection matrix of the Levi-Civita connection of the Z-metric 
$g^{\ZZ}_{\lambda}$. 

It follows from (3-13) and the first structure equation that $dX^i$ ($i=0,1,2,3$) 
is expressed as the sum of terms of the form $X^i\wedge (\text{\rm something})$ ($i=0,1,2,3$). 
Therefore the $2$-dimensional distribution on $\Cal Z$ defined by the system 
of equations $X^i=0$ ($i=0,1,2,3$) is integrable. 
Moreover this distribution is holomorphic w.r.to the canonical complex structure of $\Cal Z$. 
The expression (3-13), the first structure equation and the 
formula $d(2\alpha_2)=4\alpha_3\wedge\alpha_1+4({}^tX^2\wedge X^0+{}^tX^3\wedge X^1)$ 
imply that every integral submanifold of this distribution 
has constant curvature $4$. Therefore any integral submanifold of the distribution 
$X^i=0$ ($i=0,1,2,3$) is a holomorphic curve in $\Cal Z$ isomorphic to $\PP^1$ and more 
precisely a deformation of the twistor line\footnote{\,\,In the definition of $L_{J(z)}'$ we have chosen 
$1\in H^0(\Cal P^1,\Cal O)$. However, we can start with any element from $H^0(\PP^1,\Cal O)$ and 
develop a similar theory and if we start with a constant from $H^0(\PP^1,\Cal O)$ 
which is very small in absolute value, then the distribution $\Cal D_z$ should be very close to 
the horizontal distribution and therefore the resulting analogue of a Z-metric with $\lambda=1$ 
should be very close to the basic canonical deformation metric. 
Therefore the integral manifold under consideration must be a deformation of a twistor line.}. 

The curvature form $\Om^{\ZZ}_{\lambda}$ 
is computed from the second structure equation $\Om^{\ZZ}_{\lambda}
=d\Gamma^{\ZZ}
+\Gamma^{\ZZ}_{\lambda}\wedge \Gamma^{\ZZ}_{\lambda}$. 
As $(a_i,b_i,c_i,d_i)$ ($i=1,3$) are all ``invisible'' at $z\in\Cal Z$, we conclude that 
the $\Gamma^{\ZZ}_{\lambda}\wedge \Gamma^{\ZZ}_{\lambda}$-part of the 
curvature form $\Om^{\ZZ}_{\lambda}$ is the same as the corresponding part 
of the curvature form of the original quaternion K\"ahler manifold $(M,g)$. 
On the other hand, we can compute the $d\Gamma^{\ZZ}_{\lambda}$-part 
in the following way. We first observe that if we impose the condition (3-2$'$) 
on $\nabla\xi_0$, we conclude that the differentiation of the ``invisible'' terms 
gives rise to 1-forms such exactly like $X^i$, modulo $\alpha_1$ and $\alpha_3$ : 
$$
\split
& \alpha_1(\xi_0) \overset{\text{\rm differentiation}}\to\longmapsto X^1\,\,,\,\,
\alpha_1(\xi_2) \overset{\text{\rm differentiation}}\to\longmapsto -X^3\,\,,\,\,\\
& \alpha_3(\xi_0)  \overset{\text{\rm differentiation}}\to\longmapsto X^3\,\,,\,\,
\alpha_3(\xi_2)  \overset{\text{\rm differentiation}}\to\longmapsto X^1\,\,.
\endsplit
$$
For instance, we have 
$$\split
d(\alpha_1(\xi_0))&=(\nabla\alpha_1)(\xi_0)+\alpha_1(\nabla\xi_0)\\
&=(2\alpha_2\otimes\alpha_3-{}^tX^0\otimes X^1+{}^tX^1\otimes X^0+\cdots)(\cdot,\xi_0)\\
&\quad +\alpha_1(\Gamma_0\otimes\xi_0+(\Gamma_1+\alpha_1)\otimes\xi_1
+(\Gamma_2+\alpha_2)\otimes\xi_2+(\Gamma_3+\alpha_3)\otimes \xi_3)\\
&={}^tX^1\,\,\, \text{\rm modulo}\,\,\, \alpha_1\,\,\text{\rm and}\,\, \alpha_3\,\,.
\endsplit
$$ 
ignoring ``invisible'' terms. We have used the condition (3-2$'$) 
on the $\Gamma$-part in $\nabla\xi_0$. We can check other cases similarly. 
Secondly, we combine the above observation with the general formula
$$\Ric (e_i,e_j)=\sum_{k=1}^{\dim}g(\Om^j_k(e_i,e_k)e_j,e_j)\,\,.$$
We see that all terms in the connection matrix $\Gamma^{\ZZ}_{\lambda}$ 
whose origin stem from the $a_i(\xi_j)\Gamma_k\wedge X^l$-part 
($i=1,3$ and $j,k,l=0,2$) in the derivation formulae (3-11) and (3-12), 
i.e., all terms which appears with one of $(a_i,b_i,c_i,d_i)$'s ($i=1,3$), 
do not contribute to the Ricci tensor. This means that the Ricci form 
of the Z-metric $g^{\ZZ}_{\lambda}$ has the form 
of the sum of $\hat\alpha_1^2+\hat\alpha_3^2$ and 
$\sum_{i=0}^3{}^tX^i\cdot X^i$ with some coefficients. 
For the explicit computation, we introduce 
$\{\lambda^{-1}\xi_{-2},\lambda^{-1}\xi_{-1},\xi_0,\xi_1,\xi_2,\xi_3\}$ 
which is the orthonormal frame dual to the orthonormal coframe
$\{\lambda\hat\alpha_1,\lambda\hat\alpha_3,X^0,X^1,X^2,X^3\}$. 
As the first example, the non-trivial contribution to the Ricci tensor in 
the $\xi_{-2}$-direction comes from $d\alpha_2$. 
Recall that $d\alpha_2$ is given by the formula
$$d\alpha_2=2\alpha_3\wedge\alpha_1+2({}^tX^2\wedge X^0+{}^tX^3\wedge X^1)\,\,.$$
We have $X^j(\xi_i)=0$ for $i<0$ and $j\geq 0$ and moreover from the definition 
of $\hat\alpha_i$ we have $\hat\alpha_i\equiv\alpha_i$ modulo $X^j$ ($i=1,3$ and $j\geq 0$). 
Therefore we have
$$
\split
& \Ric^{\ZZ}_{\lambda}(\lambda^{-1}\xi_{-2},\lambda^{-1}\xi_{-2})
=g_{\lambda}^{\ZZ}((\Om^{\ZZ}_{\lambda})^{-2}_{-1}(\lambda^{-1}\xi_{-2},
\lambda^{-1}\xi_{-1})\lambda^{-1}\xi_{-2},\lambda^{-1}\xi_{-2})\\
&=g^{\ZZ}_{\lambda}((4\alpha_1\wedge \alpha_3)(\lambda^{-1}\xi_{-2},
\lambda^{-1}\xi_{-1})\lambda^{-1}\xi_{-2},
\lambda^{-1}\xi_{-2})=\frac{4}{\lambda^2}\,\,.
\endsplit
$$
As the second example, the non-trivial contribution to the Ricci tensor in the 
$X^0$-direction comes from the curvature form of the original metric $g$ : 
$$
\split
& \Ric^{\ZZ}_{\lambda}(\xi_0,\xi_0)=g^{\ZZ}_{\lambda}(\sum_{i=3,1}(\Om^{\ZZ}_{\lambda})^0_{-i}
(\xi_0,\lambda^{-1}\xi_i)\xi_0+\sum_{k=1}^3(\Om^{\ZZ}_{\lambda})^0_k(\xi_0\xi_k)\xi_0,\xi_0)\\
&=g(\sum_{k=1}^3\Om^0_k(\xi_0,\xi_k)\xi_0,\xi_0)=\Ric(\xi_0,\xi_0)=4n+8\,\,.
\endsplit
$$
Similarly, we have
$$\split
& \Ric^{\ZZ}_{\lambda}(\lambda^{-1}\xi_{-2},\lambda^{-1}\xi_{-2})
=\Ric^{\ZZ}_{\lambda}(\lambda^{-1}\xi_{-1},\lambda^{-1}\xi_{-1})=\frac4{\lambda^2}\,\,,\\
& \Ric^{\ZZ}_{\lambda}(\xi_0,\xi_0)=\Ric^{\ZZ}_{\lambda}(\xi_1,\xi_1)
=\Ric^{\ZZ}_{\lambda}(\xi_2,\xi_2)=\Ric^{\ZZ}_{\lambda}(\xi_3,\xi_3)=4n+8\,\,\\
& \text{\rm all other components of}\,\,\Ric^{\ZZ}_{\lambda}\,\,=0\,\,.
\endsplit
$$
Summing up the computations in \S3, we have the following Proposition. 

\proclaim{Proposition 3.1} The Ricci tensor of the Z-metric 
$$g^{\ZZ}_{\lambda}=\lambda^2(\hat\alpha_1^2+\hat\alpha_3^2)+\sum_{i=0}^3{}^tX^i\cdot X^i$$
on the twistor space $\Cal Z$ is given by the formula
$$\Ric^{\ZZ}_{\lambda}=\frac{4}{\lambda^2}\lambda^2(\hat\alpha_1^2+\hat\alpha_3^2)
+(4n+8)\sum_{i=0}^3{}^tX^i\cdot X^i\,\,.\tag3-14$$
In particular, the Z-metric $g^{\ZZ}_{\lambda}$ is Einstein 
if and only if $\lambda^2=\frac1{n+2}$. 
\endproclaim

Proposition 2.4 implies that the family of canonical deformation metrics $g^{\can}_{\lambda}$ 
contains two Einstein metrics, i.e., those for $\lambda^2=1$ and $\lambda^2=\frac1{n+1}$. 
From Proposition 3.1, the Z-metric $g^{\ZZ}_{\lambda}$ for $\lambda^2=\frac1{n+2}$ 
turns out to be the third Einstein metric.  This turns out to be non-isometric to 
any Einstein metric in the family of canonical deformation metrics.  

\proclaim{Remark 3.2} {\rm The Einstein metric $g^{\ZZ}_{\lambda}$ with 
$\lambda^2=\frac1{n+2}$ in the family of Z-metrics on $\Cal Z$ 
is not isometric to the Einstein metrics $g_{\lambda}^{\can}$ with $\lambda^2=1$ 
or $\lambda^2=\frac1{n+1}$ in the family of canonical deformation metrics.
}\endproclaim
\noindent
{\it Proof.} It is clear from the construction that, if $(M,g)$ is a Wolf space, 
then both canonical deformation metrics 
$g^{\can}_{\lambda}$ and Z-metrics $g^{\ZZ}_{\lambda}$ 
are homogeneous w.r.to the action of the isometry group of $(M,g)$ 
on the twistor space $\Cal Z$. The assertion is then easy to check 
in the model case $\PP^3(\C) \rightarrow \PP^1(\H)$. Indeed, the difference of 
these three Einstein metrics are visible in the relationship between the 
orthonormal basis of $T_z\Cal Z$,  the fixed complex line $L_{J(z)}$ 
(cf. discussion before (3-1)) in the horizontal subspace (of the twistor 
fibration $\PP^3(\C) \rightarrow \PP^1(\H)$) and the tangent space 
of the twistor line $\PP^1$ (the $\Sp(1)$-orbit). \qed

\proclaim{Remark 3.3}{\rm (1) Homogeneous Einstein metrics on 
$\PP^{2n+2}(\C)$ were classified by Ziller [Z]. The canonical deformation metrics 
$g^{\can}_1$ and $g^{\can}_{\frac1{n+1}}$ are the only homogeneous Einstein metrics 
on $\PP^{2n+1}$ up to homothety. 
Therefore the Z-metric $g^{\ZZ}_{\frac1{n+2}}$ on the twistor space 
of $M=\PP^n(\H)$ is not a $\Sp(n+1)$-homogeneous Einstein metric on the twistor space 
$\displaystyle \Cal Z=\frac{\Sp(n+1)}{\Sp(n)\times \SO(2)}=\PP^{2n+1}(\C)$. 
The reason why $g^{\ZZ}_{\lambda}$ on $\PP^{2n+1}(\C)$ is not $\Sp(n+1)$-homogeneous 
can be seen from the $\Sp(n+1)$-orbit decomposition of the Grassmannian 
$\text{\rm Grass}(\C^2,\C^{2n+2})$ of all $\PP^1$'s (lines's) in $\PP^{2n+1}$, 
which is the complexfication of $\PP^n(\H)$. Indeed, there exists only one $(4n)$-dimensional 
$\Sp(n+1)$-orbit, which is $\PP^n(\H)$ corresponding to the twistor fibration. It would be 
interesting to determine the group of isometries of $g_{\frac1{n+2}}^{\ZZ}$. 
Is $g^{\ZZ}_{\frac1{n+2}}$ of cohomogeneity one ?
\medskip

(2) The fiber of the $\Sp(1)$-principal bundle 
$\wcalZ \rightarrow M$ is the parameter space of all ``framed'' $L_{J(z)}'$'s ($z$ lying on a fixed 
twistor line) and therefore the extended twistor space $\wcalZ$ itself is identified 
with the space of all ``framed'' $L_{J(z)}'$ ($z\in\Cal Z$). 
Since the distribution defined by $L_{J(z)}'$ is integrable\footnote{\,\,We have shown that 
the distribution $\Cal D^{\perp}$ defined by $X^i=0$ ($i=0,1,2,3$) is integrable 
(see the arguments after (3.13)). The statement follows because $\Cal D^{\perp}$ 
is of the same type as the distribution defined by $L_{J(z)}'$ and therefore the same 
reasoning works in the proof of the integrability.}, we can introduce an equivalent relation 
on $\wcalZ$ where two points of $\wcalZ$ are equivalent if and only if 
these lie on the same integral submanifold (extended by the $S^1$-fiber 
of the fibration $\wcalZ\rightarrow\Cal Z$) of the distribution $\Cal D^{\perp}$. 
The same argument works if we start with the distribution $\Cal D^{\perp}$ 
instead of the distribution defined by $L_{J(z)}'$ and therefore we can 
identify $\wcalZ$ with the space of ``framed'' $\Cal D^{\perp}$. 
The quotient space obtained from $\wcalZ$ (moduli space of ``framed'' $\Cal D^{\perp}$) 
and the above equivalence relation (defined by the integral submanifold of $\Cal D^{\perp}$) 
turns out to be a realization of $M$. The resulting quotient map induces a $\PP^1$-fibration 
$\Cal Z \rightarrow M$ different from the twistor fibration. 
This is a Riemannian submersion w.r.to any Z-metric on $\Cal Z$ 
and the original quaternion K\"ahler metric on $M$. However the fibers are not 
totally geodesic (see Remark 3.4). 
 }
\endproclaim

\proclaim{Remark 3.4}{\rm In the case of the canonical deformation metrics, the formula 
$d\alpha_2=2\alpha_3\wedge \alpha_1+2({}^tX^2\wedge X^0+{}^tX^3\wedge X^1)$ 
implies that the vertical distribution is integrable and the twistor lines are totally geodesic. 
In the case of the Z-metrics,  
the formula $d\alpha_2=2\alpha_3\wedge \alpha_1+2({}^tX^2\wedge X^0+{}^tX^3\wedge X^1)
=2(\hat\alpha_3+\alpha_3(\xi_i)X^i)\wedge )\hat\alpha_1+\alpha_1(\xi^j)X^j)
+({}^tX^2\wedge X^0+{}^tX^3\wedge X^1)=2\hat\alpha_3\wedge\hat\alpha_1
+\hat\alpha_3\wedge\alpha_1(\xi_j)X^j
+\alpha_3(\xi_i)X^i\wedge\hat\alpha_1+\alpha_3(\xi_i)X^i\wedge \alpha_1(\xi_j)X^j
+2({}^tX^2\wedge X^0+{}^tX^3\wedge X^1)$ implies that the distribution 
defined by the equations $X^i=0$ ($i=0,1,2,3$) is integrable. 
This formula involves non-trivial mixed terms like 
$\alpha_3(\xi_i)X^i\wedge\hat\alpha_1$ and so on. Therefore the distribution defined 
by the equations $X^i=0$ ($i=0,1,2,3$) is not totally geodesic.}
\endproclaim

\proclaim{Remark 3.5 (Comparison with the orbifold case)}{\rm 
We can construct locally irreducible positive quaternion K\"ahler orbifolds 
which are uniformized by one of the Wolf spaces. 
On the other hand, many examples of non locally symmetric positive quaternion 
K\"ahler orbifolds are constructed in [G-L]. 
Here we remark that the moving frame construction of Z-metrics in \S3 does not necessarily 
generalize to positive quaternion K\"ahler orbifold case.  
Here we explain the reason. 
If we take a local uniformization of the orbifold along the locus of orbifold singularities, 
we locally get a non-singular irreducible quaternion K\"ahler manifold with a finite group 
$G$ acting isometrically preserving the local quaternion K\"ahler structure and therefore 
$G$ operates on the local holonomy reduction $\Cal P_{\loc}$ of the oriented 
orthonormal frame bundle. We can construct the orbifold version of the twistor space. 
To see what happens to the construction of the orbifold version of the Z-metrics, 
we work on the local uniformization level. 
The essential step in the construction is to determine the complex line in each $\H$-linear 
subspace explained in the beginning of \S3. 
We recall this step. Let $J$ be the orthogonal complex structure of $\H$ represented by 
$z\in\PP^1_m$. Then this defines a $S^1$-subgroup in $\Sp(1)=\Sp(1)_r$ 
acting from the right on $T_mM$. Identifying $\Sp(1)_r$ with $\Sp(1)_l$ we get a $S^1$ 
subgroup of $\Sp(1)_l$. This $S^1$-subgroup determines a complex line $L_J$ in $\H$ 
(corresponding to the axis of the rotation of the action induced on $\PP^1=\PP(\C^2_J)$). 
Suppose that $z\in\Z$ is fixed by a non-trivial subgroup of $G$. 
Then the linear isotropy representation defines a finite subgroup of $\Sp(1)$ 
acting on the $\H$-line. The orbifold version of Z-metrics is well-defined if and only if 
the association $\PP^1_m \ni z \mapsto L_J \in \PP^1(\C^2_J)$ is preserved 
by the action of $G$. In other words, the complex line $L_J$ is not necessarily fixed 
by this action and if not fixed, the very beginning of the construction of  Z-metrics 
do not work. 
Therefore we cannot define Z-metric in 
the equivariant way and this implies that the orbifold Z-metric is not defined in 
general (the case where the orbifold version of the Z-metric is defined 
corresponds to orbifolds uniformized by the Wolf spaces). }
\endproclaim

\bigskip
\noindent
{\bf \S4. Ricci Flow on the Twistor Space of a Positive Quatermion K\"ahler Manifold.}
\medskip

In \S4, we study the behavior of the canonical deformation metrics $\Cal F^{\can}$ 
and the Z-metrics $\Cal F^{\ZZ}$ on the twistor space $\Cal Z$ of a positive quaternion 
K\"ahler manifold $(M,g)$. 
Let $(M^{4n},g)$ ($n\geq 2$) be a compact quaternion K\"ahler manifold with positive 
scalar curvature. Then $(M,g)$ is positive Einstein satisfying $\Ric(g)=(4n+8) g$. 
Then the solution of the 
Ricci flow equation $\p_t g=-2\Ric$ with the initial condition $g(-1)=(8n+16)g$ at time $t=-1$ 
is given by the homothety $g(t)=-(8n+16)tg$ ($-\infty<t<0$). Indeed, putting $g(t)=\lambda(t)g$, 
the above initial value problem of the Ricci flow equation becomes $\lambda'(t)=-(8n+16)$, 
$\lambda(-1)=8n+16$. Its solution is $\lambda(t)=-(8n+16)t$. 
Therefore the metric $g$ is just the fixed point modulo homothety of the Ricci flow and 
is not so interesting as itself. 
However, we have more freedom in the twistor space $\Cal Z$. 
Indeed, as in \S2 and \S 3, we can construct, from the original quarternion K\"ahler metric, 
two Einstein metrics $g^{\can}_1$ and $g^{\ZZ}_{\sqrt{\frac1{n+2}}}$ 
on $\Cal Z$ and two families of metrics 
$\Cal F^{\can}=\{g_{\lambda}^{\can}\}$ and $\Cal F^{\ZZ}=\{g_{\lambda}^{\ZZ}\}$ 
on $\Cal Z$ containing one of these Einstein metrics. 

\proclaim{Proposition 4.1} The homothetically extended family 
$\{\rho g_{\lambda}^{\can}\}_{\rho,\lambda>0}$ of the canonical deformation metrics 
on the twistor space $\Cal Z$ consists of orbits of the Ricci flow, i.e., the Ricci flow equation 
preserves the family $\{\rho\,g_{\lambda}^{\can}\}_{\rho,\lambda>0}$ on $\Cal Z$. 
\endproclaim
\noindent
{\it Proof}. The conclusion of Theorem 4.1 is a consequence from the following two facts. 
Firstly, the family $\{\rho g_{\lambda}^{\can}\}_{\rho,\lambda>0}$ is closed under the convex sum. 
Secondly, although the Ricci tensor $\Ric_{\lambda}^{\can}$ of $g_{\lambda}$ 
is not necessarily positive definite, it is of the same type as the canonical deformation 
metrics on $\Cal Z$. Indeed, we have from (2-15) in Proposition 2.4 the formula
$$
\split
\Ric_{\lambda}^{\can}&=4(1+n\lambda^4)(\alpha_1^2+\alpha_3^2)
+4(n+2-\lambda^2)\sum_{i=0}^3{}^tX^i\cdot X^i\\
&=4(n+2-\lambda^2)\,g^{\can}_{\sqrt{\frac{1+n\lambda^4}{n+2-\lambda^2}}}\,\,.
\endsplit
$$
Here we must assume that $\lambda^2<n+2$. 
Combining these two facts, we infer that the Ricci flow equation $\p_t g=-2\,\Ric_g$ 
preserves the family $\{\rho g^{\can}_{\lambda}\}_{\rho,\lambda>0}$ if $\lambda^2<n+2$. 
\qed
\bigskip

Proposition 4.2 implies that the Ricci flow with initial metric chosen from the family of 
canonical deformation metrics reduces to a system of ODE's : 

\proclaim{Proposition 4.2} The Ricci flow equation $\p_tg=-2\,\Ric_g$ on the twistor space 
$\Cal  Z$ with initial metric in the homothetically extended family of the canonical deformation 
metrics reduces to the system of ordinary differential equations
$$
\left\{\aligned
& \frac{d}{dt}(\rho(t)\lambda^2(t))=-8(1+n\lambda(t)^4)\,\,,\\
& \frac{d}{dt}\rho(t)=-8(n+2-\lambda(t)^2)\,\,.
\endaligned\right.\tag4-1
$$
\endproclaim

We examine the behavior of the solutions using the equation
$$
\rho\frac{d\lambda^2}{dt}
=-4\{(n+1)\lambda^2-1\}(\lambda^2-1)\,\,.\tag4-2
$$
A solution to the system of ODE's (4-1) corresponds to a curve in $(\lambda,\rho)$-plane 
(where $\lambda,\rho>0$). 
The curve $\lambda^2=1$ and $\lambda^2=\frac1{n+1}$ 
correspond to two Einstein metrics. The solution with initial metric $g_{\lambda}^{\can}$ 
$\lambda^2$ slightly larger than $1$ 
corresponds to a curve $(\lambda(t),\rho(t))$ with the property that both $\lambda(t)$ 
and $\rho(t)$ decrease as $t$ increases. 
The right hand side of (4-2) is strictly negative if $\lambda^2>1$. 
This observation implies that the solution converges to the K\"ahler-Einstein metric 
$g^{\can}_1$ (modulo scaling). Similarly, the solution with initial metric $g_{\lambda}^{\can}$ 
with $\lambda^2$ slightly smaller than $1$ corresponds to a curve $(\lambda(t),\rho(t))$ 
with the property that $\lambda(t)$ increases and $\rho(t)$ decreases as $t$ increases. 
Moreover (4-2) implies that $\rho\frac{d\lambda^2}{dt}>0$ which implies that $\rho>0$ 
if $\frac1{n+1}<\lambda^2<1$. This implies that the solution converges to 
the K\"ahler-Einstein metric $g^{\can}_1$ (modulo scaling). Next we look at the solution 
with initial metric $g_{\lambda}^{\can}$ $\lambda^2$ being slightly larger than $\frac1{n+1}$. 
In this case $\lambda$ increases as $t$ increases. Therefore the solution approaches 
to the Einstein metric $g^{\can}_{\sqrt{\frac1{n+1}}}$ (modulo scaling) as $t\to\infty$, 
i.e., the solution is an ancient solution. However, as $t$ increases, the solution 
becomes extinct in finite time and after scaling approaches to the K\"ahler-Einstein metric 
$g^{\can}_{\sqrt{\frac1{n+1}}}$. 
Similarly, the solution with initial metric $g_{\lambda}^{\can}$ $\lambda^2$ being slightly 
smaller than $\frac1{n+1}$ is also an ancient solution. To examine the behavior when 
$t$ increases, we need the defining equation of the trajectory. Eliminating $t$ from (4-1) 
we have $\log\rho-c=(1+\frac1n)\log|\lambda^2-1|-\frac1{n+1}(n+\frac1n+3)\log|(n+1)\lambda^2-1|$ 
for $\exists c\in\R$. 
It follows from this that the solution becomes totally singular in finite time 
and realizes the collapse corresponding to the limit $\lambda\to 0$ 
and $\rho\to e^c>0$, i.e., $\Cal Z$ collapses to $M$. 
Because of these properties when $t$ increases, the analysis in \S 5 
does not apply in this situation. 

We next look at the homothetically extended family of Z-metrics. 

\proclaim{Proposition 4.3} The homothetically extended family 
$\{\rho g_{\lambda}^{\ZZ}\}_{\rho,\lambda>0}$ of the Z-metrics 
on the twistor space $\Cal Z$ consists of orbits of the Ricci flow, i.e., the Ricci flow equation 
preserves the family $\{\rho\,g_{\lambda}^{\ZZ}\}_{\rho,\lambda>0}$ on $\Cal Z$. 
\endproclaim
{\it Proof}. The conclusion of Proposition 4.3 is a consequence from the following two facts. 
Firstly, the family $\{\rho g_{\lambda^{\ZZ}}\}_{\rho,\lambda>0}$ is closed under the 
convex sum. 
Secondly, the Ricci tensor $\Ric_{\lambda}^{\ZZ}$ of $g_{\lambda}$ 
is positive definite and is of the same type as Z-metrics s on $\Cal Z$. 
Indeed, we have from (3-14) in Proposition 3.1 the formula
$$
\Ric_{\lambda}^{\ZZ}=4(\alpha_1^2+\alpha_3^2)+(4n+8)\sum_{i=0}^3{}^tX^i\cdot X^i
=(4n+8)\,g^{\ZZ}_{\sqrt{\frac1{n+2}}}\,\,.
$$
Combining these two facts, we infer that the Ricci flow equation $\p_t g=-2\,\Ric_g$ 
preserves the family $\{\rho g^{\ZZ}_{\lambda}\}_{\rho,\lambda>0}$. 
\qed

\proclaim{Theorem 4.4} (1) The Ricci flow equation $\p_tg=-2\,\Ric_g$ on the twistor space 
$\Cal  Z$ with initial metric in the homothetically extended family of the 
Z-metrics reduces to the system of ordinary differential equations
$$
\left\{\aligned
& \frac{d}{dt}(\rho(t)\lambda^2(t))=-8\,\,,\\
& \frac{d}{dt}\rho(t)=-8(n+2)\,\,.
\endaligned\right.\tag4-3
$$

(2) For any initial metric with $\lambda^2=\frac1{n+2}$ at time $t=0$ 
in the homothetically extended family of Z-metrics on $\Cal Z$, the system 
of ODE's (4-2) has a solution
$$
\lambda^2\equiv\frac1{n+2}\,\,.
$$
This corresponds to the Einstein metric $g^{\ZZ}_{\sqrt{\frac1{n+2}}}$. 
\medskip

(3) For any initial metric $\rho_0g^{\ZZ}_{\lambda_0}$ with $\lambda_0^2\not=\frac1{n+2}$ at time $t=0$ 
in the homothetically extended family of Z-metrics on $\Cal Z$, the system 
of ODE's (4-3) has a solution
$$
\split
& \rho(t)=\rho_0-8(n+2)t\,\,,\\
& \rho\lambda^2(t)=\rho_0\lambda_0^2-8t\,\,.
\endsplit
$$
This implies that the Ricci flow solution $g(t)$ is given by
$$
g(t)=(\rho_0\lambda_0^2-8t)(\hat\alpha_1^2+\hat\alpha_3^2)+(\rho_0-8(n+2)t)\sum_{i=0}^3{}^tX^i\cdot X^i\,\,.$$

The solution is an ancient solution. Moreover its asymptotic soliton 
is the Einstein metric $g^{\ZZ}_{\lambda}$ with $\lambda^2=\frac1{n+2}$, i.e., 
modulo scaling, the solution $g(t)$ converges to the Einstein metric $g^{\ZZ}_{\sqrt{\frac1{n+2}}}$ 
as $t\to-\infty$. 

If $\lambda_0^2>\frac1{n+2}$ at $t=0$, then the solution becomes extinct in finite 
time and after scaling approaches to a Carnot-Carath\'eodory metric corresponding 
to $g^{\ZZ}_{\lambda}$ with $\lambda\to\infty$. 

If $\lambda_0^2<\frac1{n+2}$ at $t=0$, then the solution becomes totally singular 
in finite time and after scaling realizes the ``collapse'' of $\Cal Z$ corresponding to 
$\lambda\to 0$ and $\rho \to \rho_0\{1-(n+2)\lambda_0^2\}>0$. 
\endproclaim 
\noindent
{\it Proof}. The assertions (1), (2) and (3) follow from Proposition 4.3. 
On the other hand, the solution to (4-3) described in (3) implies
$$\lambda^2(t)=\frac{\rho_0\lambda_0^2-8t}{\rho_0-8(n+2)t}\tag4-4$$
from which the rest of the assertion (3) follows. 
\qed
\medskip

The remarkable difference from Proposition 4.3 is the following, which explains 
the meaning of Theorem 4.4 (3). 
A solution of (4-3) with $\lambda_0^2$ slightly larger than $\frac1{n+1}$ 
is an ancient solution having the Einstein metric $g^{\can}_{\sqrt{\frac1{n+1}}}$ 
as its asymptotic soliton when $t\to-\infty$ and it approximates the K\"ahler-Einstein 
metric $g^{\can}_1$ just before the extinction (after scaling). 
A solution of (4-3) with $\lambda^2$ slightly larger than $\frac1{n+2}$ is also 
an ancient solution having the Einstein metric $g^{\ZZ}_{\sqrt{\frac1{n+2}}}$ as its 
asymptotic soliton. However, the behavior of the solution just before the extinction 
is essentially different. Namely we can scale the solution so that the $\Cal D$-direction 
survives in the limit toward the extinction time. We will exploit this difference 
in \S5. 
\medskip

Eliminating $t$ from (4-3) we have
$$\rho=\frac{\rho_0\biggl(\lambda_0^2-\frac1{n+2}\biggr)}{\lambda^2-\frac1{n+2}}\,\,.\tag4-5$$
This is the equation of the trajectory of the Ricci flow solution with initial metric 
corresponding to $(\lambda_0,\rho_0)$. 
\medskip

The 2-dimensional family $\Cal F=\{\rho g_{\lambda}\}_{\rho>0,\lambda>0}$ 
constitutes a Ricci flow unstable cell in the sense that the family $\Cal F$ is foliated by 
the trajectories of the Ricci flow solutions and each Ricci flow trajectory is an ancient solution 
whose asymptotic soliton (in the sense of [P,\S11]) corresponds to a K\"ahler-Einstein metric. 

\proclaim{Example 4.5}{\rm Pick a trajectory defined by the equation 
$$\rho=\frac{\rho_0\biggl(\lambda_0^2-\frac1{n+2}\biggr)}{\lambda^2-\frac1{n+2}}$$ 
where $\rho_0>0$ and $\lambda_0^2>\frac1{n+2}$ in the $(\lambda,\rho)$-plane 
identified with the family $\R_+\cdot \Cal F^{\ZZ}$. 
Along this trajectory we have $\Ric_{\lambda}^{\ZZ}=4(\alpha_1^2+\alpha_3^2)
+(4n+8)\sum_{i=0}^3{}^tX^i\cdot X^i=(4n+8)\,g^{\ZZ}_{\sqrt{\frac1{n+2}}}$. 
Therefore the scalar curvature along the trajectory is
$$\text{\rm Scal}(\rho g_{\lambda}^{\ZZ})=\frac{8}{\rho\lambda^2}+\frac{16(n+2)n}{\rho}\,\,.$$
If we set $u=\text{\rm constant}$ determined by $\int_MudV=1$, i.e., $u=1/\Vol(g_{ij}(t))$, 
$g_{ij}(t)$ being the solution, we get a solution $u(t,x)$ ($t$-dependent constant function on $M$) 
to the conjugate heat equation 
$\p_t u=-\triangle u+Ru$. 
Since $\displaystyle \Vol(\rho g^{\ZZ}_{\lambda})
=\frac{\rho_0^{2n+1}\biggl(\lambda_0^2-\frac1{n+2}\biggr)^{2n+1}}
{\biggl(\lambda^2-\frac1{n+2}\biggr)^{2n+1}}\lambda^2\Vol(M,g)$, we have 
$$u=\frac{\biggl(\lambda^2-\frac1{n+2}\biggr)^{2n+1}}
{\rho_0^{2n+1}\biggl(\lambda_0^2-\frac1{n+2}\biggr)^{2n+1}}\,\,.$$
Eliminating $\rho$ from (4-3) we have
$$\tau:=-t=\frac1{8(n+2)}\biggl\{\frac{\rho_0\biggl(\lambda_0^2-\frac1{n+2}\biggr)}
{\lambda^2-\frac1{n+2}}-\rho_0\biggr\}\,\,.$$
We observe that 
$$\lambda^2\to\frac1{n+2}+0 \Longleftrightarrow \tau\to\infty \Longleftrightarrow t\to-\infty\,\,.$$ 
The function $W(g_{ij},f,\tau)$ ($W$ being Perelman's $W$-functional) 
is monotone increasing along the Ricci flow trajectory passing through a metric 
$\rho_0 g_{\lambda_0}^{\ZZ}$ with $\lambda_0^2>\frac1{n+2}$, 
which is determined by the triple 
$(\rho g_{\lambda},f,\tau)$ where $\rho$, $\lambda$, $\tau$ are given as above, 
$\lambda^2 \in (\frac1{n+2},\infty)$ increases to from $\frac1{n+2}$ to $\infty$ 
when $\tau$ decreases from $\infty$ to $0$), 
and $f$ is determined by setting $u=(4\pi \tau)^{-(2n+1)}e^{-f}$ with $u$ and $\tau$ 
given as above.} 
\endproclaim

\bigskip
\noindent
{\bf \S5. Uniformization of Positive Quaternion K\"ahler Manifolds.} 
\medskip

Let $M^{4n}$ be a compact quaternion K\"ahler manifold with positive scalar curvature 
and $\Cal Z$ its twistor space. Let 
$\R_+\cdot\Cal F^{\ZZ}=\{\rho g^{\ZZ}_{\lambda}\}_{\rho,\lambda>0}$ 
be the homothetically extended family of Z-metrics on the twistor space $\Cal Z$. 
Consider the ``Ricci map" from the space of Riemannian metrics on $\Cal Z$ to the 
space of symmetric $(0,2)$-tensors on $\Cal Z$ defined by $g \mapsto \Ric_g$. 
Then we have shown in Section 3 and 4 (see Theorem 4.4) that 
\medskip

1) the space of homothetically extended family of Z-metrics $\R_+\cdot\Cal F^{\ZZ}$ 
contains the half-line consisting the scalings of the Einstein metric 
$g^{\ZZ}_{\sqrt{\frac1{n+2}}}$ on $\Cal Z$. 
\medskip

2) The space of homothetically extended family of Z-metrics 
$\R_+\cdot\Cal F^{\ZZ}$ 
on $\Cal Z$ consists of orbits of the Ricci flow, more precisely, the family 
$\R_+\cdot\Cal F^{\ZZ}$ is invariant under the "Ricci map" and the Ricci flow. 
\medskip

It is natural to regard the existence of a special family of ancient solutions whose 
asymptotic soliton is a homothetical family of a fixed Einstein metric as an extension 
of the notion of a single Einstein metric. 
The following Theorem 5.1 is an evidence for the usefulness of such an extended notion. 
Indeed, we obtain strong information on the Einstein metric from the analysis 
of the ancient solutions. 
Our strategy is to apply Bando-Shi's gradient estimate ([B], [Sh1,2]) for the Ricci flow 
to the ancient solution in Theorem 4.4 (3) for $\lambda_0^2>\frac1{n+2}$. 
To do so, we need to know the behavior of the full curvature tensor 
along the solution. Using the expression (3-13) of the Levi-Civita connection 
$\Gamma^{\ZZ}_{\lambda}$ of $g^{\ZZ}_{\lambda}$ and the second structure equation, 
we can estimate the norm of the curvature tensor 
of $g_{\ZZ}^{\lambda}$. 
The second structure equation says that the curvature form is defined by 
$\Om^{\ZZ}_{\lambda}:=d\Gamma^{\ZZ}_{\lambda}+\Gamma^{\ZZ}_{\lambda}\wedge 
\Gamma^{\ZZ}_{\lambda}$. As $a_i,\dots,d_i$'s ($i=1,3$) are all ``invisible'' 
at $z\in\Cal Z$, we have only to examine $d\Gamma^{\ZZ}_{\lambda}$. 
We now estimate the norm of the curvature tensor by directly estimating all entries 
of the $d\Gamma^{\ZZ}_{\lambda}$-part of the curvature form. 
The reason we must do so is the following. Although estimating all sectional curvatures 
is equivalent to estimating the norm of the curvature tensor, estimating all sectional 
curvatures of the form $K(e_i,e_j)$ where $\{e_i\}$ is the orthonormal basis under 
consideration is not enough to estimate the norm of the curvature tensor. 
Now we look at (3-13). The trouble would be that the uncontrollable quantities 
of size $\lambda$ or $\lambda^2$ appear in the curvature form as $\lambda$ 
becomes large. So we estimate the norm of the curvature form modulo $O(1)$ 
when $\lambda$ becomes large. 
The only trouble from this view point which may occur when $\lambda$ 
becomes large stems from the exterior differential of the term like 
$\lambda^2(-a_1\alpha_1-a_3\alpha_3
+a_1\alpha_1(\xi_1)X^1+a_1\alpha_1(\xi_3)X^3+a_3\alpha_3(\xi_1)X^1
+a_3\alpha_3(\xi_3)X^3)$ and so on. The definition of $\hat\alpha_i$ ($i=1,3$) 
implies that 
the above quantity is equal to $\lambda^2(-a_1\hat\alpha_1-a_3\hat\alpha_3
-a_1\alpha_1(\xi_0)X^0-a_1\alpha_1(\xi_2)X^2-a_3\alpha_3(\xi_0)X^0-a_3(\xi_2)X^2)$. 
Here the terms like $a_1\alpha_1(\xi_0)$ are a product of  terms ``invisible'' at $z$ 
and so on. Therefore, we have only to estimate the exterior differential of the quantity 
$a_1\hat\alpha_1$ and $a_3\hat\alpha_3$ to estimate the norm of the exterior differential 
of the original quantity.  As $a_1$ and $a_3$ are ``invisible'' at $z$, or, 
as $d\hat\alpha_i$ ($i=1,3$) are ``invisible'' at $z$, we have only to 
consider $da_i\wedge \alpha_i$ ($i=1,3$). The condition (3-2$'$) implies 
that the exterior differential $da_1$ consists of $X^i$ modulo $\alpha_1$ and $\alpha_3$. 
Therefore the exterior differential of $a_1\hat\alpha_1$ and $a_3\hat\alpha_3$ consists 
of $\hat\alpha_1\wedge\hat\alpha_3$ and $X^i\wedge\hat\alpha_j$ ($i,j=1,3$). 
Note that the norm w.r.to the metric $g^{\ZZ}_{\lambda}$ 
of $\lambda^2\hat\alpha_1\wedge\hat\alpha_ 3$ 
(resp. $\lambda^2X^i\wedge\hat\alpha_j$) is $1$ (resp. $\lambda$). 
The same is true for those quantity like 
$\lambda^2(-b_1\alpha_1-b_3\alpha_3
+b_1\alpha_1(\xi_1)X^1+b_1\alpha_1(\xi_3)X^3+b_3\alpha_3(\xi_1)X^1
+b_3\alpha_3(\xi_3)X^3)$ and so on. By 
gathering all these terms with norm $\lambda$ arising this way, 
we can specify a special part of the curvature tensor. The sum of the norms of this 
special part is  $\lambda$ times some constant (depending only on $n$). 

\proclaim{Theorem 5.1} Let $(M^{4n},g)$ be a compact quaternion K\"ahler manifold 
with positive scalar curvature. Then the family $\{g^{\ZZ}_{\lambda}\}_{\lambda>1}$ of 
Z-metrics on the twistor space of $M$ satisfies the limit formula
$$\lim_{\lambda\to\infty}|\nabla^{g^{\ZZ}_{\lambda}}\Rm^{g^{\ZZ}_{\lambda}}
|_{g^{\ZZ}_{\lambda}}=0\,\,.$$
 \endproclaim
\noindent
{\it Proof.} We follow the proof of Bando-Shi's derivative estimate for the curvature tensor 
under the Ricci flow ([B], [Sh1,2], see also [C-K, Chapt.7]). We consider a sequence 
$\{\lambda_k,\rho_k\}_{k=1}^{\infty}$ 
of pairs of positive numbers satisfying $\forall \rho_k=1$, $1<\forall \lambda_k$ and 
$\lambda_k \to \infty$ as $k\to\infty$. 
This represents a sequence $\{g^{\ZZ}_{\lambda_k}\}_{k=1}^{\infty}$ of the Z-metrics 
on the twistor space $\Cal Z$. Pick one $g^{\ZZ}_{\lambda_k}$ and write it as $g_k$. 
We consider the trajectory of the Ricci flow passing through $g_k$ (with the maximal 
time interval in the past and future) and we write $g_k(t)$ for the Ricci flow solution 
where the initial metric $g_k(0)$ is taken from the trajectory 
as $\rho_kg^{\ZZ}_{\lambda(0)}$ ($\rho_k$ being large) such that $\lambda(0)$ 
is slightly larger than $1$, say, $\lambda^2(0)=\frac1{n+2}+\delta$, 
and $T_k$ is the time when the Ricci flow trajectory passes through the given 
$g_{\lambda_k}^{\ZZ}:=g_k$ : $g_k(T_k)=g_k$. 
Let $\nabla$ and $\Rm$ the Levi-Civita connection and the Riemann curvature tensor of 
the Ricci flow solution $g_k(t)$. 
Let $A*B$ denote any quadratic quantity obtained from $A\otimes B$ whose meaning 
is explained in [C-K, p. 227]. In the following computation, the norm at time $t$ should 
be understood to be computed with respect to the metric $g_k(t)$. We have
$$
\frac{\p}{\p t}|\Rm|^2=\triangle|\Rm|^2-2|\nabla\Rm|^2+(\Rm)^{*3}\tag5-1
$$
and
$$
\frac{\p}{\p t}|\nabla\Rm|^2
=\triangle|\nabla\Rm|^2-2|\nabla^2\Rm|^2+\Rm*(\nabla\Rm)^{*2}\,\,.\tag5-2
$$
The strategy of Bando-Shi derivative estimate is to make the best use of the good term 
$-2|\nabla\Rm|^2$ in (5-1) to kill the bad term $\Rm*(\nabla\Rm)^{*2}$ in 
(5.2). We also use (5-1) and (5-2) but in somewhat different way because in our case 
we can use the special properties of the curvature tensor of $g^{\ZZ}_{\lambda}$. 
We have specified the special part of the curvature tensor of $g^{\ZZ}_{\lambda}$. 
This was characterized by $R_{\eta\mu\nu\kappa}$ 
where only one index is negative (i.e., $-1$ or $-2$) and satisfies 
$|R_{\eta\mu\nu\kappa}|=\lambda$. 
Here, we are considering components of the curvature tensor 
w.r.to the orthonormal coframe $(\eta\hat\alpha_1,\lambda\hat\alpha_3,X^0,X^1,X^2,X^3)$. 
Let $\wt{\Rm}$ denote the part of the curvature tensor obtained by subtracting the special 
part specified above, i.e., if $R_{\eta\mu\nu\kappa}$ is a component contained in 
$\wt{\Rm}$ then the indices contain no negative number or more than two indices are negative 
(in fact the components with at least three indices are negative are $0$). 
The advantage of introducing $\wt\Rm$ is the following. 
If we replace $\Rm$ by $\wt\Rm$ in (5-1), the first two terms in the RHS of (5-1) do not 
change, i.e., $\triangle|\Rm|^2=\triangle|\wt\Rm|^2$ and $|\nabla\Rm|^2=|\nabla\wt\Rm|^2$ hold. 
Indeed, $|\Rm|^2$ and $|\wt\Rm|^2$ differ only by a constant and $\triangle|\Rm|^2=\triangle|\wt\Rm|^2$ 
follows. On the other hand, the portion of $\nabla\Rm$ 
stemming from the above specified part of the curvature tensor vanish 
and therefore $\nabla\Rm=\nabla\wt\Rm$ follows. 
The third term $(\Rm)^{*3}$ contains all terms in the direct computation of 
$\p/\p t|\wt\Rm|^2$ obtained by differentiating the metric tensor w.r.to $t$. 
It therefore follows from the Ricci flow equation and Proposition 4.3 
that the term corresponding to the third term in the RHS of (5-1) in the computation of 
$\p/\p t|\wt\Rm|^2$ is not larger than $c_n|\wt\Rm|^2$ where $c_n$ is a positive 
constant depending only on $n$. Therefore we have
$$
\frac{\p}{\p t}(e^{-c_nt}|\wt\Rm|^2)\leq \triangle (e^{-c_nt}|\wt\Rm|^2) -2e^{-c_nt}|\nabla\wt\Rm|^2\,\,.
\tag5-1$'$
$$
Moreover, if we replace $\Rm$ by $\wt\Rm$ in (5-2), the first two terms in the RHS of (5-1) do not 
change because $\nabla\Rm=\nabla\wt\Rm$. 
The third term $\Rm*(\nabla\Rm)^{*2}$ contains all terms in the direct computation of 
$\p/\p t|\nabla\wt\Rm|^2$ obtained by differentiating the metric tensor w.r.to $t$. 
Therefore, It follows from the same reason as above that 
the term corresponding to the third term in the RHS of (5-2) in the computation of 
$\p/\p t|\nabla\wt\Rm|^2$ is not larger than $c_n|\nabla\wt\Rm|^2$ where $c_n$ is the same 
positive constant as above. 
Therefore we have
$$
\frac{\p}{\p t}(e^{-c_nt}|\nabla\wt\Rm|^2)\leq \triangle (e^{-c_nt}|\nabla\wt\Rm|^2) 
-2e^{-c_nt}|\nabla^2\wt\Rm|^2\,\,.
\tag5-2$'$
$$
To prove Theorem 5.1, suppose the contrary, i.e., the maximum over $\Cal Z$ 
of $|\nabla\wt\Rm|^2=|\nabla\Rm|^2$ at time $T_k$ is uniformly (w.r.to $k$) bounded below 
by a positive constant $D$. Note that $|\wt\Rm|^2$ is uniformly bounded below by 
a positive constant for all $T_k$. This is because $g$ is a 
positive Einstein metric on $M$ and it follows from the second structure equation and (3-13) 
that  the curvature tensor of $M$ directly contributes to the norm of $\Om^{\ZZ}_{\lambda}$. 
Thus we can conclude that there exists a positive constant $C$ such that
$$
2\,\ep\,|\nabla\wt\Rm|^2\, \geq\, \ep\,C\,(|\nabla\wt\Rm|^2+\ep\,|\wt\Rm|^2) \tag5-3
$$
holds at the point where $|\nabla\Rm|$ takes its maximum. If $\ep<<1$ then $C\approx 2$. 
Take a small positive constant $\ep$ to be determined later. 
From (5-1$'$), (5-2$'$) and (5-3) we have
$$
\split
\frac{\p}{\p t}\{e^{-c_nt}(|\nabla\wt\Rm|^2+\ep\,|\wt\Rm|^2)\} \,
&  \leq -C\,e^{-c_nt}\,(|\nabla\wt\Rm|^2+\ep\,|\wt\Rm|^2)
-2\,e^{-c_nt}\,|\nabla^2\wt\Rm|^2\\
& +(\text{\rm error term depending on $\lambda$ and $\ep$})
\endsplit
$$
at the point where $|\nabla\wt\Rm|$ takes its maximum. Here, 
the error term stems from the quantity which bound the Laplacian term 
$\triangle \{(e^{-c_nt}(|\nabla\wt\Rm|^2+\ep\,|\wt\Rm|^2)\}$ from above. 
We note that $\triangle|\Rm|^2$ involves the curvature tensor 
of the metric $g_k$, we should take $\ep$ small compared to $\lambda=\lambda_k$ 
so that the Laplacian term $\triangle \{(e^{-c_nt}(|\nabla\wt\Rm|^2+\ep\,|\wt\Rm|^2)\}$ 
is uniformly of order, say, $\sqrt{\ep}$. For this purpose we should take $\ep$ 
comparable to or smaller than $1/\lambda_k^2$. 
So, there exists a positive constant $K$, which is taken to be uniform 
for all $\ep>0$ sufficiently small in the sense that $\ep \approx O(1/\lambda_k^2)$ 
holds, which satisfies the estimate
$$
\frac{\p}{\p t}\{e^{-c_nt}(|\nabla\wt\Rm|^2+\ep\,|\wt\Rm|^2)-K\sqrt{\ep}\,t\} \,
 \leq -\ep\,C\,e^{-c_nt}\,(|\nabla\wt\Rm|^2+\ep\,|\wt\Rm|^2-K\sqrt{\ep}\,t)
-2\,e^{-c_nt}\,|\nabla^2\wt\Rm|^2
\tag5-4
$$
at the point where $|\nabla\wt\Rm|^2$ takes its maximum. 
The parabolic maximum principle implies that we have
$$
\split
&\quad \{e^{-c_nt}(|\nabla\wt\Rm|^2+\ep\,|\wt\Rm|^2-K\sqrt{\ep}\,t)\}|_{t=T} \leq 
\{e^{-c_nt}(|\nabla\wt\Rm|^2+\ep\,|\wt\Rm|^2-K\sqrt{\ep}\,t)\}|_{t=0}\exp(-\ep\,CT)\\
& \leq \frac{\text{\rm const.}}{\rho_0^3}\exp(-\ep\,CT)
\endsplit
\tag5-5
$$
for $0\leq \forall T \leq T_k$.  From (5-5) we have
$$
|\nabla\wt\Rm|^2(T)+\ep\,|\wt\Rm|^2(T)
-K\sqrt{\ep}\,T\, \leq \, \frac{\text{\rm const.}}{\rho_0^3}\exp\{(c_n-\ep\,C)T\}\,\,.\tag5-6
$$
Now we take $\ep>0$ so that 
$$K\sqrt{\ep}T\leq \frac12 D$$
holds 
($D$ being the assumed uniform lower bound 
of $\max_{z\in \Cal Z}|\nabla\Rm|^2$ when $k\to\infty$). 
Once we choose $\ep>0$ this way, (5-6) implies that, in the above discussion, 
we can replace $c_n$ by $c_n-\ep\,C$ and repeat the same argument. 
Thus we can repeat the same argument with $c_n$ replaced by $c_n-\ep\,C, c_n-2\ep\,C, \dots$. 
There exists a positive integer $N$ such that $c_n-N\ep\,C<0$. Therefore we have
$$
|\nabla\wt\Rm|^2(T)+\ep\,|\wt\Rm|^2(T)-K\sqrt{\ep}\,T\, \leq \, 
\frac{\text{\rm const.}}{\rho_0^3}\exp\{(\underbrace{c_n-N\ep\,C}_{<0})T\}\,\,.\tag5-7
$$
For instance, we choose $\displaystyle \ep \leq \frac{D^2}{4K^2T_k^2}$ and 
$\displaystyle N>\frac{8c_nK^2T^2}{D^2}$. Then from (5-7) we have
$$
|\nabla\wt\Rm|^2(T_k)+O(T_k^{-2}) \leq \frac{\text{\rm const.}}{\rho_0^3}\exp(-T_k)\,\,.\tag5-8
$$
As we have assumed that $\max_{z\in\Cal Z}|\nabla\Rm|^2(T_k)$ is uniformly (w.r.to $k$) 
bounded below by a positive constant $D$, (5-8) implies
$$
D\leq |\nabla\wt\Rm|^2(T_k)+O(T_k^{-2}) \leq \frac{\text{\rm const.}}{\rho_0^3}\exp(-T_k)$$
as $k\to\infty$ ($T_k\to\infty$), which is clearly a contradiction. 
Therefore
$$
\lim_{k\to\infty}\max_{z\in\Cal Z}|\nabla^{g^{\ZZ}_{\lambda_k}}
\Rm(g^{\ZZ}_{\lambda_k})|_{g^{\ZZ}_{\lambda}}(z,T_k)=0
$$
must be the case and the assertion of Theorem 5.1 follows. \qed

\proclaim{Remark 5.2}{\rm (1) Although the proof of Theorem 5.1 is based 
on the local  computation in \S3 and \S4, Theorem 5.1 is a global result. 
Indeed, the parabolic maximum principle, which does not make sense (and not correct) 
in local situation, played an essential role in its proof. 
\medskip

(2) The 2-parameter family $\{\rho g^{\ZZ}_{\lambda}\}_{\rho,\lambda>0}$ 
of homothetically extended Z-metrics on the twistor space $\Cal Z$ 
is foliated by the trajectories of the Ricci flow solution. We have two facts (i) each Ricci flow solution 
in this family is an ancient solution (in the sense of [H,\S19]) and (ii) it realizes the ``collapse'' 
where the base ($\{\xi_i\}_{i=0}^3$-) direction shrinks faster when $t$ becomes large. 
These two facts are essential in the proof of Theorem 5.1. Indeed, because each trajectory 
corresponds to an ancient solution, for any sequence $\lambda_k$ such that 
$\lim_{k\to\infty}\lambda_k=\infty$, 
we can find a Ricci flow solution defined on $[0,T_k]$ with initial metric 
with $\lambda^2=1+\delta$ and $g(T_k)=g_k$, where $g_k=g_{\lambda_k}$. 
This argument does not work for the ancient solution in 
Proposition 4.2 having the Einstein 
metric $g^{\can}_{\sqrt{\frac1{n+1}}}$ as its asymptotic soliton. }
\endproclaim

\proclaim{Theorem 5.3} Any locally irreducible positive quaternion K\"ahler 
manifold $(M^{4n},g)$ is isometric to one of the Wolf spaces. \endproclaim
\noindent
{\it Proof}. The goal is to prove 
$\nabla^g\Rm(g)\equiv 0$, i.e., . $(M,g)$ is a Riemannian symmetric space. 
In order to do so we compute the covariant derivative 
$\nabla^{g^{\ZZ}_{\lambda}} \Rm(g^{\ZZ}_{\lambda})$ 
(we will write this simply as $\nabla\Rm^{\ZZ}_{\lambda}$ for brevity) 
using (3-13) and compare the result with the limit formula in Theorem 5.1 
(note that the curvature form of $g$ is contained in the curvature form 
of $g^{\ZZ}_{\lambda}$). 
The covariant derivative $\nabla\Rm^{\ZZ}_{\lambda}$ is computed from 
$d\Om^{\ZZ}_{\lambda}$ and the product of $\Gamma^{\ZZ}_{\lambda}$ 
and $\Om^{\ZZ}_{\lambda}$. From the second structure equation 
(or the second Bianchi identity) $d\Om^{\ZZ}_{\lambda}$ is equal to 
$$
(*):=d\Gamma^{\ZZ}_{\lambda}\wedge\Gamma^{\ZZ}_{\lambda}
-\Gamma^{\ZZ}_{\lambda}\wedge d\Gamma^{\ZZ}_{\lambda}\,\,.
$$
On the other hand $\Gamma^{\ZZ}_{\lambda}$ is ``invisible'' at the point 
$z\in\Cal Z$ where the computation in \S3 is performed, 
except for the Levi-Civita connection form of $g$. 
We use the matrix (3-13) of $\Gamma^{\ZZ}_{\lambda}$ to 
list the parts of $d\Gamma^{\ZZ}_{\lambda}$ which do not come from 
the Levi-Civita connection form of $g$ : 
\newline\noindent
(i) $d\alpha_2=2\alpha_3\wedge\alpha_1+2(X^2\wedge X^0+X^3\wedge X^1)$. 
In the product ($*$) this is coupled with the invisible term (i.e., $\alpha_2$) 
and therefore it does not survive in ($*$). 
\newline\noindent
(ii) the 2-forms obtained by taking the exterior differential of 
$\lambda(a_1X^1+b_1X^3)$ and so on in the first two raws of (3-13). 
These consist of the linear combination of $\lambda \alpha_i\wedge X^j$ and 
$\lambda X^i\wedge X^j$ ($i,j=1,3$) at $z$. 
\newline\noindent
(iii) the 2-forms obtained by taking the exterior differential of 
$\lambda^2(-a_1\alpha_1-a_3\alpha_3
+a_1\alpha_1(\xi_1)X^1+a_1\alpha_1(\xi_3)X^3+a_3\alpha_3(\xi_1)X^1
+a_3\alpha_3(\xi_3)X^3)=\lambda^2(-a_1\hat\alpha_1-a_3\hat\alpha_3
-a_1\alpha_1(\xi_0)X^0-a_1\alpha_1(\xi_2)X^2-a_3\alpha_3(\xi_0)X^0-a_3(\xi_2)X^2)$ 
and so on in (3-13). 
These consist of the linear combination of $\lambda^2X^i\wedge \hat\alpha_j$ 
and $\lambda^2\alpha_i\wedge\hat\alpha_j$ ($i,j=1,3$) at $z$. 

In the computation of $\nabla\Rm^{\ZZ}_{\lambda}$, the 3-forms appearing 
from different entries of the matrix representation of ($*$) are coupled with 
the tensor product of the members of the orthonormal coframe 
$(\hat\alpha_1,\hat\alpha_3,X^0,X^1,X^2,X^3)$ which are mutually orthogonal. 
Therefore, (ii) and (iii) from the above listed 2-forms contribute with coefficients 
$\lambda$ and $\lambda^2$ (if they survive after computing ($*$)). 
Therefore, ignoring the ``invisible'' terms, we have the expansion
$$
|\nabla\Rm^{\ZZ}_{\lambda}|^2=A_0+\lambda A_1+ \lambda^2 A_2\,\,.
$$
The limit formula in Theorem 5.1 implies $A_1\equiv A_2\equiv 0$. 
The $A_0$-part consists only of the curvature form of the Levi-Civita connection 
of $g$. Therefore $A_0=|\nabla^g\Rm(g)|_g\equiv 0$ which 
completes the proof of Theorem 5.3.  \qed

\enddocument